\documentclass[11pt]{amsart}
\usepackage{amsbsy,amssymb,amscd,amsfonts,latexsym,amstext,delarray,
amsmath,amsthm, graphicx}

\input xypic

\usepackage{hyperref}

\newtheorem{thm}{Theorem}[section]
\newtheorem{prop}[thm]{Proposition}
\newtheorem{cor}[thm]{Corollary}
\newtheorem{lem}[thm]{Lemma}

\newtheorem{defn}[thm]{Definition}
\newtheorem{rem}[thm]{Remark}

\numberwithin{equation}{section}

\def\bC{{\mathbb C}}

\def\bL{{\mathbb L}}

\def\bR{{\mathbb R}}

\def\bZ{{\mathbb Z}}

\def\A{{\mathbb A}}
\def\C{{\mathbb C}}

\def\F{{\mathbb F}}

\def\N{{\mathbb N}}

\def\Q{{\mathbb Q}}
\def\R{{\mathbb R}}
\def\Z{{\mathbb Z}}
\def\K{{\mathbb K}}

\def\cA{{\mathcal A}}
\def\cB{{\mathcal B}}
\def\cC{{\mathcal C}}

\def\cG{{\mathcal G}}
\def\cH{{\mathcal H}}

\def\cJ{{\mathcal J}}
\def\cK{{\mathcal K}}
\def\cL{{\mathcal L}}

\def\cO{{\mathcal O}}
\def\cP{{\mathcal P}}

\def\cS{{\mathcal S}}

\def\cU{{\mathcal U}}
\def\cV{{\mathcal V}}

\def\Fr{{\rm Fr}}
\def\Gal{{\rm Gal}}
\def\GL{{\rm GL}}

\def\Ker{{\rm Ker}}

\def\Out{{\rm Out}}

\def\Sp{{\rm Spec}}

\def\Tor{{\rm Tor}}

\def\Tr{{\rm Tr}}

\def\n{{\mathfrak n}}

\def\fE{{\mathfrak E}}
\def\urep{\vartheta_a}
\def\vrep{\vartheta_m}

\newcommand{\ie}{{\it i.e.\/}\ }
\newcommand{\eg}{{\it e.g.\/}\ }
\newcommand{\cf}{{\it cf.\/}\ }

\def\qqq{\,,\quad \forall}
\title[Weil's proof and adeles classes]
{The Weil proof and the geometry of the adeles class space}

\author[Connes]{Alain Connes}
\author[Consani]{Caterina Consani}
\author[Marcolli]{Matilde Marcolli}
\address{A.~Connes: Coll\`ege de France \\
3, rue d'Ulm \\ Paris, F-75005 France}
\email{alain\@@connes.org}
\address{C.~Consani: Mathematics Department \\ Johns Hopkins
University \\ Baltimore, MD 21218 USA}
\email{kc\@@math.jhu.edu}
\address{M.~Marcolli: Max--Planck Institut f\"ur Mathematik  \\
Vivatsgasse 7 \\ Bonn, D-53111 Germany}
\email{marcolli\@@mpim-bonn.mpg.de}

\begin{document}

\maketitle

\begin{center}
{\em  Dedicated to Yuri Manin on the occasion of his 70th birthday}
\end{center}

\medskip
\begin{verse}
O simili o dissimili che sieno questi mondi \\
non con minor raggione sarebe bene a l'uno  \\
l'essere che a l'altro
\\
\smallskip
Giordano Bruno -- {\em De l'infinito, universo e mondi}
\end{verse}

\tableofcontents

\section{Introduction}

This paper explores analogies between the Weil proof of the Riemann
Hypothesis for function fields and the geometry of the adeles class
space, which is the noncommutative space underlying the spectral
realization of the zeros of the Riemann zeta function constructed in
\cite{Co-zeta}. Our purpose is to build a dictionary between the
algebro-geometric setting of algebraic curves, divisors, the
Riemann--Roch formula, and the Frobenius map, around which the Weil
proof is built, and the world of noncommutative spaces, cyclic
cohomology and KK-theory, index formulae, and the thermodynamical
notions of quantum statistical mechanics, which, as we already
argued in \cite{CCM}, provide an analog of the Frobenius in
characteristic zero via the scaling action on the dual system.

The present work builds upon several previous results. The first
input is the spectral realization of \cite{Co-zeta}, where the
adeles class space was first identified as the natural geometric
space underlying the Riemann zeta function, where the Weil explicit
formula acquires an interpretation as a trace formula. In
\cite{Co-zeta} the analytic setting is that of Hilbert spaces, which
provide the required positivity, but the spectral realization only
involves the critical zeros. In \cite{CCM}, we provided a
cohomological interpretation of the trace formula, using cyclic
homology. In the setting of \cite{CCM}, the analysis is as developed
by Ralph Meyer in \cite{Meyer} and uses spaces of rapidly decaying
functions instead of Hilbert spaces. In this case, all zeros
contribute to the trace formula, and the Riemann Hypothesis becomes
equivalent to a positivity question. This mirrors more closely the
structure of the two main steps in the Weil proof, namely the
explicit formula and the positivity $\Tr(Z*Z')>0$ for
correspondences (see below). The second main building block we need
to use is the theory of endomotives and their quantum statistical
mechanical properties we studied in \cite{CCM}. Endomotives are a
pseudo-abelian category of noncommutative spaces that naturally
generalize the category of Artin motives. They are built from
semigroup actions on projective limits of Artin motives. The
morphisms in the category of endomotives generalize the notion of
correspondence given by algebraic cycles in the product used in the
theory of motives to the setting of \'etale groupoids, to account
naturally for the presence of the semigroup actions. Endomotives
carry a Galois action inherited from Artin motives and they have
both an algebraic and an analytic manifestation. The latter provides
the data for a quantum statistical mechanical system, via the
natural time evolution associated by Tomita's theory to a
probability measure carried by the analytic endomotive. The main
example that is of relevance to the Riemann zeta function is the
endomotive underlying the Bost--Connes quantum statistical
mechanical system of \cite{BC}. One can pass from a quantum
statistical mechanical system to the ``dual system'' (in the sense
of the duality of type III and type II factors in \cite{Co-th},
\cite{Tak}), which comes endowed with a scaling action induced by
the time evolution. A general procedure described in \cite{CCM}
shows that there is a ``restriction map'' (defined as a morphism in
the abelian category of modules over the cyclic category) from the
dual system to a line bundle over the space of low temperature KMS
states of the quantum statistical mechanical system. The cokernel of
this map is not defined at the level of algebras, but it makes sense
in the abelian category and carries a corresponding scaling action.
We argued in \cite{CCM} that the induced scaling action on the
cyclic homology of this cokernel may be thought of as an analog of
the action of Frobenius on \'etale cohomology. This claim is
justified by the role that this scaling action of $\R^*_+$, combined
with the action of $\hat\Z^*$ carried by the Bost--Connes
endomotive, has in the trace formula, see \cite{Co-zeta}, \cite{CCM}
and \S 4 of \cite{CMbook}. Further evidence for the role of the
scaling action as Frobenius is given in \cite{ConsMar-ff}, where it
is shown that, in the case of function fields, for a natural quantum
statistical mechanical system that generalizes the Bost--Connes
system to rank one Drinfeld modules, the scaling action on the dual
system can be described in terms of the Frobenius and inertia
groups.

In the present paper we continue along this line of thought. We
begin by reviewing the main steps in the Weil proof for function
fields, where we highlight the main conceptual steps and the main
notions that will need an analog in the noncommutative geometry
setting. We conclude this part by introducing the main entries in
our still tentative dictionary. The rest of the paper discusses in
detail some parts of the dictionary and provides evidence in support
of the proposed comparison. We begin this part by recalling briefly
the properties of the Bost--Connes endomotive from \cite{CCM}
followed by the description of the ``restriction map'' corresponding
to the inclusion of the ideles class group $C_\K=\A_\K^*/\K^*$ in
the noncommutative adeles class space $X_\K=\A_\K/\K^*$. We discuss
its relation to the exact sequence of Hilbert spaces of
\cite{Co-zeta} that plays a crucial role in obtaining the spectral
realization as an ``absorption spectrum''.

We then concentrate on the geometry of the adeles class space over
an arbitrary global field and the restriction map in this general
setting, viewed as a map of cyclic modules. We introduce the actions
$\vartheta_a$ and $\vartheta_m$ (with $a$ and $m$ respectively for
additive and multiplicative) of $\A_\K^*$ on suitable function
spaces on $\A_\K$ and on $C_\K$ and the induced action on the
cokernel of the restriction map in the category of cyclic modules.
We prove the corresponding general form of the associated Lefschetz
trace formula, as a cohomological reformulation of the trace formula
of \cite{Co-zeta} using the analytical setting of \cite{Meyer}.

The form of the trace formula and the positivity property that is
equivalent, in this setting, to the Riemann Hypothesis for the
corresponding $L$-functions with Gr\"ossencharakter, suggest by
comparison with the analogous notions in the Weil proof a natural
candidate for the analog of the Frobenius correspondence on the
curve. This is given by the graph of the scaling action. We can also
identify the analog of the degree and co-degree of a correspondence,
and the analog of the self intersection of the diagonal on the
curve, by looking at the explicit form of our Lefschetz trace
formula. We also have a clear analog of the first step in the Weil
proof of positivity, which consists of adjusting the degree by
multiples of the trivial correspondences. This step is possible,
with our notion of correspondences, due to a subtle failure of
Fubini's theorem that allows us to modify the degree by adding
elements in the range of the ``restriction map'', which play in this
way the role of the trivial correspondences. This leaves open the
more difficult question of identifying the correct analog of the
principal divisors, which is needed in order to continue the
dictionary.

We then describe how to obtain a good analog of the algebraic points
of the curve in the number field case (in particular in the case of
$\K=\Q$), in terms of the thermodynamical properties of the system.
This refines the general procedure described in \cite{CCM}. In fact,
after passing to the dual system, one can consider the periodic
orbits. We explain how, by the result of \cite{Co-zeta}, these are
the noncommutative spaces where the geometric side of the Lefschetz
trace formula concentrates. We show that, in turn, these periodic
orbits carry a time evolution and give rise to quantum statistical
mechanical systems, of which one can consider the low temperature
KMS states. To each periodic orbit one can associate a set of
``classical points'' and we show that these arise as extremal low
temperature KMS states of the corresponding system. We show that, in
the function field case, the space obtained in this way indeed can
be identified, compatibly with the Frobenius action, with the
algebraic points of the curve, albeit by a non-canonical
identification. Passing to the dual system is the analog in
characteristic zero of the transition from $\F_q$ to its algebraic
closure $\bar \F_q$. Thus, the procedure of considering periodic
orbits in the dual system and classical points of these periodic
orbits can be seen as an analog, for our noncommutative space, of
considering points defined over the extensions $\F_{q^n}$ of $\F_q$
in the case of varieties defined over finite fields (\cf \cite{CCM}
and \S 4 of \cite{CMbook}).

We analyze the behavior of the adeles class space under field
extensions and the functoriality question. We then finish the paper
by sketching an analogy between some aspects of the geometry of the
adeles class spaces and the theory of singularities, which may be
useful in adapting to this context some of the techniques of
vanishing and nearby cycles.

\section{A look at the Weil proof}\label{WeilproofSect}

In this preliminary section, we briefly review some aspects of
the Weil proof of the Riemann Hypothesis for function fields, with an
eye on extending some of the basic steps and concept to a
noncommutative framework, which is what we will be doing in the rest
of the paper.

\smallskip

In this section we let $\K$ be a global field of positive
characteristic $p>0$. One knows that, in this case, there
exists a smooth projective curve over a finite field $\F_q$, with
$q=p^r$ for some $r\in \N$, such that
\begin{equation}\label{functionfieldK}
\K =\F_q(C)
\end{equation}
is the field of functions of $C$. For this reason, a global field of
positive characteristic is called a function field.

\smallskip

We denote by $\Sigma_\K$ the set of places of $\K$. A place $v\in
\Sigma_\K$ is a Galois orbit of points of $C(\bar\F_q)$. The degree
$n_v=\deg(v)$ is its cardinality, namely the number of points in the
orbit of the Frobenius acting on the fiber of the natural map
from points to places
\begin{equation}\label{ptstoplaces}
C(\bar \F_q) \to \Sigma_\K .
\end{equation}
This means that the fiber over $v$ consists of $n_v$ conjugate
points defined over $\F_{q^{n_v}}$, the residue field of the local
field $\K_v$.

\smallskip

The curve $C$ over $\F_q$ has a zeta function of the form
\begin{equation}\label{zetaCT}
Z_C(T)=\exp\left( \sum_{n=1}^\infty \frac{\# C(\F_{q^n})}{n}\, T^n
\right),
\end{equation}
with $\log Z_C(T)$ the generating function for the number of points of
$C$ over the fields $\F_{q^n}$. It is customary to use the notation
\begin{equation}\label{zetaCs}
\zeta_\K(s)=\zeta_C(s)=Z_C(q^{-s}).
\end{equation}
It converges for $\Re(s)>1$. In terms of Euler product expansions one
writes
\begin{equation}\label{zetaKEuler}
\zeta_\K(s)=\prod_{v\in \Sigma_\K} (1-q^{-{n_v}s})^{-1}.
\end{equation}
In terms of divisors of $C$, one has equivalently
\begin{equation}\label{sumNDs}
\zeta_\K(s)=\zeta_C(s)= \sum_{D\geq 0} N(D)^{-s},
\end{equation}
where the norm of the divisor $D$ is $N(D)=q^{\deg(D)}$.

\smallskip

The Riemann--Roch formula for the curve $C$ states that
\begin{equation}\label{CRRthm}
\ell(D)-\ell(\kappa_C-D)=\deg(D)-g+1,
\end{equation}
where $\kappa_C$ is the canonical divisor on $C$, with degree
$\deg(\kappa_C)=2g-2$ and $h^0(\kappa_C)=g$, and $\ell(D)$ the
dimension of $H^0(D)$. Both $\deg(D)$ and $N(D)$ are well defined on the
equivalence classes obtained by adding principal divisors, that is,
\begin{equation}\label{modprincdiv}
D \sim D' \ \ \ \Longleftrightarrow \ \ \ D-D' = (f),
\end{equation}
for some $f\in \K^*$.
The Riemann--Roch formula \eqref{CRRthm} also implies that the zeta
function $\zeta_\K(s)$ satisfies the functional equation
\begin{equation}\label{zetaKfuncteq}
\zeta_\K(1-s)=q^{(1-g)(1-2s)} \zeta_\K(s).
\end{equation}
The zeta function $\zeta_\K(s)$ can also be written as a rational
function
\begin{equation}\label{prodzeta3}
Z_C(T) =\frac{P(T)}{(1-T)(1-qT)}\,,\ \ \ T=q^{-s}\,,
\end{equation}
where $P(T)$ is a polynomial of degree $2g$ and integer coefficients
\begin{equation}\label{PTzeta}
P(T) =\prod_{j=1}^{2g} (1-\lambda_j T).
\end{equation}
In particular, one has
\begin{equation}\label{countptsPT}
\# C(\F_{q^n})= q^n +1 - \sum_{j=1}^{2g} \lambda_j.
\end{equation}

Another important reformulation of the zeta function can be given in
terms of \'etale cohomology. Namely, the coefficients $\# C(\F_{q^n})$ that
appear in the zeta function can be rewritten as
\begin{equation}\label{fixFrobn}
\# C(\F_{q^n}) = \# {\rm Fix}(\Fr^n: \bar C \to \bar C),
\end{equation}
with $\bar C=C\otimes_{\F_q}\bar \F_q$.
The {\em Lefschetz fixed point formula} for \'etale cohomology then
shows that
\begin{equation}\label{fixFrobnLef}
\# C(\F_{q^n}) = \sum_{i=0}^2 (-1)^i \Tr\left(\Fr^n |\, H^i_{et}(\bar
C,\Q_\ell)\right).
\end{equation}
Thus, the zeta function can be written in the form
\begin{equation}\label{zetadetFrob}
\zeta_\K(s)=\prod_{i=0}^2 \left( \exp\left( \sum_{n=1}^\infty
\Tr(\Fr^n |\, H^i_{et}(\bar C,\Q_\ell)) \frac{q^{-sn}}{n} \right)
\right)^{(-1)^i}.
\end{equation}

\smallskip

The analog of the Riemann hypothesis for the zeta functions
$\zeta_\K(s)$ of function fields was stated in 1924 by E.~Artin as the
property that the zeros lie on the line $\Re(s)=1/2$. Equivalently, it
states that the complex numbers $\lambda_j$ of \eqref{PTzeta}, which
are are the eigenvalues of the Frobenius acting on
$H^1_{et}(\bar C,\Q_\ell)$, are algebraic numbers satisfying
\begin{equation}\label{RHff}
| \lambda_j |=\sqrt{q}.
\end{equation}

\smallskip

The Weil proof can be formulated either using
\'etale cohomology, or purely in terms of the Jacobian of the curve, or
again (equivalently) in terms of divisors on $C\times C$. We follow
this last viewpoint and we recall in detail some of the main steps in
the proof.

\subsection{Correspondences and divisors}\label{CorrdivSect}

Correspondences $Z$, given by (non-vertical) divisors on $C\times C$, form a ring
under composition
\begin{equation}\label{Z1circZ2}
Z_1\star Z_2 =(p_{13})_* (p_{12}^* Z_1
\bullet p_{23}^* Z_2),
\end{equation}
with $p_{ij}: C\times C\times C \to C\times C$ the projections, and
$\bullet$ the intersection product. The ring has an involution
\begin{equation}\label{adjZcorr}
Z'= \sigma(Z)
\end{equation}
where $\sigma$ is the involution that exchanges the two copies of $C$
in the product $C\times C$.

\smallskip

The degree $d(Z)$ and the codegree $d'(Z)$ are defined as the
intersection numbers
\begin{equation}\label{2degs}
d(Z)=Z \bullet (P\times C)\ \ \ \text{ and } \ \ \  d'(Z)=Z \bullet
(C\times P), \ \ \ \forall P\in C.
\end{equation}
They satisfy the relations
\begin{equation}\label{degstarcorr}
d(Z')=d'(Z), \ \ \text{ and } \ \  d(Z_1\star Z_2)=d(Z_1)d(Z_2).
\end{equation}

\smallskip

The correspondences $P\times C$ and $C\times P$ are called {\em
trivial correspondences}. One can consider the abelian group
${\rm Div}_{tr}(C\times C)$ generated by these trivial correspondence
and take the quotient
\begin{equation}\label{corrmodtriv}
\cC(C) := {\rm Div}(C\times C) / {\rm Div}_{tr}(C\times C) .
\end{equation}
It is always possible to change the degree and codegree of a
correspondence $Z$ by adding a multiple of the trivial correspondences
$P\times C$ and $C\times P$, so that, for any element in $\cC$ we find
a representative $Z\in {\rm Corr}$ with
\begin{equation}\label{deg0}
d(Z)=d'(Z)=0.
\end{equation}
One also wants to consider correspondences up to linear equivalence,
\begin{equation}\label{Z12equivf}
Z_1 \sim Z_2 \Longleftrightarrow Z_1 - Z_2 = (f),
\end{equation}
where $(f)$ is a principal divisor on $C\times C$.
Thus, one can consider $${\rm Pic}(C\times C)={\rm Div}(C\times C)
/\sim$$ and its quotient $\cP(C)$ modulo the classes of the trivial
correspondences.

\smallskip

A correspondence $Z$ is effective if it is given by an effective
divisor on $C\times C$, namely if it is a combination $Z=\sum_i n_i
Z_i$ of curves $Z_i\subset C\times C$ with coefficients $n_i\geq 0$.
We write $Z\geq 0$ to mean its effectiveness.
An effective correspondence $Z\geq 0$ that is nonempty can be viewed 
as a multivalued map
\begin{equation}\label{multivmapZ}
Z:C\to C, \ \ \ \ P\mapsto Z(P),
\end{equation}
with $Z(P) = proj_C(Z\bullet (P\times C))$,
of which the divisor $Z\subset C\times C$ is the graph, 
and with the product \eqref{Z1circZ2} given by the composition.

\smallskip

The trace of a correspondence $Z$ on $C\times C$ is the expression
\begin{equation}\label{TraceCorr}
\Tr(Z)=d(Z)+d'(Z)-Z\bullet \Delta,
\end{equation}
with $\Delta$ the diagonal (identity correspondence) and $\bullet$ the
intersection product. This is well defined on $\cP(C)$ since
$\Tr(Z)=0$ for principal divisors and trivial correspondences.

\smallskip

Consider a correspondence of degree $g$ of the form $Z=\sum a_n
\Fr^n$, given by a combination of powers of the Frobenius. Then $Z$
can be made effective by adding a principal correspondence
which is defined over $\F_q$ and which commutes with $\Fr$.

\smallskip

This can be seen as follows. The Riemann--Roch theorem ensures that
one can make $Z$
effective by adding a principal correspondence, over the field
$k(P)$, where $k$ is the common field of definition of the
correspondence $Z$ and of the curve (\cf \cite{SD}) and $P$
is a generic point. A correspondence of the form $Z=\sum a_n \Fr^n$
is in fact defined over $\F_q$ hence the principal correspondence is
also defined over $\F_q$. As such it automatically commutes with $\Fr$
(\cf \cite{Weil-lett}, p. 287).

\smallskip

Notice however that, in general, it is not possible to modify a
divisor $D$ of degree one on $C$ to an effective divisor in such a way
that the added pricipal divisor has support on the same Frobenius
orbit. An illustrative example is given in Chapter 4 of
\cite{CMbook}.

\subsection{The explicit formula}\label{expformulaSect}

The main steps in the Weil proof of RH for function fields are
\begin{enumerate}
\item The explicit formula
\item Positivity
\end{enumerate}

Let $\K$ be a global field and let $\A_\K$ denote its ring of
adeles. Let $\Sigma_\K$ denote the set of places of $\K$.
Let $\alpha$ be a non-trivial character of $\A_\K$ which
is trivial on $\K\subset \A_\K$. We write
\begin{equation}\label{alphacharv}
\alpha=\prod_{v\in \Sigma_\K} \alpha_v,
\end{equation}
for the decomposition of $\alpha$ as a product of its restrictions
to the local fields $\alpha_v=\alpha|_{\K_v}$.

Consider the bicharacter
\begin{equation}\label{bicharpair}
\langle z ,\lambda \rangle :=\lambda^z,
\ \ \ \text{ for} \ (z,\lambda)\in \C \times \R^*_+ .
\end{equation}
Let $N$ denote the range of the norm $|\cdot|:C_\K \to
\R^*_+$.
Then $N^\bot\subset \C$ denotes the subgroup
\begin{equation}\label{Nbot}
N^\bot := \{ z\in \C |\, \lambda^z=1, \ \forall \lambda\in N \}.
\end{equation}

Consider then the expression
\begin{equation}\label{spectralside}
 \sum_{\rho\in \C/N^\bot |\,  L(\tilde\chi,
\rho) =0} \hat f (\tilde\chi,\rho),
\end{equation}
with $L(\tilde\chi,\rho)$ the $L$-function with Gr\"ossencharakter
$\chi$, where $\tilde\chi$ denotes the extension to $C_\K$ of the
character $\chi\in \widehat{C_{\K,1}}$, the Pontrjagin dual of $C_{\K,1}$.
Here $\hat f (\tilde\chi,\rho)$ denotes the Fourier transform
\begin{equation}\label{rhoFourier}
\hat f (\tilde\chi,\rho) = \int_{C_\K} f(u) \tilde\chi(u) \,
|u|^\rho \, d^*u 
\end{equation}
of a test function $f$ in the Schwartz space $\cS(C_\K)$.

In the case where $N=q^\Z$ (function fields), the subgroup
$N^\bot$ is nontrivial and given by
\begin{equation}\label{ffNbot}
N^\bot = \frac{2\pi i}{\log q} \Z.
\end{equation}
Since in the function field case the $L$-fuctions are functions of
$q^{-s}$, there is a periodicity by $N^\bot$, hence we need to
consider only $\rho\in \C/N^\bot$. In the number
field case this does not matter, since $N=\R^*_+$ and $N^\bot$ is
trivial.

The Weil explicit formula is the remarkable identity \cite{weilexplicit}
\begin{equation}\label{Weilexpl}
\hat h (0) + \hat h (1) - \sum_{\rho\in \C/N^\bot |\,  L(\tilde\chi,
\rho) =0} \hat h (\tilde\chi,\rho) = \sum_{v\in \Sigma_\K}
\int_{(\K^*_v,\alpha_v)}' \frac{h(u^{-1})}{|1-u|}\, d^* u .
\end{equation}
Here the Fourier transform $\hat h$ is as in \eqref{rhoFourier}. The
test function $h$ has compact support and belongs to the Schwartz
space $\cS(C_\K)$. As soon as $h(1)\neq 0$ the integrals in the
right hand side are singular so that one needs to specify how to
take their principal value. This was done in \cite{weilexplicit} and
it was shown in \cite{Co-zeta} that the same principal value can in
fact be defined in the following unified way.

\begin{defn}
For a  local field $K$ and a given (non-trivial) additive character
$\beta$ of $K$, one lets $\varrho_\beta$ denote the unique
distribution extending $d^*u$ at $u=0$, whose Fourier transform
\begin{equation}\label{Fourierbeta}
\hat\varrho(y)=\int_K \varrho(x)\beta(xy)\,dx
\end{equation}
satisfies the vanishing condition $\hat\varrho (1) = 0$.
\end{defn}
Then by definition the principal value $\int'$ is given by
\begin{equation}\label{princvalrhobeta}
\int_{(K,\beta)}' \frac{f(u^{-1}) }{|1-u|} \, d^* u
=\langle \varrho_\beta,g\rangle, \ \ \ \text{ with } \ \
g(\lambda)=\frac{f((\lambda+1)^{-1})}{|\lambda+1|},
\end{equation}
where $\langle \varrho_\beta,g\rangle$ denotes the pairing of the
distribution $\varrho_\beta$ and the function $g(\lambda)$. This
makes sense provided the support of $f$ is compact which implies
that $g(\lambda)$ vanishes identically in a neighborhood of
$\lambda=-1$.

The Weil explicit formula is a far reaching generalization of the
relation between primes and zeros of the Riemann zeta function,
originally due to Riemann \cite{Riemann}.

\subsection{Riemann--Roch and positivity}\label{PositSect}

Weil positivity is the statement that, if $Z$ is a nontrivial
correspondence in $\cP(C)$ (\ie as above a correspondence on $C\times
C$ modulo trivial ones and up to linear equivalence), then
\begin{equation}\label{posTrZZ}
\Tr (Z\star Z') >0.
\end{equation}

This is proved using the Riemann--Roch
formula on $C$ to show that one can achieve effectivity. In fact,
using trivial correspondences to adjust the degree one can assume
that $d(Z)=g$. Then the Riemann--Roch formula \eqref{CRRthm} shows
that if $D$ is a divisor on $C$ of degree $\deg(D)=g$
then there are effective representatives in the
linear equivalence class of $D$.
The intersection of $Z\subset C\times C$ with $P\times C$ defines a
divisor $Z(P)$ on $C$ with $$\deg(Z(P))=d(Z)=g.$$ Thus, the argument
above shows that there exists $f_P\in \K^*$ such that $Z(P)+(f_P)$ is
effective. This determines an effective divisor $Z+(f)$ on $C\times C$.
Thus, we can assume that $Z$ is effective, hence 
we can write it as a multivalued function
\begin{equation}\label{Zmultifunct}
P\mapsto Z(P)=Q_1 +\cdots + Q_g.
\end{equation}

\smallskip

The product $Z\star Z'$ is of the form
\begin{equation}\label{DeltaYprod}
Z\star Z' = d'(Z) \Delta + Y,
\end{equation}
where $\Delta$ is the diagonal in $C\times C$ and $Y$ is the
effective correspondence such that $Y(P)$ is the divisor on $C$ given
by the sum of points in $$\{ Q\in C |\, Q=Q_i(P)=Q_j(P), \, i\neq j
\}.$$

One sees this from the description
in terms of intersection product that it is given by the multivalued
function
\begin{equation}\label{ZstarZprime}
(Z\star Z')(Q)= \sum_{i,j} \sum_{P\in \cU_{ij}(Q)} P,
\end{equation}
where
$$ \cU_{ij}(Q)=\{ P\in C |\, Q_i(P)=Q_j(P)=Q \}. $$
One can separate this out in the contribution of the locus where
$Q_i=Q_j$ for $i\neq j$ and the part where $i=j$,
$$  (Z\star Z')(Q)= U(Q) + Y(Q). $$
Notice that
\begin{equation}\label{degprimeQij}
\# \{ P\in C |\, Q=Q_i(P), \text{ for some }  i=1,\ldots, g \} = d'(Z).
\end{equation}
Thus, for $i=j$ we obtain that the divisor $U(Q)=\sum_i
\sum_{P\in\cU_{ii}(Q)} P$ is just $d'(Z)\Delta(Q)$, while for $i\neq
j$ one obtains the remaining term $Y$ of \eqref{DeltaYprod}.

\smallskip

In the case $g=1$, the effective correspondence $Z(P)=Q(P)$ is single
valued and the divisor $(Z\star Z')(P)$ of \eqref{ZstarZprime} reduces
to the sum of points in $$\cU(Q)=\{ P\in C |\, Q(P)=Q \}.$$ There are
$d'(Z)$ such points so one obtains
\begin{equation}\label{g1caseZZ}
Z\star Z' =d'(Z) \Delta, \ \ \ \text{ with } \ \ \ \Tr(Z\star Z')=
2d'(Z) \geq 0,
\end{equation}
since for $g=1$ one has $\Delta\bullet\Delta =0$ and $d'(Z) \geq 0$
since $Z$ is effective.

\smallskip

In the case of genus $g>1$, the Weil proof proceeds as follows. Let
$\kappa_C$ be a choice of an effective canonical divisor for $C$
without multiple points, and let $\{ f_1, \ldots, f_g \}$ be a basis
of the space $H^0(\kappa_C)$. One then considers the function
$C\to M_{g\times g}(\F_q)$ to $g\times g$ matrices
\begin{equation}\label{matrfQ}
 P \mapsto M(P), \ \ \ \text{ with } \ \ M_{ij}(P)= f_i(Q_j(P)).
\end{equation}
and the function $K: C\to \F_q$ given by
\begin{equation}\label{detfQ}
K(P)= \det(M(P))^2.
\end{equation}

\smallskip

The function $P\mapsto K(P)$ of \eqref{detfQ} is a rational function
with $(2g-2)d'(Z)$ double poles. In fact, $K(P)$ is a symmetric
function of the $Q_j(P)$, because of the squaring of the
determinant. The composition $P\mapsto (Q_j(P))\mapsto
K(P)$ is then a rational function of $P\in C$.
The poles occur (as double poles) at those points $P\in C$ for which
some $Q_i(P)$ is a component of $\kappa_C$. The canonical divisor
$\kappa_C$ has degree $2g-2$. This means that there are $(2g-2)d'(Z)$
such double poles.

\smallskip

For $Z\star Z' = d'(Z) \Delta + Y$ as above, the
intersection number $Y\bullet \Delta$ satisfies the estimate
\begin{equation}\label{YintDelta}
Y\bullet \Delta \leq (4g-4)\, d'(Z).
\end{equation}
In fact, the rational function $K(P)$ of
\eqref{detfQ} has a number of zeros equal to $(4g-4)\, d'(Z)$.
On the other hand, $Y\bullet \Delta$ counts the number of times that
$Q_i=Q_j$ for $i\neq j$. Since each point $P$ with $Q_i(P)=Q_j(P)$ for
$i\neq j$ produces a zero of $K(P)$, one sees that $Y\bullet \Delta$
satisfies the estimate \eqref{YintDelta}.
Notice that, for genus $g>1$ the self intersection of the diagonal is
the Euler characteristic
\begin{equation}\label{gDeltaDelta}
\Delta \bullet \Delta = 2-2g = \chi(C).
\end{equation}

\smallskip

Moreover, we have
\begin{equation}\label{ZZdd}
d(Z\star Z')=d(Z)d'(Z)= g\, d'(Z) = d'(Z\star Z').
\end{equation}
Thus, using again the decomposition \eqref{DeltaYprod} and the
definition of the trace of a correspondence \eqref{TraceCorr},
together with \eqref{gDeltaDelta} and \eqref{ZZdd} one obtains
\begin{equation}\label{positZZ}
\begin{array}{c}
\Tr(Z\star Z')= 2g\, d'(Z) + (2g-2)\, d'(Z) - Y\bullet \Delta \\[3mm]
\geq (4g-2)\, d'(Z)-(4g-4)\, d'(Z) = 2d'(Z) \geq 0.
\end{array}
\end{equation}
This gives the positivity \eqref{posTrZZ}.

\medskip

In the Weil proof of RH for function fields, one concentrates on a
particular type of correspondences, namely those that are of the form
\begin{equation}\label{FrDeltaCorr}
Z_{n,m}= m \Delta + n\Fr,
\end{equation}
for $n,m\in \Z$, with $\Fr$ the Frobenius correspondence.

\smallskip

Notice that, while the correspondence depends linearly on $n,m\in \Z$,
the expression for the trace gives
\begin{equation}\label{TrnmFr}
\Tr(Z_{n,m}\star Z_{n,m}')= 2 g m^2+2(1+q-N)mn+2g q n^2,
\end{equation}
where $N=\# C(\F_q)$. In particular, \eqref{TrnmFr} depends
quadratically on $(n,m)$.  In the process of passing from a
correspondence of degree $g$ to an effective correspondence, this
quadratic dependence on $(n,m)$ is contained in the multiplicity
$d'(Z)$. Notice, moreover, that the argument does not depend on the torsion
part of the ring of correspondences.

\subsection{A tentative dictionary}

In the rest of the paper we illustrate some steps towards the creation
of a dictionary relating the main steps in the Weil proof described
above to the noncommutative geometry of the adeles class space of a
global field. The noncommutative geometry approach has the advantage
that it provides (see \cite{Co-zeta}, \cite{Meyer}, \cite{CCM}) a
Lefschetz trace formula interpretation for the Weil explicit formula
and that it gives a parallel formulation for both function fields and
number fields. Parts of the dictionary sketched below are very
tentative at this stage, so we mostly concentrate, in the rest of the
paper, on illustrating what we put in the first few lines of the
dictionary, on the role of the scaling correspondence as Frobenius and
its relation to the explicit formula.

\bigskip

\begin{center}
\begin{tabular}{|c|c|}
\hline
& \\
Frobenius correspondence &  $Z(f)=\int_{C_\K} f(g)\, Z_g \, d^*g$ \\
&  \\
\hline
& \\
Trivial correspondences & Elements of the range $\cV$ \\
& \\
\hline
& \\
Adjusting the degree  & Fubini step  \\
by trivial correspondences & on the test functions \\
& \\
\hline
& \\
Correspondences & Bivariant elements $Z(f)\Rightarrow \Gamma(f)$ \\
& \\
\hline
& \\
Degree of a correspondence & Pointwise index  \\
& \\
\hline
& \\
Riemann--Roch & Index theorem \\
&  \\
\hline
& \\
Effective correspondences & Epimorphism of $C^*$-modules \\
& \\
\hline
& \\
$\deg Z(P)\geq g \Rightarrow$ $Z+(f)$ effective & $d(\Gamma)>0
\Rightarrow$ $\Gamma+K$ onto  \\
& \\
\hline
& \\
Lefschetz formula & Bivariant Chern  of $\Gamma(f)$ \\
& (by localization on the graph $Z(f)$) \\
& \\
\hline
\end{tabular}
\end{center}

\section{Quantum statistical mechanics and arithmetic}\label{QSMsect}

The work of Bost--Connes \cite{BC} first revealed the presence of an
interesting interplay between quantum statistical mechanics and
Galois theory. More recently, several generalizations \cite{CM},
\cite{CMR}, \cite{CMR2}, \cite{CCM}, \cite{ConsMar-ff},
\cite{HaPau}, \cite{Jacob} have confirmed and expanded this
viewpoint. The general framework of interactions between
noncommutative geometry and number theory described in \cite{Man1},
\cite{Man2}, \cite{CMbook}, \cite{Mar} recast these phenomena into a
broader picture, of which we explore in this paper but one of many
facets.

The basic framework that combines quantum statistical mechanics and
Galois theory can be seen as an extension, involving noncommutative
spaces, of the category of Artin motives. In the setting of pure
motives (see \cite{Man-mot}), Artin motives correspond to the
subcategory generated by zero dimensional objects, with morphisms
given by algebraic cycles in the product (in this case without the
need to specify with respect to which equivalence relation).
Endomotives were introduced in \cite{CCM} as noncommutative spaces
of the form
\begin{equation}\label{algdef}
\cA_\K=\,A\rtimes S ,
\end{equation}
where $A$ is an inductive limit of reduced finite dimensional
commutative algebras over the field $\K$, \ie a projective limit of
Artin motives, and $S$ is a unital abelian semigroup of algebra
endomorphisms $\rho:  A \to A$. These have the following properties:
the algebra $A$ is unital; the image $e=\rho(1)\in A$ is an
idempotent, for all $\rho \in S$; each $\rho\in S$ is an isomorphism
of $A$ with the compressed algebra $eAe$. A general construction
given in \cite{CCM} based on self maps of algebraic varieties
provides a large class of examples over different fields $\K$. We
are mostly interested here in the case where $\K$ is a number field
and for part of our discussion below we will concentrate on a
special case (the Bost--Connes endomotive) over the field $\K=\Q$.

Endomotives form a pseudo-abelian category where morphisms are
correspondences given by $\cA_\K$--$\cB_\K$-bimodules that are
finite and projective as right modules. These define morphisms in
the additive KK-category and in the abelian category of cyclic
modules. In fact, in addition to the algebraic form described above,
endomotives also have an analytic structure given by considering,
instead of the $\K$-algebra \eqref{algdef} the $C^*$-algebra
\begin{equation}\label{andef}
C(X)\rtimes S,
\end{equation}
where $X$ denotes the totally disconnected Hausdorff space $X=X(\bar
\K)$ of algebraic points of the projective limit of Artin motives.
There is a canonical action of the Galois group $G=\Gal(\bar\K/\K)$
by automorphisms of the $C^*$-algebra \eqref{andef} globally
preserving $C(X)$. We refer the reader to \cite{CCM} for a more
detailed discussion of algebraic and analytic endomotives and the
properties of morphisms in the corresponding categories.

If the endomotive is ``uniform'' in the sense specified in
\cite{CCM}, the space $X$ comes endowed with a probability measure
$\mu$ that induces a state $\varphi$ on the $C^*$-algebra
\eqref{andef}. The general Tomita theory of modular automorphism
groups in the context of von Neumann algebras \cite{tt} shows that
there is a natural time evolution for which the state $\varphi$ is
KMS$_1$. One can then consider the set $\Omega_\beta$ of low
temperature KMS states for the resulting quantum statistical
mechanical system. One also associates to the system $(\cA,\sigma)$
of the $C^*$-algebra with the time evolution its dual system
$(\hat\cA,\theta)$, where the algebra $\hat\cA=\cA\rtimes_\sigma \R$
is obtained by taking the crossed product with the time evolution
and $\theta$ is the scaling action of $\R^*_+$
\begin{equation}\label{dualaction}
\theta_\lambda(\int x(t)\,U_t\,dt)=\int \lambda^{it}\,x(t)\,U_t\,dt.
\end{equation}
One then constructs an $\R^*_+$ equivariant map
\begin{equation}\label{cool}
\pi : \hat \cA_\beta \to C(\tilde \Omega_\beta, \cL^1),
\end{equation}
from a suitable subalgebra $\hat\cA_\beta\subset \hat\cA$ of the
dual system to functions on a principal $\R^*_+$-bundle
$\tilde\Omega_\beta$ over the low temperature KMS states of the
system, with values in trace class operators. Since traces define
morphisms in the cyclic category, the map \eqref{cool} can be used
to construct a morphism $\delta=(\Tr\circ\pi)^\natural$ at the level of cyclic modules
\begin{equation}\label{cyclcool}
\hat\cA_\beta^\natural
\stackrel{(\Tr\circ\pi)^\natural}{\longrightarrow}
C(\tilde\Omega_\beta)^\natural .
\end{equation}
This map can be loosely thought of as a ``restriction map''
corresponding to the inclusion of the ``classical points'' in the
noncommutative space. One can then consider the cokernel of this map
in the abelian category of cyclic modules. In \cite{CCM} we denoted
the procedure described above ``cooling and distillation'' of
endomotives. We refer the reader to \cite{CCM} for the precise
technical hypotheses under which this procedure can be performed.
Here we only gave an impressionistic sketch aimed at recalling
briefly the main steps involved.

\subsection{The Bost--Connes endomotive}

The main example of endomotive we will consider here in relation to
the geometry of the adeles class space is the Bost--Connes system.
This can be constructed as an endomotive over $\K=\Q$, starting from
the projective system $X_n=\Sp(A_n)$, with $A_n=\Q[\Z/n\Z]$ the
group ring of $\Z/n\Z$. The inductive limit is the group ring
$A=\Q[\Q/\Z]$ of $\Q/\Z$. The endomorphism $\rho_n$ associated to an
element $n\in S$ of the (multiplicative) semigroup $S=\N=\Z_{>0}$ is
given on the canonical basis $e_r \in \Q[\Q/\Z]$, $r\in \Q/\Z$, by
\begin{equation} \label{rhoBC}
\rho_n(e_r)=\,\frac{1}{n}\,\sum_{ns=r}\,e_s
\end{equation}
The corresponding analytic endomotive is the crossed product
$C^*$-algebra $$\cA=C^*(\Q/\Z)\rtimes \N.$$ The Galois action is given by
composing a character $\chi: A_n \to \bar\Q$ with an element $g$ of
the Galois group $G=\Gal(\bar\Q/\Q)$. Since $\chi$ is determined by
the $n$-th root of unity $\chi(e_{1/n})$, one obtains the cyclotomic
action.

In the case of the Bost--Connes endomotive, the state $\varphi$ on
$\cA$ induced by the measure $\mu$ on $X=\hat\Z$ is of the form
\begin{equation}\label{BCstateKMS1}
\varphi(f) = \int_{\hat\Z} f(1,\rho)\, d\mu(\rho),
\end{equation}
and the modular automorphism group restricts to the $C^*$-algebra as
the time evolution of the BC system, \cf \cite{BC}, \cite{CCM} and
\S 4 of \cite{CMbook}.

The dual system of the Bost--Connes system is best described in
terms of commensurability classes of $\Q$-lattices. In \cite{CM} the
Bost--Connes system is reinterpreted as the noncommutative space
describing the relation of commensurability for 1-dimensional
$\Q$-lattices up to scaling. One can also consider the same
equivalence relation without dividing out by the scaling action. If
we let $\cG_1$ denote the groupoid of the commensurability relation
on 1-dimensional $\Q$-lattices and $\cG_1/\R^*_+$ the one obtained
after moding out by scaling, we identify the $C^*$-algebra of the
Bost--Connes system with $C^*(\cG_1/\R^*_+)$ (\cf \cite{CM}). The
algebra $\hat\cA$ of the dual system is then obtained in the
following way (\cf \cite{CCM}). There is a $C^*$-algebra isomorphism
$\iota: \hat\cA \to C^*(\cG_1)$ of the form
\begin{equation}\label{iotamap}
\iota( X )(k,\rho,\lambda)=\int_\R x(t)(k,\rho)\, \lambda^{it}\, dt
\end{equation}
for $(k,\rho,\lambda)\in \cG_1$ and $X =\int x(t)\, U_t\, dt \in
\hat\cA$.

\subsection{Scaling as Frobenius in characteristic
zero}\label{Frob0sect}

In the general setting described in \cite{CCM} one denotes by
$D(\cA,\varphi)$ the cokernel of the morphism \eqref{cyclcool},
viewed as a module in the cyclic category. The notation is meant to
recall the dependence of the construction on the initial data of an
analytic endomotive $\cA$ and a state $\varphi$. The cyclic module
$D(\cA,\varphi)$ inherits a scaling action of $\R^*_+$ and one can
consider the induced action on the cyclic homology
$HC_0(D(\cA,\varphi))$. We argued in \cite{CCM} that this cyclic
homology with the induced scaling action plays a role analogous to
the role played by the Frobenius action on \'etale cohomology in the
algebro-geometric context. Our main supporting evidence is the
Lefschetz trace formula for this action that gives a cohomological
interpretation of the spectral realization of the zeros of the
Riemann zeta function of \cite{Co-zeta}. We return to discuss the
Lefschetz trace formula for the more general case of global fields
in \S \ref{cohTrSect} below.

\medskip

The main results of \cite{Co-zeta} show that we have the following
setup. There is an exact sequence of Hilbert spaces
\begin{equation}\label{seqL2spaces}
0\to L^2_\delta(\A_\Q/\Q^*)_0\to L^2_\delta(\A_\Q/\Q^*) \to \C^2 \to 0,
\end{equation}
which defines the subspace $L^2_\delta(\A_\Q/\Q^*)_0$ by imposing
the conditions $f(0)=0$ and $\hat f(0)=0$ and a suitable decay
condition imposed by the weight $\delta$. The space
$L^2_\delta(\A_\Q/\Q^*)_0$ fits into another exact sequence of Hilbert
spaces of the form
\begin{equation}\label{mapfEagain}
 0 \to L^2_\delta(\A_\Q/\Q^*)_0\stackrel{\fE}{\to}
L^2_\delta(C_\Q)\to \cH \to 0
\end{equation}
where the map $\fE$ is defined by
\begin{equation}\label{mapfEagain2}
 \fE(f)(g)=|g|^{1/2} \sum_{q\in \Q^*} f(qg), \ \  \ \ \forall g\in
C_\Q=\A_\Q^*/\Q^* .
\end{equation}
The map is equivariant with respect to the actions of $C_\Q$ \ie
\begin{equation}\label{CQact}
\fE\circ \urep(\gamma)=|\gamma|^{1/2}\vrep(\gamma)\circ \fE
\end{equation}
where $(\urep(\gamma)\xi)(x)=\xi(\gamma^{-1}x)$ for $\xi \in L^2_\delta(\A_\Q/\Q^*)_0$ and similarly $\vrep(\gamma)$ is the regular representation of $C_\K$.

\smallskip

We showed in \cite{CCM} that the map $\fE$, translated from the
context of Hilbert spaces to that of nuclear spaces as in
\cite{Meyer}, has a natural interpretation in terms of the ``cooling
and distillation process'' for the BC endomotive. In fact, we showed
in \cite{CCM} that, if $(\cA,\sigma)$ denotes the BC system,
then the following properties hold.
\begin{enumerate}
\item For $\beta >1$ there is a canonical isomorphism
\begin{equation}\label{OmegaCQ}
\tilde \Omega_\beta \simeq \hat \Z^*\times \R_+^* \simeq C_\Q
\end{equation}
of $\tilde\Omega_\beta$ with the space of invertible 1-dimensional
$\Q$-lattices.
\item For $ X\in \hat \cA$ and $f=\iota( X) \in C^*(\cG_1)$,
the cooling map \eqref{cyclcool} takes the form
\begin{equation}\label{deltaBC}
\delta(X)(u,\lambda)=\sum_{n\in \N=\Z_{>0}}  f (1, n u,n\lambda), \ \ \
\forall (u,\lambda)\in C_\Q\simeq \tilde \Omega_\beta .
\end{equation}
\end{enumerate}

One can compare directly the right hand side of \eqref{deltaBC} with
the map $\fE$ (up to the normalization by $|j|^{1/2}$) written as in
\eqref{mapfEagain2} by considering a function
$f(\rho,v)=f(1,\rho,v)$ and its unique extension $\tilde f$ to
adeles  where $f$ is extended by $0$ outside $\hat \Z\times \R^*$
and one requires the parity
\begin{equation}\label{parit}
\tilde f(-u,-\lambda)=f(u,\lambda)\,.
\end{equation}
 This
gives then
\begin{equation}\label{estimate1}
\sum_{n\in \N}  f (1, n u,n\lambda)=\frac 12 \;\sum_{q\in \Q^*}
\tilde f (q \,j), \ \ \ \ \text{ where } j=(u,\lambda)\in C_\Q.
\end{equation}

\section{The adeles class space}\label{adclspSect}

Let $\K$ be a global field, with $\A_\K$ its ring of adeles.

\begin{defn}\label{adeleclassdef}
The adeles class space of a global field $\K$ is the quotient
$\A_\K/\K^*$.
\end{defn}

When viewed from the classical standpoint this  is a ``bad
quotient'' due to the ergodic nature of the action which makes the
quotient ill behaved topologically. Thus, following the general
philosophy of noncommutative geometry, we describe it by a
noncommutative algebra of coordinates, which allows one to continue
to treat the quotient as a ``nice quotient'' in the context of
noncommutative geometry.

A natural choice of the algebra is the cross product
\begin{equation}\label{crossKCAK}
C_0(\A_\K)\rtimes \K^* \ \ \text{ with the smooth subalgebra } \ \
\cS(\A_\K)\rtimes \K^*.
\end{equation}

A better description can be given in terms of groupoids.

Consider the groupoid law $\cG_\K=\K^*\ltimes \A_\K$ given by
\begin{equation}\label{gpdAK}
(k,x)\circ (k',y)= (k k',y), \ \ \ \forall k,k'\in
\K^*,\ \ \text{ and } \forall x,y\in \A_\K \ \text{ with } x=k' y,
\end{equation}
with the composition \eqref{gpdAK} defined
whenever the source $s(k,x)= x$ agrees with the range
$r(k',y)=k' y$.

\begin{lem}\label{groupoidAK}
The algebras \eqref{crossKCAK} are, respectively, the groupoid
$C^*$-algebra $C^*(\cG_\K)$ and its dense subalgebra $\cS(\cG_\K)$.
\end{lem}

\proof The product in the groupoid algebra is given by the associative
convolution product
\begin{equation}\label{convolGKalg}
 (f_1 * f_2) \, (k,x)= \sum_{s\in \K^*} f_1(k\,s^{-1}, s\,x)
f_2(s,x),
\end{equation}
and the adjoint is given by $f^*(k,x)=\overline{f(k^{-1},k\,x)}$.

The functions (on the groupoid) associated to $f \in \cS(\A_K)$ and
$U_k$ are given, respectively, by
\begin{equation}\label{crossgrpd}
\begin{array}{llll}
f(1,x)=f(x) & \text{and} & f(k,x)=0 & \forall k\neq 1 \\[2mm]
U_k(k,x)=1 &  \text{and} & U_g(k,x)=0 & \forall g\neq k .
\end{array}
\end{equation}
The product $f\,U_k$ is then the convolution product of the
groupoid.

The algebra $\cS(\cG_\K)$ is obtained by considering finite sums of the form
\begin{equation}
\sum_{k\in \K^*} f_k\,U_k, \ \ \text{ for } \  f_k \in \cS(\A_\K).
\end{equation}
The product is given by the convolution product
\begin{equation}\label{UkfUk}
(U_k\,f\,U_k^*)(x)= f(k^{-1} x),
\end{equation}
for $f \in \cS(\A_\K)$, $k\in \K^*$, and $x\in \A_\K$.
\endproof

\subsection{Cyclic module}\label{XKcyclSect}

We can associate to the algebra $\cS(\cG_\K)$ of the adeles class
space an object in the category of $\Lambda$-modules. This means
that we consider the cyclic module $\cS(\cG_\K)^\natural$ and the
two cyclic morphisms
\begin{equation}\label{varepsilonj}
\varepsilon_j : \cS(\cG_\K)^\natural \to \C
\end{equation}
given by
\begin{equation}\label{varepsjdeg0}
\varepsilon_0(\sum\,f_k\,U_k)= f_1(0) \ \ \ \text{ and } \ \ \
\varepsilon_1(\sum\,f_k\,U_k)=\int_{\A_\K} f_1(x)\,dx
\end{equation}
and in higher degree by
\begin{equation}\label{varepsjdeg}
\varepsilon_j^\natural(a^0\otimes \cdots \otimes a^n)=\,
\varepsilon_j (a^0\, \cdots \, a^n).
\end{equation}

The morphism $\varepsilon_1$ is given by integration on $\A_\K$ with
respect to the additive Haar measure. This is $\K^*$ invariant,
hence it defines a trace on $\cS(\cG_\K)$. In the case of $\K=\Q$,
this corresponds to the dual trace $\tau_\varphi$ for the
KMS$_1$-state $\varphi$ associated to the time evolution of the BC
system. The morphism $\varepsilon_0$ here takes into account the
fact that we are imposing a vanishing condition at $0\in \A_\K$ (\cf
\cite{CCM} and \cite{CMbook} Chapter 4). In fact, the
$\Lambda$-module we associate to $\cS(\cG_\K)$ is given by
\begin{equation}\label{cyclA0mod}
\cS(\cG_\K)^\natural_0 := \Ker\, \varepsilon_0^\natural \cap \Ker\,
\varepsilon_1^\natural.
\end{equation}
Note that since $\cS(\cG_\K)$ is non-unital, the cyclic module
$\cS(\cG_\K)^\natural$ is obtained using the adjunction of a unit to
$\cS(\cG_\K)$.

\subsection{The restriction map}\label{XKCKressect}

Consider the ideles $\A_\K^*=\GL_1(\A_\K)$ of $\K$ with their
natural locally compact topology induced by the map
\begin{equation}\label{gginvmap}
\A^*_\K \ni g \mapsto (g,g^{-1}).
\end{equation}
We can see the ideles class group $C_\K=\A_\K^*/\K^*$ as a subspace of
the adeles class space $X_\K=\A_\K/\K^*$ in the following way.

\begin{lem}\label{idelesubgrpd}
The pairs $((k,x),(k',y))\in \cG_\K$ such that both $x$ and $y$ are in
$\A^*_\K$ form a full subgroupoid of $\cG_\K$ which is isomorphic to
$\K^* \ltimes \A^*_\K$.
\end{lem}

\proof Elements of $\A_\K$ whose orbit under the $\K^*$ action
contains an idele are also ideles. Thus, we obtain a groupoid that is
a full subcategory of $\cG_\K$.
\endproof

This implies the existence of a restriction map. Consider the map
\begin{equation}\label{rho0restr}
\rho: \cS(\A_K) \ni f \mapsto f|_{\A^*_\K}.
\end{equation}
We denote by $C_\rho(\A^*_\K)\subset C(\A^*_\K)$ the range of $\rho$.

\begin{cor}\label{restrAKCK}
The restriction map $\rho$ of \eqref{rho0restr}
extends to an algebra homomorphism
\begin{equation}\label{rhoAKCKrestr}
\rho: \cS(\cG_\K) \to C_\rho(\A^*_\K)\rtimes \K^*.
\end{equation}
\end{cor}

\proof The map \eqref{rho0restr} induced by the inclusion
$\A_\K^*\subset \A_\K$ is continuous and $\K^*$ equivariant hence the
map
$$ \rho(\sum_{k\in \K^*} f_k\,U_k)=\sum_{k\in \K^*} \rho(f_k)\,U_k $$
is an algebra homomorphism.
\endproof

 The action of $\K^*$ on $\A_\K^*$ is free and proper so that
 we have an equivalence
of the locally compact groupoids $\K^*\ltimes \A_\K^*$ and
$\A_\K^*/\K^*=C_\K$. We use the exact sequence of locally compact
groups
\begin{equation}\label{akbykstar}
1\to \K^*\to \A_\K^*\stackrel{p}\to C_\K\to 1
\end{equation}
to parameterize the orbits of $\K^*$ as the fibers $p^{-1}(x)$ for
$x\in C_\K$. By construction the Hilbert spaces
\begin{equation}\label{hilbfiber}
\cH_x=\ell^2(p^{-1}(x)) \qqq x\in C_K
\end{equation}
form a continuous field of Hilbert spaces over $C_\K$. We let
$\cL^1(\cH_x)$ be the Banach algebra of trace class operators in
$\cH_x$, these form a continuous field over $C_\K$.

\begin{prop}\label{restrAKCK1}
The restriction map $\rho$ of \eqref{rho0restr} extends to an
algebra homomorphism
\begin{equation}\label{rhoAKCKrestrtr}
\rho: \cS(\cG_\K) \to C(C_\K,\cL^1(\cH_x))\,.
\end{equation}
\end{prop}
\proof Each $p^{-1}(x)$ is globally invariant under the action of
$\K^*$ so the crossed product rules in $C_\rho(\A^*_\K)\rtimes \K^*$
are just   multiplication of operators in $\cH_x$. To show that the
obtained operators are in $\cL^1$ we just need to consider monomials
$f_k\,U_k$. In that case the only non-zero matrix elements
correspond to $k=x y^{-1}$. It is enough to show that, for any $f\in
\cS(\A_\K)$, the function $k\mapsto f(k\,b)$ is summable. This
follows from the discreteness of $b\,\K\subset \A_\K$ and the
construction of the Bruhat--Schwartz space $\cS(\A_\K)$, \cf
\cite{Co-zeta}. In fact the associated operator is of finite rank
when $f$ has compact support. In general what happens is that the
sum will look like the sum over $\Z$ of the values $f(nb)$ of a
Schwartz function $f$ on $\R$.
\endproof

 In general  the exact sequence \eqref{akbykstar} does not split
 and one does not have a natural $C_\K$-equivariant trivialization of the
 continuous field $\cH_x$. Thus it is important in the general case to keep
 the nuance between the algebras $C(C_\K,\cL^1(\cH_x))$ and
 $C(C_\K)$. We shall first deal with the special case $\K=\Q$ in
 which this issue does not arise.

 \smallskip

\subsection{The Morita equivalence and cokernel for $\K=\Q$}\label{XKCKsect}

The exact sequence \eqref{akbykstar} splits for $\K=\Q$ and admits a
natural continuous section which corresponds to the open and closed
fundamental domain $\Delta_\Q=\hat\Z^*\times \R^*_+\subset \A_\Q^*$
for the action of $\Q^*$ on ideles. This allows us to construct a
cyclic morphism between the cyclic modules associated, respectively, to the algebra
$C_\rho(\A^*_\Q)\rtimes \Q^*$ and to a suitable algebra
$C_\rho(C_\Q)$ of functions on $C_\Q$.

\begin{lem}\label{equivgrpdsCK}
The composition $d_\Q \circ e_\Q$ of the maps
\begin{equation}\label{eqgrpds}
e_\Q: (k,hb) \mapsto (b,(k,h)), \ \ \text{ and } \ \
d_\Q(k,h)=(kh,h)
\end{equation}
with $b\in \Delta_\Q$ and $k,h\in \Q^*$, gives an isomorphism of the
locally compact groupoids
\begin{equation}\label{eqgrpds2}
\Q^*\ltimes \A_\Q^* \simeq \Delta_\Q \times \Q^* \times \Q^*.
\end{equation}
\end{lem}

\proof The map $e_\Q$ realizes an isomorphism between the locally
compact groupoids
$$  \Q^*\ltimes \A_\Q^* \simeq \Delta_\Q \times (\Q^* \ltimes \Q^*), $$
where $\Q^* \ltimes \Q^*$ is the groupoid of the action of $\Q^*$ on
itself by multiplication. The latter is isomorphic to the trivial
groupoid $\Q^* \times \Q^*$ via the map $d_\Q$.
\endproof

We then have the following result.

\begin{prop}\label{rhoMinftymap}
The map
\begin{equation}\label{Mbmap}
\sum_{k\in \Q^*} f_k\, U_k \mapsto  M_b(x,y)=f_{x y^{-1}}(x\,b),
\end{equation}
for $x,y\in \Q^*$ with $k=x y^{-1}$ and $b\in \Delta_\Q$, defines an
algebra homomorphism
$$ C_\rho(\A_\Q^*)\rtimes \Q^* \to C(\Delta_\Q, M_\infty(\C)) $$
to the algebra of matrix valued functions on $\Delta_\Q$. For any
$f\in \cS(\cG_\Q)$ the element $M_b$ obtained in this way is of
trace class.
\end{prop}

\proof We use the groupoid isomorphism \eqref{eqgrpds} to write $k=x
y^{-1}$ and $k h b= x b$, for $x=kh$ and $y=h$. The second statement
follows from Proposition \ref{restrAKCK1}.
\endproof

Let $\pi=M\circ \rho: \cS(\cG_\K)\to C(\Delta_\Q, M_\infty(\C))$ be the 
composition of the restriction map $\rho$ of
\eqref{rhoAKCKrestr} with the algebra morphism \eqref{Mbmap}. Since
the trace $\Tr$ on $M_\infty(\C)$ gives a cyclic morphism one can
use this to obtain a  morphism of cyclic modules  $(\Tr\circ
\pi)^\natural$, which we now describe explicitly. We let, in the
number field case,
\begin{equation}\label{SCK}
{\bf S\,}(C_\K)=\,\cap_{\beta \in \R}\,\mu^\beta \cS(C_\K),
\end{equation}
where $\mu \in C(C_\K)$ is the module morphism from $C_\K$ to
$\R_+^*$. In the function field case one can simply use for ${\bf
S\,}(C_\K)$ the Schwartz functions with compact support.

\begin{prop}\label{Trrhocyclmap}
The map $\Tr\circ \pi$ defines a morphism $(\Tr\circ \pi)^\natural$
of cyclic modules from $\cS(\cG_\Q)^\natural_0$ to the cyclic
 submodule ${\bf S\,}^\natural(C_\Q)\subset
C(C_\Q)^\natural$ whose elements are continuous functions whose
restriction to the main diagonal belongs to ${\bf S\,}(C_\Q)$.
\end{prop}

\proof By Proposition \ref{rhoMinftymap} the map $\pi$ is an algebra
homomorphism from $\cS(\cG_\Q)$ to $C(\Delta_\Q,\cL^1)\sim
C(C_\Q,\cL^1)$. We need to show that the corresponding cyclic
morphism using $\Tr^\natural$ lands in the cyclic submodule
$S^\natural(C_\Q)$.

For simplicity we can just restrict to the case of monomials, where
we consider elements of the form
\begin{equation}\label{fU0n}
Z=f_{k_0}\,U_{k_0} \otimes f_{k_1}\,U_{k_1} \otimes \cdots \otimes
f_{k_n}\,U_{k_n}.
\end{equation}
The matrix valued functions associated to the monomials
$f_{k_j}\,U_{k_j}$ as in Proposition \ref{rhoMinftymap} have matrix
elements at a point $b\in \Delta_\Q$ that are non zero only for
$x_{j+1}=x_j k_j^{-1}$ and are of the form
\begin{equation}\label{MbfU0n}
f_{k_j}\,U_{k_j}\mapsto M_b(x_j,x_{j+1})= f_{k_j}(x_j b).
\end{equation}
Composing with the cyclic morphism $\Tr^\natural$ gives
\begin{equation}\label{cyclicmapL1}
(\Tr\circ \pi)^\natural(Z)(b_0,b_1,\ldots,b_n)=\sum \prod
M_{b_j}(x_j,x_{j+1})
\end{equation}
where the $x_j\in \K^*$ and $x_{n+1}=x_0$.
 Let $\gamma_0=1$ and
$\gamma_{j+1}=k_j \gamma_j$. Then we find that $(\Tr\circ
\pi)^\natural(Z)=0$, unless $\prod_j k_j=1$, \ie $\gamma_{n+1}=1$.
In this case we obtain
\begin{equation}\label{TrMbfj}
\Tr\circ \pi(Z)(b_0,b_1,\cdots,b_n)=\sum_{k\in \Q^*} \prod_{j=0}^n
f_{k_j} (\gamma_j^{-1} k b_j), \ \ \ \forall b_j \in \Delta_\Q.
\end{equation}

For $n=0$ the formula \eqref{TrMbfj} reduces to
\begin{equation}\label{TrMbf0}
\Tr\circ \pi(f)(b)=\sum_{k\in \Q^*} f (k b), \ \ \ \forall b \in
\Delta_\Q,\ \ \ \forall f\in \cS(\A_\Q)_0,
\end{equation}
where $\cS(\A_\Q)_0 =\Ker \varepsilon_0 \cap \Ker \varepsilon_1
\subset \cS(\A_\Q)$. This gives an element of ${\bf S\,}(C_\Q)$, by
Lemma 2 Appendix 1 of \cite{Co-zeta}. In general, \eqref{TrMbfj}
gives a continuous function of $n+1$ variables on $C_\Q$, and  its
restriction to the main diagonal belongs to ${\bf S\,}(C_\Q)$.
\endproof

Since the category of cyclic modules is an abelian category, we can
consider the cokernel in the category of $\Lambda$-modules of the
cyclic morphism $(\Tr\circ \pi)^\natural$, with $\pi$ the composite
of \eqref{rhoAKCKrestr} and \eqref{Mbmap}. This works nicely for
$\K=\Q$ but makes use of the splitting of the exact sequence
\eqref{akbykstar}.

\medskip

\subsection{The cokernel of $\rho$ for general global fields}\label{H1motAKSect}

To handle the general case in a canonical manner one just needs to
work directly with $C(C_\K,\cL^1(\cH_x))$ instead of
 $C(C_\K)$ and express at that level the decay condition of the restrictions to the diagonal in the
 cyclic submodule ${\bf S\,}^\natural(C_\Q)$ of Proposition \ref{Trrhocyclmap}.

\begin{defn}\label{strongS}
We define ${\bf S\,}^\natural(C_\K, \cL^1(\cH_x))$ to be the cyclic
submodule of the cyclic module $C(C_\K, \cL^1(\cH_x))^\natural$,
whose elements are continuous
functions such that the trace of the restriction to the main diagonal
belongs to ${\bf S\,}(C_\K)$.
\end{defn}

Note that for $T\in C(C_\K, \cL^1(\cH_x))^\natural$ of degree $n$,
$T(x_0,\ldots,x_n)$ is an operator in $\cH_{x_0}\otimes \ldots
\otimes \cH_{x_n}$. On the diagonal, $x_j=x$ for all $j$, the trace
map corresponding to $\Tr^\natural$ is given by
\begin{equation}\label{Trnatural}
\Tr^\natural( T_0\otimes T_1\otimes\ldots\otimes T_n)=\Tr( T_0\,
T_1\,\ldots\, T_n)\,.
\end{equation}
This makes sense since on the diagonal all the Hilbert spaces
$\cH_{x_j}$ are the same.

\smallskip
The argument of Proposition \ref{Trrhocyclmap} extends to the
general case and shows that the cyclic morphism $\rho^\natural$ of
the restriction map $\rho$ lands in ${\bf S\,}^\natural(C_\K,
\cL^1(\cH_x))$.

\begin{defn}\label{motH1AKCK}
We define $\cH^1(\A_\K/\K^*,C_\K)$ to be the cokernel of the cyclic
morphism
$$
\rho^\natural \ :\ \cS(\cG_\K)^\natural_0 \to  {\bf
S\,}^\natural(C_\K, \cL^1(\cH_x))
$$
\end{defn}

Moreover, an important issue arises, since the ranges of continuous
linear maps are not necessarily closed subspaces. In order to
preserve the duality between cyclic homology and cyclic cohomology
we shall define the cokernel of a cyclic map $T: \cA^\natural \to
\cB^\natural$ as the quotient of $\cB^\natural$ by the closure of
the range of $T$. In a dual manner, the kernel of the transposed map
$T^t :  \cB^\sharp \to \cA^\sharp$ is automatically closed and is
the dual of the above.

The choice of the notation $\cH^1(\A_\K/\K^*,C_\K)$ is explained by
the fact that we consider this a first cohomology group, in the
sense that it is a cokernel in a sequence of cyclic homology groups
for the inclusion of the ideles class group in the adeles class
space (dually for the restriction map of algebras), hence we can
think of it as giving rise to an $H^1$ in the relative cohomology
sequence of an inclusion of  $C_\K$ in the noncommutative space
$\A_\K/\K^*$. We can use the result of \cite{CoExt}, describing the
cyclic (co)homology in terms of derived functors in the category of
cylic modules, to write the cyclic homology as
\begin{equation}\label{HCnTorC}
HC_n(\cA) = \Tor_n (\C^\natural,\cA^\natural) .
\end{equation}
Thus, we obtain a cohomological realization of the  cyclic module
$\cH^1(\A_\K/\K^*,C_\K)$ by setting
\begin{equation}\label{TorH1mot}
H^1(\A_\K/\K^*,C_\K):= \Tor(\C^\natural,\cH^1(\A_\K/\K^*,C_\K)).
\end{equation}
We think of this as an $H^1$ because of its role as a relative term
in a cohomology exact sequence of the pair $(\A_\K/\K^*,C_\K)$.

We now show that   $H^1(\A_\K/\K^*,C_\K)$ carries an action of
$C_\K$, which we can view as the abelianization $W_\K^{ab}\sim C_\K$
of the Weil group. This action is induced by the multiplicative
action of $C_\K$ on $\A_\K/\K^*$ and on itself. This generalizes to
global fields the action of $C_\Q=\hat\Z^*\times \R^*_+$ on
$HC_0(D(\cA,\varphi))$ for the Bost--Connes endomotive (\cf
\cite{CCM}).

\begin{prop}\label{CKactionH1}
The cyclic modules $\cS(\cG_\K)^\natural_0$ and ${\bf
S\,}^\natural(C_\K, \cL^1(\cH_x))$ are endowed with an action of
$\A^*_\K$ and the morphism $\rho^\natural$ is $\A^*_\K$-equivariant.
This induces an action of $C_\K$ on $H^1(\A_\K/\K^*,C_\K)$.
\end{prop}

\proof For $\gamma\in \A^*_\K$ one defines an action by
automorphisms of the algebra $\cA=\cS(\cG_\K)$ by setting
\begin{equation}\label{actCKAKbis}
\vartheta_a(\gamma)(f)(x):= f(\gamma^{-1} x), \ \ \text{ for } f\in
\cS(\A_\K),
\end{equation}
\begin{equation}\label{actCKAK}
\vartheta_a(\gamma)(\sum_{k\in \K^*} f_k\,U_k):=\sum_{k\in \K^*}
\vartheta_a(\gamma)(f_k)\,U_k\,.
\end{equation}
This action is inner for $\gamma\in \K^*$ and induces an outer
action
\begin{equation}\label{actCKAKout}
C_\K\to \Out(\cS(\cG_\K))\,.
\end{equation}
 Similarly, the continuous field $\cH_x=\ell^2(p^{-1}(x))$  over
 $C_\K$
 is $\A^*_\K$-equivariant for the action of $\A^*_\K$ on $C_\K$ by
 translations, and the equality
 \begin{equation}\label{actAKltwo}
(V(\gamma)\xi)(y):=\xi(\gamma^{-1}\,y) \qqq y\in p^{-1}(\gamma
x)\,,\ \xi \in \ell^2(p^{-1}(x))\,,
\end{equation}
defines an isomorphism $\cH_x\stackrel{V(\gamma)}\longrightarrow
\cH_{\gamma x}$.
One obtains then an action of $\A^*_\K$ on $C(C_\K, \cL^1(\cH_x))$
by setting
\begin{equation}\label{actCKCK}
\vartheta_m(\gamma)(f)(x):=V(\gamma)\,f(\gamma^{-1}\,x)\,V(\gamma^{-1}),
\ \ \ \forall f\in C(C_\K, \cL^1(\cH_x))\,.
\end{equation}
The morphism $\rho$ is $\A^*_\K$-equivariant, so that one obtains an
induced action on the cokernel $\cH^1(\A_\K/\K^*,C_\K)$. This action
is inner for $\gamma\in \K^*$ and thus induces an action of $C_\K$
on $H^1(\A_\K/\K^*,C_\K)$.
\endproof

We  denote by
\begin{equation}\label{psigactH1}
C_\K \ni \gamma\mapsto {\underline\vartheta_m}(\gamma)
\end{equation}
the induced action on   $H^1(\A_\K/\K^*,C_\K)$.

We have a non-canonical isomorphism
\begin{equation}\label{normrangeCK1}
C_\K \simeq C_{\K,1}\times N ,
\end{equation}
where $N\subset \R_+^*$ is the range of the norm $|\cdot|: C_\K \to
\R^*_+$. For number fields this is $N=\R_+^*$, while for function
fields in positive characteristic $N \simeq \Z$ is the subgroup
$q^{\Z} \subset \R_+^*$ with $q=p^{\ell}$ the cardinality of the field
of constants.
We denote by $\widehat{C_{\K,1}}$ the group of characters of
the compact subgroup $C_{\K,1}\subset C_\K$, \ie the Pontrjagin dual
of $C_{\K,1}$.
Given a character $\chi$ of $C_{\K,1}$, we let $\tilde\chi$ denote the
unique extension of $\chi$ to $C_\K$ which is equal to one on $N$.

One obtains a   decomposition of   $H^1(\A_\K/\K^*,C_\K)$ according
to projectors associated to characters of $C_{\K,1}$.

\begin{prop}\label{charCK1sum}
Characters $\chi\in \widehat{C_{\K,1}}$ determine a canonical direct
sum decomposition
\begin{equation}\label{H1sumchi}
\begin{array}{c}
H^1(\A_\K/\K^*,C_\K) =\bigoplus_{\chi \in \widehat{C_{\K,1}}}
H^1_{\chi}(\A_\K/\K^*,C_\K) \\[3mm]
H^1_{\chi}(\A_\K/\K^*,C_\K)= \{ \xi |\,
{\underline\vartheta_m}(\gamma) \, \xi = \chi(\gamma) \, \xi, \,
\forall \gamma \in C_{\K,1} \}.
\end{array}
\end{equation}
where ${\underline\vartheta_m}(\gamma)$ denotes the induced action
\eqref{psigactH1} on $H^1(\A_\K/\K^*,C_\K)$.
\end{prop}

\proof The action of $\A^*_\K$ on $\cH^1(\A_\K/\K^*,C_\K)$ induces a
corresponding action of $C_\K$ on $H^1(\A_\K/\K^*,C_\K)$. \endproof

We can then reformulate the result of \cite{CCM} based on the trace
formula of \cite{Co-zeta} in the formulation of \cite{Meyer} in
terms of the cohomology $H^1(\A_\K/\K^*,C_\K)$ in the following way.

\begin{prop}\label{specrealH1mot}
The induced representation of $C_\K$ on $H^1_\chi(\A_\K/\K^*,C_\K)$
gives the spectral realization of the zeros of the $L$-function with
Gr\"ossencharakter $\chi$.
\end{prop}

This result is a variant of Corollary 2 of \cite{Co-zeta}, the proof
is similar and essentially reduces to the result of \cite{Meyer}.
There is a crucial difference with \cite{Co-zeta} in that all zeros
(including those not located on the critical line) now appear due to
the choice of the function spaces. To see what happens it is simpler
to deal with the dual spaces \ie to compute the cyclic cohomology
$HC^0$. Its elements are cyclic morphisms $T$ from
$\cH^1(\A_\K/\K^*,C_\K)$ to $\C^\natural$ and they are determined by
the map $T^0$ in degree $0$. The cyclic morphism property then shows
that $T^0$ defines a trace on  ${\bf S\,}^\natural(C_\K,
\cL^1(\cH_x))$ which vanishes on the range of $\rho^\natural$. The
freeness of the action of $\K^*$ on $\A^*_\K$ then ensures that
these traces are given by continuous linear forms on ${\bf
S\,}(C_\K)$ which vanish on the following subspace  of ${\bf
S\,}(C_\K)$ which  is the range of the restriction map, defined as
follows.

\begin{defn}\label{rangeTrpi}
Let $\cV\subset {\bf S\,}(C_\K)$ denote the range of the map
$\Tr\circ \rho$, that is,
\begin{equation}\label{rangecVdef}
\cV=\{ h\in {\bf S\,}(C_\K)|\, h(x)=\sum_{k\in \K^*} \xi(kx), \text{
with } \xi\in \cS(\A_\K)_0\},
\end{equation}
where $\cS(\A_\K)_0=\Ker \varepsilon_0 \cap \Ker\varepsilon_1
\subset \cS(\A_\K)$.
\end{defn}

We have seen above in the case $\K=\Q$ (\cf \cite{Co-zeta}) that the
range of $\Tr\circ \rho$ is indeed contained in ${\bf S\,}(C_\K)$.

Moreover, we have the following results about the action
${\underline\vartheta_m}(\gamma)$, for $\gamma\in C_\K$, on $H^1(\A_\K/\K^*,C_\K)$.
Suppose given  $f \in {\bf S\,} (C_\K)$. We define a corresponding
operator
\begin{equation}\label{varthetafH1}
{\underline\vartheta_m}(f) = \int_{C_\K} f(\gamma)\, {\underline\vartheta_m}(\gamma) \, d^*\gamma,
\end{equation}
acting on the complex vector space $H^1(\A_\K/\K^*,C_\K)$. Here
$d^*\gamma$ is the multiplicative Haar measure on $C_\K$.
 We  have the
following description of the action of ${\underline\vartheta_m}(f)$.

\begin{lem}\label{varthetaH1V}
 For $f\in {\bf S\,}(C_\K)$, the action of the operator
${\underline\vartheta_m}(f)$ of \eqref{varthetafH1} on
$H^1(\A_\K/\K^*,C_\K)$ is the action induced on the quotient of
${\bf S\,}(C_\K)$ by $\cV\subset {\bf S\,}(C_\K)$ of the action of
$\vartheta_m(f)$ on ${\bf S\,}(C_\K)$ by convolution product
\begin{equation}\label{varthetafconvol}
\vartheta_m(f)\xi(u)=\int_{C_\K} \xi(g^{-1} u)  f(g) \, d^*g =( f \star
\xi)(u).
\end{equation}
\end{lem}

\proof One first shows that one can lift $f$ to a function $\tilde
f$ on $\A^*_\K$ such that $$ \sum_{k\in \K^*} \tilde f(kx)=f(x) $$
and that  convolution by $\tilde f$ \ie
$$
\int \tilde f(\gamma)\urep(\gamma)d^*\gamma
$$
leaves $\cS(\cG_\K)$ globally invariant. This means showing that
 that $\cS(\A_\K)_0$ is stable under
convolution by the lift of ${\bf S\,}(C_\K)$. Then
\eqref{varthetafconvol}  follows directly from the definition of the
actions \eqref{actCKCK}, \eqref{actCKAK}, \eqref{psigactH1} and the
operator \eqref{varthetafH1}.
\endproof

For $f \in {\bf S\,}(C_\K)$ and $\tilde\chi$ the extension of a character
$\chi\in \widehat{C_{\K,1}}$ to $C_\K$ and $\hat f (\tilde\chi,\rho)$
the Fourier transform \eqref{rhoFourier}, the operators ${\underline\vartheta_m}(f)$
of \eqref{varthetafH1} satisfy the spectral side of the trace
formula. Namely, we have the following result.

\begin{thm}\label{varthetagTrH1}
For any  $f \in {\bf S\,}(C_\K)$, the operator ${\underline\vartheta_m}(f)$ defined in
\eqref{varthetafH1} acting on $H^1(\A_\K/\K^*,C_\K)$ is of trace
class. The trace is given by
\begin{equation}\label{Trvarthetaf}
\Tr ({\underline\vartheta_m}(f)| H^1(\A_\K/\K^*,C_\K)) = \sum_{\rho\in \C/N^\bot |\,
L\left(\tilde\chi, \rho \right) =0} \hat f (\tilde\chi,\rho),
\end{equation}
with $\hat f (\tilde\chi,\rho)$ the Fourier transform
\eqref{rhoFourier}.
\end{thm}

\proof  Due to the different normalization of the summation map, the
representation ${\underline\vartheta_m}(\gamma)$ considered here
differs from the action $W(\gamma)$ considered in \cite{Co-zeta} by
\begin{equation}\label{Wvarthetag}
{\underline\vartheta_m}(\gamma)=\vert \gamma\vert^{1/2}\,W(\gamma).
\end{equation}
This means that we have
\begin{equation}\label{varthetaWf}
{\underline\vartheta_m}(f) = \int_{C_\K}  f(\gamma)\, {\underline\vartheta_m}(\gamma)  \, d^*\gamma
=\int_{C_\K} h(\gamma)\, W(\gamma) \, d^*\gamma,
\end{equation}
where
\begin{equation}\label{hfrelg}
h(\gamma)=\,\vert \gamma\vert^{1/2}\,f(\gamma).
\end{equation}
We then have, for $W(h)=\int_{C_\K} h(\gamma)\, W(\gamma) \, d^*
\gamma$,
\begin{equation}\label{TrWf}
\Tr W(h) = \sum_{\rho \in \C/N^\bot |\, L(\tilde\chi,\frac{1}{2}+\rho)=0}
\hat h (\tilde\chi,\rho).
\end{equation}
Note that unlike in \cite{Co-zeta} all zeros contribute, including
those that might fail to be on the critical line, and they do with
their natural multiplicity. This follows from the choice of function
space as in \cite{Meyer}. The Fourier transform $\hat
h(\tilde\chi,\rho)$ satisfies
\begin{equation}\label{hathchirho}
\hat h (\tilde\chi ,\rho) = \int_{C_\K} h(u) \tilde\chi(u) \, \vert
u \vert^\rho \, d^*u = \int_{C_\K} f(u) \tilde\chi(u) \,
\vert u \vert^{\rho+1/2} \, d^*u =\hat f (\tilde\chi,\rho+1/2),
\end{equation}
where $h$ and $f$ are related as in \eqref{hfrelg}.
Thus, the shift by $1/2$ in \eqref{TrWf} is absorbed in
\eqref{hathchirho} and this gives the required formula
\eqref{Trvarthetaf}.
\endproof

\subsection{Trace pairing and vanishing}

The commutativity of the convolution product  implies the following
vanishing result.

\begin{lem}\label{vanishV}
Suppose given an element $f\in \cV\subset {\bf S\,}(C_\K)$, where
$\cV$ is the range of the reduction map as in Definition
\ref{rangeTrpi}. Then one has
\begin{equation}\label{restrH1Vzero}
{\underline\vartheta_m}(f)|_{H^1(\A_\K/\K^*,C_\K)} =0.
\end{equation}
\end{lem}

 \proof The result follows by showing that, for $f\in
\cV$, the operator ${\underline\vartheta_m}(f)$ maps any element
$\xi\in {\bf S\,}(C_\K)$ to an element in $\cV$, hence the induced
map on the quotient of ${\bf S\,}(C_\K)$ by $\cV$ is trivial. Since
$\cV$ is a submodule of ${\bf S\,}(C_\K)$ for the action of ${\bf
S\,}(C_\K)$ by convolution we obtain
$$ {\underline\vartheta_m}(f)\xi = f\star\xi = \xi \star f \in \cV, $$
where $\star$ is the convolution product of \eqref{varthetafconvol}.
\endproof

\smallskip

This makes it possible to define a trace pairing as follows.

\begin{rem}\label{tracepairTrH1}
The pairing
\begin{equation}\label{tracepairing}
f_1 \otimes f_2 \mapsto \langle f_1, f_2 \rangle_{H^1} :=
\Tr({\underline\vartheta_m}(f_1\star f_2) | H^1(\A_\K/\K^*,C_\K) )
\end{equation}
descends to a well defined pairing on $H^1(\A_\K/\K^*,C_\K)\otimes
H^1(\A_\K/\K^*,C_\K)$.
\end{rem}

\bigskip

\section{Primitive cohomology}\label{motivicSect}

The aim of this section is to interpret the motivic
construction described in the previous paragraph as the
noncommutative version of a classical construction in algebraic
geometry. In motive theory, realizations of (mixed) motives appear
frequently in the form of kernels/cokernels of relevant
homomorphisms. The {\it primitive cohomology} is the example we
shall review hereafter.\vspace{.1in}

If $Y$ is a compact K\"ahler variety, a K\"ahler cocycle class
$[\omega]\in H^2(Y,\bR)$ determines the Lefschetz operator ($i\in
\bZ_{\ge 0}$):
\[
L: H^i(Y,\bR) \to H^{i+2}(Y,\bR),\quad L(a):= [\omega]\cup a.
\]
Let $n = \dim Y$. Then, the {\it primitive cohomology} is defined as
the kernel of iterated powers of the Lefschetz operator
\[
H^i(Y,\bR)_{prim} := {\rm Ker}(L^{n-i+1}: H^i(Y,\bR) \to
H^{2n-i+2}(Y,\bR)).
\]
In particular, for $i=n$ we have
\[
H^n(Y,\bR)_{prim} := {\rm Ker}(L: H^n(Y,\bR) \to H^{n+2}(Y,\bR)).
\]
Let assume, from now on, that $j: Y \hookrightarrow X$ is a smooth
hyperplane section of a smooth, projective complex algebraic variety
$X$. Then, it is a classical result of geometric topology that $L =
j^*\circ j_*$, where
\[
j_*: H^i(Y,\bR) \to H^{i+2}(X,\bR)
\]
is the Gysin homomorphism: the Poincar\'e dual of the restriction
homomorphism
\[
j^*: H^{2n-i}(X,\bR) \to H^{2n-i}(Y,\bR).
\]
In fact, because the class of $L$ comes from an integral class, the
equality $L = j^*\circ j_*$ holds already in integral cohomology.
For $i=n$, the above description of the Lefschetz operator together
with the Lefschetz theorem of hyperplane sections imply that
\[
H^n(Y,\bR)_{prim} \cong {\rm Ker}(j_*: H^n(Y,\bR) \to
H^{n+2}(X,\bR))=: H^n(Y,\bR)_{van}
\]
where by $H^i(Y,\bR)_{van}$ we denote the {\it vanishing cohomology}
\[
H^i(Y,\bR)_{van}:={\rm Ker}(j_*: H^i(Y,\bR) \to H^{i+2}(X,\bR)).
\]
Now, we introduce the theory of mixed Hodge structures in this
set-up.

Let $U := X\smallsetminus Y$ be the open space which is the complement of
$Y$ in $X$ and let denote by $k: U\hookrightarrow X$ the
corresponding open immersion. Then, one knows that $R^ij_*\bZ = 0$
unless $i=0,1$ so that the Leray spectral sequence for $j$:
\[
E_2^{p,q} = H^q(X,{\text R}^pk_*\bZ) \Rightarrow H^{p+q}(U,\bZ)
\]
coincides with the long exact sequence (of mixed Hodge structures)
\[
\ldots\stackrel{\partial}{\to} H^{i-2}(Y,\bZ)(-1)
\stackrel{j_*}{\to} H^i(X,\bZ) \to H^i(U,\bZ)
\stackrel{\partial}{\to}\ldots
\]
The boundary homomorphism $\partial$ in this sequence is known to
coincide (\cite{De}, \S~9.2) with the residue homomorphism
\[
Res: H^{i+1}(U,\bZ) \to H^i(Y,\bZ)(-1)
\]
whose description, with complex coefficients, is derived from a
corresponding morphism of filtered complexes (Poincar\'e residue
map). This morphism fits in the following exact sequence of filtered
complexes of Hodge modules
\[
0 \to \Omega_X^\cdot \to \Omega_X^\cdot(\log Y) \stackrel{{\rm
res}}{\to}j_*\Omega_Y^\cdot[-1] \to 0
\]
\[
{\rm res}(\alpha\wedge\frac{dt}{t}) = \alpha_{|Y}.
\]
One knows that Res is a homomorphism of Hodge structures, hence the
Hodge filtration on $H^{n+1}(U,\bC)\cong \mathbb
H^{n+1}(X,\Omega_X^\cdot({\rm log}Y))$ determines a corresponding
filtration on the (twisted) vanishing cohomology
\[
H^n(Y,\bC)(n)_{van} = {\rm Ker}(j_*: H^n(Y,\bC)(n) \to
H^{n+2}(X,\bC)(n+1))\cong H^{n+1}(U)(n+1).
\]

In degree $i=n$, one also knows that the excision exact sequence (of
Hodge structures) becomes the short exact sequence
\[
0\to H^n(X,\bC)\stackrel{j^*}{\to}H^n(Y,\bC)\to H^{n+1}_c(U,\bC)\to
0.
\]
Therefore, it follows by the Poincar\'e duality isomorphism
\[
H^{n+1}_c(U,\bC)^{*} \cong H^{n+1}(U,\bC)(n+1)
\]
that
\begin{equation}\label{mainiso}
({\rm Coker}(j^*: H^n(X,\bC)\to H^n(Y,\bC)))^* \cong
H^{n+1}_c(U,\bC)^{*} \cong H^n(Y,\bC)(n)_{van}.
\end{equation}\vspace{.1in}

When $j: Y\hookrightarrow X$ is a singular hypersurface or a divisor
in $X$ with (local) normal crossings ({\it i.e.}: $Y = \bigcup_i
Y_i$, $\dim Y_i = n = \dim X -1$, $Y$ locally described by an
equation $x_{i_1}\cdots x_{i_r} = 0$, $\{i_1,\ldots
i_r\}\subseteq\{1,\ldots n+1\}$, $\{x_1,\ldots x_{n+1}\}=$ system of
local coordinates in $X$), the notion of the Gysin homomorphism is
lost. One then replaces the vanishing cohomology by the primitive
cohomology, whose definition extends to this general set-up and is
given, in analogy to \eqref{mainiso}, as
\[
H^n(Y,\bC)_{prim} := {\rm Coker}(j^*: H^n(X,\bC) \to
H^n(Y,\bC))\subseteq H^{n+1}_c(U,\bC).
\]

One also knows that the primitive cohomology is motivic (\cf
\cite{GM} and \cite{BEK} for interesting examples). Following the
classical construction that we have just reviewed, we like to argue
now that the definition of the cyclic module
$\cH^1(\A_\K/\K^*,C_\K)$ (as in Definition~\ref{motH1AKCK}), which is based
on a noncommutative version of a restriction map ``from adeles to
ideles'' defined in the category of $\Lambda$-modules, should be interpreted
as the noncommutative analogue of a primitive motive (a cyclic
primitive module). The cohomological realization of such motive
({\it i.e.} its cyclic homology) is given by the group
$H^1(\A_\K/\K^*,C_\K) = \Tor(\C^\natural,\cH^1(\A_\K/\K^*,C_\K))$ (\cf
\eqref{TorH1mot}) which therefore can be interpreted as a
noncommutative version of a primitive cohomology.

\section{A cohomological Lefschetz trace formula}\label{cohTrSect}

\subsection{Weil's explicit formula as a trace formula}

As in \S \ref{expformulaSect} above, let $\alpha$ be a non-trivial
character of $\A_\K$ which is trivial on $\K\subset \A_\K$.
It is well known (\cite{Weil} VII-2) that for such a character
$\alpha$ there exists a differental idele $a=(a_v)\in \A_\K^*$
such that
\begin{equation}\label{diffidele}
\alpha_v(x)= e_{\K_v}(a_v\,x), \ \ \ \forall  x\in \K_v ,
\end{equation}
where, for a local field $K$, the additive character $e_K$ is
chosen in the following way.
\begin{itemize}
\item If $K=\R$ then $e_\R(x)=e^{-2\pi ix}$, for all $x\in \R$.
\item If $K=\C$ then $e_\C(z)=e^{-2\pi i(z+\bar z)}$, for all $z\in \C$.
\item If $K$ is a non-archimedean local field with maximal compact
subring $\cO$, then the character $e_K$ satisfies $\Ker\, e_K=\cO$.
\end{itemize}

The notion of differental idele can be thought of as an extension of
the canonical class of the algebraic curve $C$, from the setting of
function fields $\F_q(C)$ to arbitrary global fields $\K$. For instance,
one has
\begin{equation}\label{eulerdiffidele}
|a|=q^{2-2g} \ \ \ \ \text{ or } \ \ \ \  |a|=D^{-1},
\end{equation}
respectively, for the case of a function field $\F_q(C)$ and of a number field.
In the number field case $D$ denotes the discriminant.

\medskip

In \cite{CCM} we gave a cohomological formulation of the Lefschetz
trace formula of \cite{Co-zeta}, using the version of the
Riemann--Weil explicit formula as a trace formula given in
\cite{Meyer} in the context of nuclear spaces, rather than the
semi-local Hilbert space version of \cite{Co-zeta}.

\begin{thm}\label{TraceformulaThm}
For $f\in {\bf S\,}(C_\K)$ let ${\underline\vartheta_m}(f)$ be the operator
\eqref{varthetafH1} acting on the space $H^1=H^1(\A_\K/\K^*,C_\K)$.
Then the trace is given by
\begin{equation}\label{TraceformulaH1}
\Tr({\underline\vartheta_m}(f) |H^1)=  \hat f(0) + \hat f(1) - (\log |a|) \, f(1)
-\sum_{v\in \Sigma_\K} \int_{(\K^*_v,e_{\K_v})}' \frac{f(u^{-1})}{
 |1-u|}\, d^* u .
\end{equation}
\end{thm}

The formula \eqref{TraceformulaH1} is obtained in \cite{CCM} first by
showing that the Lefschetz trace formula of \cite{Co-zeta} in the
version of \cite{Meyer} can be formulated equivalently in the form
\begin{equation}\label{Geomtrace1}
\Tr({\underline\vartheta_m}(f)|H^1)=\hat f(0) + \hat f(1) -\sum_{v\in \K_v}
\int_{\K^*_v}' \frac{f(u^{-1})}{ |1-u|}\, d^* u ,
\end{equation}
where one uses the global character $\alpha$ to fix the local
normalizations of the principal values in the last term of the
formula. We then compute this principal value using the differental
idele in the form
\begin{equation}\label{intalphae}
\int_{(\K^*_v,\alpha_v)}' \frac{f(u^{-1})}{
 |1-u|}\, d^* u=(\log |a_v|)\,f(1)+ \int_{(\K^*_v,e_{\K_v})}'  \frac{f(u^{-1})}{
 |1-u|}\, d^* u.
\end{equation}

\subsection{Weil Positivity and the Riemann
Hypothesis}\label{posTrH1motSect}

We introduce an involution for elements
$f\in {\bf S\,} (C_\K)$ by setting
\begin{equation}\label{fstarg}
f^*(g)=\overline{ f(g^{-1}) }.
\end{equation}
We also consider a one parameter group $z\mapsto \Delta^z$
of automorphisms of the convolution algebra ${\bf S\,} (C_\K)$, with the
convolution product \eqref{varthetafconvol} by setting
\begin{equation}\label{Deltazfg}
\Delta^z(f)(g)=| g |^z\,f(g),
\end{equation}
for $f\in {\bf S\,} (C_\K)$ and $z\in \C$. Since \eqref{Deltazfg} is given
by multiplication by a character, it satisfies
\begin{equation}\label{Deltazfstarh}
\Delta^z(f\star h)=\Delta^z(f)\star \Delta^z(f), \ \ \ \forall f,h \in {\bf S\,}
(C_\K).
\end{equation}
We consider also the involution
\begin{equation}\label{fsharpinv}
f\mapsto f^{\sharp}=\Delta^{-1}\, f^*, \ \ \ \text{ with } \ \
f^{\sharp}(g)=|g|^{-1}
\overline{ f(g^{-1}) }.
\end{equation}

\smallskip

The reformulation, originally due to A. Weil, of the Riemann
Hypothesis  in our setting is given by the following statement.

\begin{prop} \label{RHpositTr}
The following two conditions are equivalent:
\begin{itemize}
\item All $L$-functions with Gr\"ossencharakter on $\K$ satisfy the Riemann
Hypothesis.
\item The trace pairing \eqref{tracepairing}
satisfies the positivity condition
\begin{equation}\label{TrpositivityDelta}
\langle \Delta^{-1/2}\, f, \Delta^{-1/2}\,f^*\rangle \geq 0, \ \ \ \
\forall f\in{\bf S\,}(C_\K).
\end{equation}
\end{itemize}
\end{prop}

\proof Let $W(\gamma)=| \gamma
|^{-1/2}\,{\underline\vartheta_m}(\gamma)$. Then, by \cite{weilpos}
the RH for $L$-functions with Gr\"ossencharakter on $\K$ is
equivalent to the positivity
\begin{equation}\label{TrWfstarf}
 \Tr(W(f\star f^*))\geq 0, \ \ \ \forall f\in {\bf S\,}(C_\K).
\end{equation}
Thus, in terms of the representation ${\underline\vartheta_m}$ we
are considering here, we have
$$W(f)={\underline\vartheta_m}(\Delta^{-1/2}\, f).$$ Using the
multiplicative property \eqref{Deltazfstarh} of $\Delta^z$ we
rewrite \eqref{TrWfstarf} in the equivalent form
\eqref{TrpositivityDelta}.
\endproof

\smallskip

In terms of the involution \eqref{fsharpinv} we can reformulate Proposition
\ref{RHpositTr} in the following equivalent way.

\begin{cor}\label{RHpositTrsharp}
The following conditions are equivalent
\begin{itemize}
\item All $L$-functions with Gr\"ossencharakter on $\K$ satisfy the Riemann
Hypothesis.
\item The trace pairing \eqref{tracepairing} satisfies $\langle  f,
f^\sharp\rangle \geq 0$, for all $f\in{\bf S\,}(C_\K)$.
\end{itemize}
\end{cor}

\proof In \eqref{TrpositivityDelta} we write
$\Delta^{-1/2} f=h$. This gives $$\Delta^{-1/2} f^*=\Delta^{-1/2}
(\Delta^{1/2} h)^*=\Delta^{-1} h^*=h^\sharp$$ and the result
follows, since $\Delta^{-1/2}$ is an automorphism of ${\bf S\,}(C_\K)$.
\endproof

\smallskip

The vanishing result of Lemma \ref{vanishV}, for elements in the range
$\cV\subset {\bf S\,}(C_\K)$ of the reduction map $\Tr\circ \rho$ from
adeles, gives then the following result.

\begin{prop}\label{pairTrvanishV}
The elements $f\star f^\sharp$ considered in Corollary
\ref{RHpositTrsharp} above have the following properties.
\begin{enumerate}
\item The trace pairing $\langle f,f^\sharp\rangle$ vanishes for all
$f \in \cV$, \ie when $f$ is the restriction $\Tr\circ \rho$ of an
element of $\cS(\cG_\K)$.
\item By adding elements of $\cV$ one can make the values
\begin{equation}\label{Vepsilonfstarf}
f\star f^\sharp(1)=\int_{C_\K} |f(g)|^2\,|g|\,d^*g < \epsilon
\end{equation}
for arbitrarily small $\epsilon>0$.
\end{enumerate}
\end{prop}

\proof (1) The vanishing result of Lemma \ref{vanishV} shows that
${\underline\vartheta_m}(f)|_{H^1(\A_\K/\K^*,C_\K)}=0$ for all $f\in \cV$.
Thus, the
trace pairing satisfies $\langle f,h\rangle=0$, for $f\in \cV$ and
for all $h\in {\bf S\,}(C_\K)$. In particular this applies to the case
$h=f^\sharp$.

\smallskip

(2) This follows from the surjectivity of the map $\fE$
 for the weight $\delta=0$ (\cf Appendix 1 of \cite{Co-zeta}).
\endproof

\smallskip

Proposition \ref{pairTrvanishV} shows that the trace pairing admits a
large radical given by all functions that extend to adeles. Thus, one
can divide out this radical and work with the
cohomology $H^1(\A_\K/\K^*,C_\K)$ described above.

\section{Correspondences}

To start building the dictionary between the Weil proof and the
noncommutative geometry of the adeles class space, we begin by
reformulating the trace formula discussed above in more intersection theoretic
language, so as to be able to compare it with the setup of \S
\ref{CorrdivSect} above. We also discuss in this section the analog of
moding out by trivial correspondence.

\subsection{The scaling correspondence as Frobenius}\label{ScaleCorrSect}

To the scaling action $$\vartheta_a(\gamma)(\xi)(x)=\xi(\gamma^{-1}x)$$ one
associates the graph $Z_g$ given by the pairs $(x,g^{-1}x)$. These
should be considered as points in the product $\A_\K/\K^*\times
\A_\K/\K^*$ of two copies of the adeles class space. Thus, the
analog in our context of the correspondences $Z=\sum_n a_n \Fr^n$ on
$C\times C$ is given by elements of the form
\begin{equation}\label{Zfcorr}
Z(f)=\int_{C_\K} f(g)\, Z_g \, d^*g,
\end{equation}
for some $f\in {\bf S\,}(C_\K)$.

\smallskip

With this interpretation of correspondences, we can then make sense of
the terms in the trace formula in the following way.

\begin{defn}\label{degcodegZf}
For a correspondence of the form \eqref{Zfcorr} we define degree and
codegree by the following prescription
\begin{equation}\label{degcorrZf}
d(Z(f)):= \hat f(1) = \int_{C_\K} f(u)\, |u|\, d^*u,
\end{equation}
\begin{equation}\label{codegcorrZf}
d'(Z(f)):= d(Z(f^\sharp))= \int_{C_\K} f(u)\, d^*u =\hat f(0).
\end{equation}
\end{defn}

Here the Fourier transform $\hat f$ is as in \eqref{rhoFourier}, with
the trivial character $\chi=1$.
Notice that, with this definition of degree and codegree we find
\begin{equation}\label{degcodegZg}
d(Z_g)=|g|, \ \ \ \ \text{ and } \ \ \ \ d'(Z_g)=1.
\end{equation}

\smallskip

Thus, the term $\hat f(1) + \hat f(0)$ in the trace formula of Theorem
\ref{TraceformulaThm} match the term $d(Z)+d'(Z)$ in Weil's
formula for the trace of a correspondence as in \eqref{TraceCorr}.
The term
\begin{equation}\label{intKvalphav2}
 -\int_{(\K^*_v,\alpha_v)}' \frac{f(u^{-1})}{
|1-u|}\, d^* u
\end{equation}
of \eqref{Geomtrace1} in turn can be seen as the remaining term
$-Z\bullet \Delta$ in \eqref{TraceCorr}. In fact, the formula
\eqref{intKvalphav2} describes, using distributions, the local
contributions to the trace of the intersections between the graph
$Z(f)$ and the diagonal $\Delta$. This was proved in \cite{Co-zeta},
Section VI and Appendix III. It generalizes the analogous formula for
flows on manifolds of \cite{GS}, which in turn can be seen as a
generalization of the usual Atiyah--Bott Lefschetz formula for a
diffeomorphism of a smooth compact manifold \cite{AB}.

\smallskip

When we separate out the contribution $\log|a|\, h(1)$, as in passing
from \eqref{Geomtrace1} to \eqref{TraceformulaH1},
and we rewrite the trace formula as in Theorem
\ref{TraceformulaThm}. This corresponds to separating the intersection
$Z\bullet \Delta$ into a term that is proportional to the self
intersection $\Delta \bullet \Delta$ and a remaning terms where the
intersection is transverse.

To see this, we notice that the term $\log|a|$, for $a=(a_v)$ a
differental idele, is of the form \eqref{eulerdiffidele}. Indeed one
sees that, in the function field case the term
$$ -\log|a| = -\log q^{2-2g} =(2g-2)\log q = - \Delta\bullet\Delta \,
\log q $$
is proportional to the self intersection of the diagonal, which brings
us to consider the value $\log|a|=-\log D$ with the discriminant of a
number field as the analog in characteristic zero of the self
intersection of the diagonal.

\smallskip

In these intersection theoretic terms we can reformulate the
positivity condition (\cf \cite{EB}) equivalent to the Riemann
Hypothesis in the following way.

\begin{prop}\label{PosTrRHint}
The following two conditions are equivalent
\begin{itemize}
\item All $L$-functions with Gr\"ossencharakter on $\K$ satisfy the
Riemann Hypothesis.
\item The estimate
\begin{equation}\label{ZfbulletZf}
Z(f)\bullet_{trans} Z(f) \leq  2 d(Z(f))d'(Z(f)) - \Delta\bullet\Delta \,
f\star f^\sharp (1)
\end{equation}
holds for all $f\in {\bf S\,}(C_\K)$.
\end{itemize}
\end{prop}

\proof As in the Weil proof one separates the terms $Z\star
Z'=d'(Z)\Delta + Y$, where $Y$ has transverse intersection with the
diagonal, here we can write an identity
\begin{equation}\label{ZfZftransv}
\Tr({\underline\vartheta_m}(f\star f^\sharp)|H^1)=: Z(f)\bullet Z(f) =
\Delta\bullet\Delta \, f\star f^\sharp (1) + Z(f)\bullet_{trans}
Z(f)
\end{equation}
where the remaning term $Z(f)\bullet_{trans} Z(f)$ which represents
the transverse intersection is given by the local contributions
given by the principal values over $(\K_v^*,e_{\K_v})$ in the
formula \eqref{TraceformulaH1} for $\Tr({\underline\vartheta_m}(f\star
f^\sharp)|H^1)$.

\smallskip

The formula \eqref{TraceformulaH1} for $\Tr({\underline\vartheta_m}(f\star
f^\sharp)|H^1)$ gives a term of the form $-\log|a|\, f\star f^\sharp
(1)$, with
$$ f\star f^\sharp (1) =\int_{C_\K} |f(g)|^2\, |g|\, d^*g. $$
We rewrite this term as $- \Delta\bullet \Delta \, f\star f^\sharp
(1)$ according to our interpretation of $\log|a|$ as self-intersection
of the diagonal. This matches the term $(2g-2) d'(Z)$ in the estimate
for $\Tr(Z\star Z')$ in the Weil proof.

\smallskip

The first two terms in the formula \eqref{TraceformulaH1} for
$\Tr({\underline\vartheta_m}(f\star f^\sharp)|H^1)$ are of the form
\begin{equation}\label{terms01Tr}
\widehat{ f\star f^\sharp}(0) + \widehat{ f\star f^\sharp}(1) =
2\hat f(0)\hat f(1) = 2  d'(Z(f))d(Z(f)).
\end{equation}
This matches the term $2g d'(Z)=2d(Z)d'(Z)$ in the expression for
$\Tr(Z\star Z')$ in the Weil proof.

\smallskip

With this notation understood, we see that the positivity
$\Tr({\underline\vartheta_m}(f\star f^\sharp)|H^1)\geq 0$ corresponds indeed to the
estimate \eqref{ZfbulletZf}.
\endproof

\subsection{Fubini's theorem and the trivial correspondences}

As we have seen in recalling the main steps in the Weil proof,
a first step in dealing with correspondences is to use the freedom to
add multiples of the trivial correspondences in order to adjust the
degree. We describe an analog, in our noncommutative geometry setting, of
the trivial correspondences and of this operation of modifying the
degree.

\smallskip

In view of the result of Proposition \ref{pairTrvanishV} above, it is
natural to regard the elements $f\in \cV\subset {\bf S\,}(C_\K)$ as those
that give rise to the trivial correspondences $Z(f)$. Here, as above,
$\cV$ is the range of the reduction map from adeles.

\smallskip

The fact that it is possible to arbitrarily modify the degree
$d(Z(f))=\hat f(1)$ of a correspondence by adding to $f$ an element in
$\cV$ depends on the subtle fact that we deal with a case where the
Fubini theorem does not apply.

\smallskip

In fact, consider an element $\xi\in \cS(\A_\K)_0$. We know that it
satisfies the vanishing condition
$$ \int_{\A_\K} \xi(x)\, dx =0. $$
Thus, at first sight it would appear that, for the function on
$C_\K$ defined by $f(x)=\sum_{k\in \K^*} \xi(kx)$,
\begin{equation}\label{fubini1}
\hat f(1) = \int_{C_\K} f(g) |g| d^*g
\end{equation}
should also vanish, since we have $f(x)=\sum_{k\in \K^*} \xi(kx)$
and for local fields (but not in the global case) the relation
between the additive and multiplicative Haar measures is of the form
$dg=|g| d^*g $. This, however, is in general not the case. To see
more clearly what happens, let us just restrict to the case $\K=\Q$
and assume that the function $\xi(x)$ is of the form
$$ \xi = {\bf 1}_{\hat\Z}\otimes \eta, $$
with ${\bf 1}_{\hat\Z}$ the characteristic function of $\hat\Z$ and
with $\eta\in \cS(\R)_0$. We then have $C_\Q={\hat \Z}^*\times
\R_+^*$ and the function $f$ is of the form
\begin{equation}\label{fubini2}
f(u,\lambda)=\sum_{n\in \Z,n\neq 0} \eta(n\lambda), \ \ \ \forall
\lambda\in \R^*_+\,,\ u\in {\hat \Z}^*\,.
\end{equation}
We can thus write \eqref{fubini1} in this case as
\begin{equation}\label{fubini3}
 \hat f(1)=\int_{{\hat \Z}^*\times\R^*_+} f(u,\lambda) du\,d\lambda  =
 \int_\R \;\sum_{n \in \N}
\eta(n\lambda)\; d\lambda
\end{equation}
Moreover since $\eta\in \cS(\R)_0$ we have for all $n$,
\begin{equation}\label{fubini4}
 \int_\R
\eta(n\lambda) d\lambda=0\,.
\end{equation}
 It is however {\em
not} necessarily the case that we can apply Fubini's theorem and
write
\begin{equation}\label{nothappen}
 \int_\R \;\sum_{n \in \N}
\eta(n\lambda)\; d\lambda=  \sum_n \int_\R \eta(n\lambda) d\lambda
=0
\end{equation}
since as soon as $\eta \neq 0$ one has
$$ \sum_{n=1}^\infty \int_\R |\eta(n\lambda)| d\lambda =
(\int_\R |\eta(\lambda)| d\lambda)\, \ \sum_{n=1}^\infty \frac
1n\,=\infty
$$ so that Fubini's theorem does not apply and one cannot
interchange the integral and the sum in \eqref{nothappen}. Thus, one
can in general have $\hat f(1)\neq 0$, even though $\sum_n \int_\R
\eta(n\lambda) d\lambda =0$. In fact, we have the following result.

\begin{lem}\label{degmodifyV}
Given $f\in {\bf S\,}(C_\K)$, it is possible to change arbitrarily the value of the
degree $d(Z(f))=\hat f(1)$ by adding elements of $\cV$.
\end{lem}

\proof It suffices to exhibit an element $f\in \cV$ such that $\hat
f(1)\neq 0$, as then by linearity one can obtain the result. We only
treat the case $\K=\Q$.  We take $\eta\in \cS(\R)_0$ given by
$$
\eta(x)= \pi x^2(\pi x^2-\frac 32)\, e^{-\pi x^2}
$$
One finds   that, up to normalization, the Fourier transform $\hat
f$ is given by
$$
\hat f(is)=\int_{\R^*_+}\sum_{n \in \N}
\eta(n\lambda)\,\lambda^{is}d^*\lambda=s(s+i)\zeta^*(is)
$$
where $\zeta^*$ is the complete zeta function,
\begin{equation}\label{xi0}
\zeta^* (z) = \pi^{-z/2} \, \Gamma \left( \frac{z}{2} \right) \zeta
(z).
\end{equation}

This function has a simple pole at $z=1$ thus one gets that $\hat
f(1)\neq 0$.
\endproof

\medskip

An important question, in order to proceed and build a dictionary
that parallels the main steps in the Weil proof, is to identify the
correct notion of principal divisors. To this purpose, we show that
we have at least a good analog for the points of the curve, in terms
of states of some thermodynamical systems, that extend from the
function field setting to the number field case.

\section{Thermodynamics and geometry of the primes}

Let $\K$ be a global field, with $\A_\K$ the ring of adeles and $C_\K$
the ideles classes, as above. We denote by $C_{\K,1}\subset C_\K$ the
kernel of the norm $|\cdot|: C_\K \to \R^*_+$.

The origin (\cf \cite{Co-zeta}) of the terms in the geometric side
of the trace formula (Theorem \ref{TraceformulaThm}) comes from the
Lefschetz formula by Atiyah-Bott \cite{AB} and its adaptation by
Guillemin-Sternberg (\cf \cite{GS}) to the distribution theoretic
trace for flows on manifolds, which is a variation on the theme of
\cite{AB}. For the action of $C_\K$ on the adele class space $X_\K$
the relevant periodic points are
\begin{equation} \label{globalper}
P=\{ (x,u)\in  X_\K \times C_\K \,|\,  u\,x=x\}
\end{equation}
and one has (\cf \cite{Co-zeta})
\begin{prop}\label{preorbitckak} Let $(x,u)\in P$, with $u\neq 1$.
There exists a place $v\in \Sigma_\K$ such that
\begin{equation} \label{globalper1}
x\in X_{\K,v}=\{ x\in  X_\K \,|\,  x_v=0\}
\end{equation}
The isotropy subgroup of any $x\in X_{\K,v}$ contains the cocompact
subgroup
\begin{equation} \label{globalper2}
\K^*_v\subset C_\K\,,\  \  \ \K^*_v=\{ (k_w)\,|\, k_w=1\ \forall
w\neq v\}
\end{equation}
\end{prop}
The spaces $X_{\K,v}$ are noncommutative spaces, as such they are
described by the following noncommutative algebras:

\begin{defn}\label{AKvdefn}
Let $\A_{\K,v} \subset \A_\K$ denote the closed $\K^*$-invariant
subset of adeles
\begin{equation}
\A_{\K,v}=\{a=(a_w)_{w\in \Sigma_\K}\,|\, a_v=0\}.
\end{equation}
Let $\cG_{\K,v}$ denote the closed subgroupoid of $\cG_\K$ given by
\begin{equation}\label{cRKv}
\cG_{\K,v} =\{ (k,x)\in \cG_\K \,|\, x_v=0 \},
\end{equation}
and let $\cA_v=\cS(\cG_{\K,v})$ be the corresponding groupoid
algebra.
\end{defn}

Since the inclusion $\A_{\K,v} \subset \A_\K$ is  $\K^*$-equivariant
and proper, it extends to an algebra homomorphism
\begin{equation}\label{rhovrestr}
\rho_v\;:\; \cS(\cG_{\K})\to \cS(\cG_{\K,v})
\end{equation}
which plays the role of the restriction map to the periodic orbit
$X_{\K,v}$. We shall now determine the classical points of each of
the $X_{\K,v}$. Taken together these will form the following locus
inside the adeles class space, which  we refer to as the ``periodic
classical points'' of $X_\K=\A_\K/\K^*$.

\smallskip

\begin{defn}\label{classXiadeles}
Let $\K$ be a global field. For a place $v\in \Sigma_\K$ consider
the adele
\begin{equation}\label{avwadele}
a^{(v)}=(a^{(v)}_w), \ \ \ \text{ with } \ \ a^{(v)}_w =\left\{
\begin{array}{lr} 1 & w\neq v \\ 0 & w=v. \end{array}\right.
\end{equation}
The set of {\em periodic classical} points of the adeles class space
$\A_\K/\K^*$ is defined as the union of orbits
\begin{equation}\label{Xiorbits}
\Xi_\K := \bigcup_{v\in \Sigma_\K} C_\K a^{(v)}.
\end{equation}
\end{defn}

\smallskip

\subsection{The global Morita equivalence}\label{globMor}

In order to deal with states rather than weights, we perform a
global Morita equivalence, obtained by reducing the groupoid
$\cG_\K$ by a suitable open set. The set $\A^{(1)}_\K$ of
\eqref{cWQneighbglob} that we use
to reduce the groupoid $\cG_\K$ will only capture part of the
classical subspace $C_\K$, but since our main focus is on the
geometry of the complement of this subspace (the cokernel of the
reduction map), this will not be a problem.

\begin{lem}\label{neighb0ideles}
Let $\K$ be a global field. Let $W\subset \A_\K$ be a neighborhood
of $0\in \A_\K$. Then for $x\in  \A_\K$ one has $\K^* x\cap W\neq
\emptyset$, unless $x\in \A_\K^*$ is an idele. For $x\in \A_\K^*$,
the orbit $\K^* x$ is discrete in $\A_\K$.
\end{lem}

\proof One can assume that $W$ is of the form
$$
W=\{a=(a_w)|\ |a_w|<\varepsilon \ \forall w\in S \text{ and }
|a_w|\leq 1\ \forall w\notin S\},
$$
for $S$ a finite set of places and for some $\varepsilon>0$.
Multiplying by a suitable idele one can in fact assume that
$S=\emptyset$, so that we have
$$
W=\{a=(a_w) |\ |a_w|\leq 1\ \forall w\in \Sigma_\K\}.
$$
One has $|x_v|\leq 1$ except on a finite set $F\subset \Sigma_\K$ of
places. Moreover, if $x$ is not an idele, one can also assume that
$$
\prod_{v\in F} |x_v|<\delta
$$
for any fixed $\delta$. Thus, $- \log |x_v|$ is as large as one
wants and there exists $k\in \K^*$ such that $k\,x\in W$. This is
clear in the function field case because of the Riemann Roch formula
\eqref{CRRthm}. In the case of $\Q$ one can first multiply $x$ by an
integer to get $|x_v|\leq 1$ for all finite places, then since this
does not alter the product of all $|x_v|$ one gets $|x_\infty|<1$
and $x\in W$. The case of more general number fields is analogous.
In the case of ideles, one can assume that $x=1$ and then the second
statement follows from the discreteness of $\K$ in $\A_\K$.
\endproof

We consider the following choice of a neighborhood of zero.

\begin{defn}\label{glMorgrpd}
 Consider the open neighborhood of
$0\in \A_\K$ defined by
\begin{equation}\label{cWQneighbglob}
\A^{(1)}_\K=\prod_{w\in \Sigma_\K} \K_w^{(1)}\subset \A_\K
\end{equation}
where for any place we let $\K_w^{(1)}$ be the {\em interior} of
$\{x\in \K_w\,;\, |x|\leq 1\}$. Let $\cG_\K^{(1)}$ denote the
reduction of the groupoid $\cG_\K$ by the open subset $\A^{(1)}_\K\subset
\A_\K$ of the units and let $\cS(\cG_\K^{(1)})$
denote the corresponding (smooth) groupoid algebra.
\end{defn}

The algebra $\cS(\cG_\K^{(1)})$ is a subalgebra of $\cS(\cG_\K)$ where one
simply extends the function $f(k,x)$ by zero outside of the open
subgroupoid $\cG_\K^{(1)}\subset \cG_\K$. With this convention, the
convolution product of $\cS(\cG_\K^{(1)})$ is simply given by the
convolution product of $\cS(\cG_\K)$ of the form
$$
(f_1\star f_2)(k,x)=\sum_{h\in \K^*} f_1(k h^{-1},h x) f_2(h,x).
$$

We see from Lemma \ref{neighb0ideles} above that the only effect of
the reduction to $\cG_\K^{(1)}$ is to remove from the noncommutative
space $\A_\K/\K^*$ all the elements of $C_\K$ whose class modulo
$\K^*$ does not intersect $\cG_\K^{(1)}$ (\ie in particular those
whose norm is greater than or equal to one). We then have the
following symmetries for the algebra $\cS(\cG_\K^{(1)})$.

\begin{prop}\label{globMorita}
Let $\cJ^+$ denote the semi-group of ideles $j\in \A_\K^*$ such that
$j \A^{(1)}_\K \subset \A^{(1)}_\K$.
 The semigroup $\cJ^+$ acts on the algebra $\cS(\cG_\K^{(1)})$ by
endomorphisms obtained as restrictions of the automorphisms of
$\cS(\cG_\K)$ of the form
\begin{equation}\label{thetagauto}
\urep(j) (f) (k,x) = f(k,j^{-1}x), \ \ \  \forall (k,x)\in
\cG_\K\,,\ \ j\in \cJ^+.
\end{equation}
Let $\K=\Q$ and $C^+_\Q\subset C_\Q$ be the semigroup $C^+_\Q=\{
g\in C_\Q |\, |g|<1 \}$. The semi-group $C_\Q^+$ acts on
$\cS(\cG_\Q^{(1)})$ by the endomorphisms
\begin{equation}\label{Fthetagauto}
F(g)=\urep(\bar g)
\end{equation}
 with
$\bar g$ the natural lift of $g\in C_\Q^+$ to $\hat\Z^*\times
\R^*_+$.
\end{prop}

\proof  By construction $\urep(j)$ is an automorphism of
$\cS(\cG_\K)$. For a function $f$ with support $B$ in the open set
$\cG_\K^{(1)}$ the support of the function $\urep(j)(f)$ is
$jB=\{(k,jx)|(k,x)\in B\}\subset \cG_\K^{(1)}$ so that $\urep(j)(f)$
still has support in $\cG_\K^{(1)}$.

For $\K=\Q$ let $\bar g\in \hat \Z^*\times \R^*_+$ be the natural
lift of an element $g\in C_\Q^+$. Then the archimedean component
$\bar g_\infty$ is of absolute value less than $1$ so that $ \bar
g\in \cJ^+$. The action of $\urep(\bar g)$ by endomorphisms of
$\cS(\cG_\Q^{(1)})$ induces a corresponding action of $C_\Q^+$.
\endproof

\begin{rem}\label{globMoritabis}{\rm
For $m$ a positive integer, consider the element $g=(1,m^{-1})\in
C_\Q^+$. Both $g=(1,m^{-1})$ and $\tilde m=(m,1)$ are in $\cJ^+$ and
have the same class in the idele class group $C_\Q$, since $m\,
g=\tilde m$. Thus the automorphisms $\urep(g)$ and $\urep(\tilde m)$
of $\cS(\cG_\K)$ are inner conjugate. Since the  open set
$\A_\K^{(1)}\subset \A_\K$ is not closed its characteristic function
is not continuous and does not define a multiplier of $\cS(\cG_\K)$.
It follows  that the endomorphism $F(g)$ is inner conjugate to the
endomorphism $\urep(\tilde m)$ only in the following weaker sense.
There exists a sequence  of elements $u_n$ of $\cS(\cG_\K^{(1)})$
such that for any $f\in \cS(\cG_\K^{(1)})$ with compact support
$$
F(g)(f)= u_n \, \urep(\tilde m)(f) \, u_n^* ,
$$
 holds for all $n$ large enough.}
\end{rem}

\subsection{The valuation systems}\label{valsystSect}

We now explain why the orbits $C_\K a^{(v)}$ appear indeed as the
set of classical points, in the sense of the low temperature KMS
states, of the noncommutative spaces $X_{\K,v}$. The notion of
classical points obtained from low temperature KMS states is
discussed at length in \cite{CMR2} (\cf also \cite{CM},
\cite{CMbook}, \cite{CMR}).

\smallskip
The noncommutative space $X_{\K,v}$ is described by the the
restricted groupoid
\begin{equation}\label{resgroupoid}
\cG(v)=\K^*\ltimes \A^{(1)}_{\K,v}=\{(g,a)\in \K^*\ltimes
\A_{\K,v}\,|\, a \  {\rm and} \  ga \in \A^{(1)}_{\K,v}\}\,.
\end{equation}
We  denote by $\varphi$ the positive functional on $C^*(\K^*\ltimes
\A^{(1)}_{\K,v})$  given by
\begin{equation}\label{resgroupoidstate}
\varphi(f)=\int_{\A^{(1)}_{\K,v}}\,f(1,a)\,da
\end{equation}

\begin{prop}\label{sectperorbmodular}
The modular automorphism group of the functional $\varphi$ on the
crossed product $C^*(\cG(v))$ is given by  the time evolution
\begin{equation}\label{sigmavevolmod}
\sigma_t^v(f)(k,x)=|k|_v^{it} \,f(k,x), \ \ \   \forall t \in \R,\ \
\forall f\in C^*(\K^*\ltimes \A^{(1)}_{\K,v})\,.
\end{equation}
\end{prop}

\proof  We identify elements of $C_c(\K^*\ltimes \A^{(1)}_{\K,v})$
with functions $f(g,a)$ of elements  $g\in \K^*$ and $a\in
\A^{(1)}_{\K,v}$. The product is simply of the form
$$ f_1 * f_2 (g,a)=\sum_r f_1(g\,r^{-1},ga)f_2(r,a). $$
The additive Haar measure $da$  on $\A_{\K,v}$ satisfies the scaling
property
\begin{equation}\label{muscalehaar}
d(ka) = |k|_v^{-1} da\,, \ \ \forall k\in \K^*\,,
\end{equation}
since the product measure $da\times da_v$ on $\A_\K=\A_{\K,v}\times
\K_v$ is invariant under the scaling by $k\in \K^*$, while the
additive Haar measure $da_v$ on $\K_v$ gets multiplied by $|k|_v$, namely
$d(ka_v)=|k|_vda_v$. We then check the KMS$_1$ condition, for
$\varphi$ associated to the additive Haar measure, as follows,
$$ \varphi(f_1*f_2) =\sum_r \int_{\A^{(1)}_{\K,v}}
f_1(r^{-1},r\,a)f_2(r,\,a)\, da$$ $$= \sum_r \int_{\A^{(1)}_{\K,v}}
f_2(k^{-1},k\,b) f_1(k,\,b)\, |k|_v^{-1}\,
db=\varphi(f_2*\sigma_i(f_1))\,,
$$
using the change of variables $k=r^{-1}$, $a=kb$ and
$da=|k|_v^{-1}\, db$.
\endproof

\smallskip
It is worthwhile to observe that these automorphisms extend to the
global algebra. Let $\cG^{(1)}_\K$ be the groupoid
$\K^*\ltimes \A^{(1)}_\K$ of Definition \ref{glMorgrpd}.

\begin{lem}\label{sigmavlem}
Let $\K$ be a global field and $v\in\Sigma_\K$ a place. The map
\begin{equation}\label{dvkxmap}
d_v(k,x)=\log |k|_v \in \R
\end{equation}
defines a homomorphism of the groupoid $\cG^{(1)}_\K$ to the
additive group $\R$ and the time evolution
\begin{equation}\label{sigmavevol}
\sigma_t^v(f)(k,x)=|k|_v^{it} \,f(k,x), \ \ \   \forall t \in \R,\ \
\forall f\in \cS(\cG^{(1)}_\K)
\end{equation}
generates a 1-parameter group of automorphisms of the algebra
$\cS(\cG^{(1)}_\K)$.
\end{lem}

The following result shows that the nontrivial part of the dynamics
$\sigma^v_t$ concentrates on the  algebra $\cS(\cG(v))$ with
$\cG(v)$ as in \eqref{resgroupoid}.

\begin{prop}\label{dvcocycle}
The morphism $\rho_v$ of \eqref{rhovrestr} restricts to a
$\sigma_t^v$-equivariant morphism $\cS(\cG^{(1)}_\K)\to
\cS(\cG(v))$. Moreover, the restriction of the one
parameter group $\sigma_t^v$ to the kernel of $\rho_v$ is inner.
\end{prop}

\proof For the first statement note that the proper inclusion
$\A_{\K,v} \subset \A_\K$ restricts to a proper inclusion
$\A^{(1)}_{\K,v} \subset \A^{(1)}_\K$. For the second statement,
notice that the formula
\begin{equation}\label{multiplhv}
h_v(x)=\log|x|_v, \ \ \  \forall x\in \A^{(1)}_{\K},
\end{equation}
defines the multipliers $e^{it h_v}$ of the kernel of $\rho_v$.
Indeed $e^{it h_v}$ is a bounded continuous function on $\A^{(1)}_\K
\smallsetminus \A^{(1)}_{\K,v}$.

We can then check that the $1$-cocycle $d_v$ is the coboundary of
$h_v$. In fact, we have
\begin{equation}\label{1cocyclecoboundhv}
h_v(k\,x) - h_v(x)=d_v(k,x), \ \ \  \forall (k,x)\in \cG^{(1)}_\K
\smallsetminus \cG(v).
\end{equation}
\endproof

\medskip

We now recall that, for an \'etale groupoid like $\cG(v)$, every
unit $y\in \cG(v)^{(0)}$ defines, by
\begin{equation}\label{cGpiybis}
(\pi_y(f)\xi)(\gamma) = \sum_{\gamma_1\gamma_2=\gamma}
f(\gamma_1)\xi(\gamma_2),
\end{equation}
a representation $\pi_y$ by left convolution of the algebra of
$\cG(v)$ in the Hilbert space $\cH_y=\ell^2(\cG(v)_y)$, where
$\cG(v)_y$ denotes the set of elements of the groupoid $\cG(v)$ with
source $y$. By construction the unitary equivalence class of the
representation $\pi_y$ is unaffected when one replaces $y$ by an
equivalent $z\in \cG(v)^{(0)}$ \ie one assumes that there exists
$\gamma\in \cG(v)$ with range and source $y$ and $z$. Thus we can
think of the label $y$ of $\pi_y$ as living in the quotient space
$X_{\K,\,v}$ of equivalence classes of elements of $\cG(v)^{(0)}$.

\smallskip
The relation between $\Xi_{\K,\,v}$ and $X_{\K,\,v}$ is then the
following.

\begin{thm}\label{classpointsperorb}  For $y\in X_{\K,\,v}$, the representation $\pi_y$ is a
positive energy representation if and only if $y\in \Xi_{\K,\,v}$.
\end{thm}

\proof Let first $y \in \cG(v)^{(0)}\cap \Xi_{\K,\,v}$. Thus one has
$y\in \A^{(1)}_{\K,v}$, $y_w\neq 0$ for all $w$ and $|y_w|=1$ for
all $w\notin S$ where $S$ is a finite set of places. We can identify
$\cG(v)_y$ with the set of $k\in \K^*$  such that $k\,y \in
\A^{(1)}_{\K,v}$. We extend $y$ to the adele $\tilde y=y\times 1$
whose component at the place $v$ is equal to $1\in \K_v$. Then
$\tilde y$ is an idele. Thus by Lemma \ref{neighb0ideles} the number
of elements of the orbit $\K^*\tilde y$ in a given compact subset of
$\A_\K$ is finite. It follows that $\log |k|_v$ is lower bounded on
$\cG(v)_y$. Indeed otherwise there would exist a sequence $k_n\in
\K^*\cap \cG(v)_y$ such that $|k_n|_v\to 0$. Then $k_n \,\tilde y\in
\A^{(1)}_\K$ for all $n$ large enough and this contradicts the
discreteness of $\K^*\tilde y$. In the representation $\pi_y$ the
time evolution $\sigma_t$ is implemented by the Hamiltonian $H_y$ given by
\begin{equation}\label{hamilthy}
(H_y\, \xi)(k,y) =\log |k|_v \,\, \xi(k,y).
\end{equation}
Namely, we have
\begin{equation}\label{Hysigmatbis}
 \pi_y(\sigma_t(f))=e^{itH_y} \pi_y(f) e^{-itH_y}, \ \ \ \ \forall
f\in C_c(\cG(v))\,.
\end{equation}
Thus since $\log |k|_v$ is lower bounded on $\cG(v)_y$ we get that
the representation $\pi_y$ is a positive energy representation.

\smallskip Let then $y \in \cG(v)^{(0)}\smallsetminus \Xi_{\K,\,v}$.
We shall show that $\log |k|_v$ is not lower bounded on $\cG(v)_y$,
and thus that $\pi_y$ is not a positive energy representation. We
consider as above the adele $\tilde y=y\times 1$ whose component at
the place $v$ is equal to $1\in \K_v$. Assume that $\log |k|_v$ is
 lower bounded on $\cG(v)_y$. Then there exists $\epsilon>0$ such
 that, for $k\in \K^*$,
 $$
 k\,y\in \A^{(1)}_{\K,v}\Rightarrow |k|_v\geq \epsilon \,.
 $$
 This shows that the neighborhood of $0\in \A_\K$ defined as
 $$
W=\{ a \in \A_\K\,;\, |a_v|<\epsilon \,,\ a_w\in \K_w^{(1)}\,, \
\forall w\neq v\}
 $$
 does not intersect $\K^*\tilde  y$. Thus by Lemma
 \ref{neighb0ideles} we get that $\tilde  y$ is an idele
 and $y\in \Xi_{\K,\,v}$.
\endproof

\smallskip
The specific example of the Bost-Connes system combined with Theorem
\ref{classpointsperorb} shows that one can refine the recipe of
\cite{CMR2} (\cf also \cite{CM}, \cite{CMbook}, \cite{CMR}) for
taking ``classical points" of a noncommutative space. The latter
recipe only provides a notion of classical points that can be
thought of, by analogy with the positive characteristic case, as
points defined over the mysterious ``field with one element" $\F_1$
(see \eg \cite{Man-zetas}). To obtain instead a viable notion of the
 points defined over the maximal unramified extension $\bar
\F_1$, one performs the following sequence of operations.

\begin{equation}\label{takingpoints}
X \stackrel{Dual\, System} \longrightarrow \hat X
\stackrel{Periodic\, Orbits}\longrightarrow \cup\, {\hat X}_v
\stackrel{Classical\, Points}\longrightarrow \cup\,\Xi_v
\end{equation}

\smallskip

which make sense in the framework of endomotives of \cite{CCM}. Note
in particular that the dual system $\hat X$ is of type II and as
such does not have a non-trivial time evolution. Thus it is only by
restricting to the periodic orbits that one passes to noncommutative
spaces of type III for which the cooling operation is non-trivial.
In the analogy with geometry in non-zero characteristic, the set of
points $X(\bar\F_q)$ over $\bar\F_q$ of a variety $X$ is indeed
obtained as the union of the periodic orbits of the Frobenius.

\medskip
\begin{rem} \label{kmsbeware} {\rm Theorem \ref{classpointsperorb}
does not give the classification of KMS$_\beta$ states for the
quantum statistical system $(C^*(\K^*\ltimes
\A^{(1)}_{\K,v}),\sigma_t)$. It just exhibits extremal KMS$_\beta$
states but does not show that all of them are of this form. }
\end{rem}

\medskip

\subsection{The curve inside the adeles class
space}\label{curveinSect}

In the case of a function field $\K=\F_q(C)$, the set of {\em
periodic classical} points of the adeles class space $\A_\K/\K^*$ is
(non-canonically) isomorphic to the algebraic points $C(\bar\F_q)$.
In fact, more precisely the set of algebraic points $C(\bar\F_q)$ is
equivariantly isomorphic to the quotient $\Xi_\K/C_{\K,1}$ where
$C_{\K,1}\subset C_\K$ is the kernel of the norm $|\cdot|: C_\K \to
\R^*_+$, and $\Xi_\K$ is as in \eqref{Xiorbits}.

\begin{prop}\label{XiKalgpts}
For $\K=\F_q(C)$ a function field, the orbits of Frobenius on
$C(\bar\F_q)$ give an equivariant identification
\begin{equation}\label{XiKCbarFq}
\Xi_\K/C_{\K,1} \simeq C(\bar\F_q),
\end{equation}
between $\Xi_\K/C_{\K,1}$ with the action of $q^\Z$ and
$C(\bar\F_q)$ with the action of the group of integer powers of the
Frobenius.
\end{prop}

\proof At each place $v\in \Sigma_\K$ the quotient group of the
range $N$ of the norm $|\cdot|: C_\K \to \R^*_+$ by the range $N_v$
of $|\cdot|: \K_v \to \R^*_+$ is the finite cyclic group
\begin{equation}\label{NNvqZ}
N/N_v = q^\Z/ q^{n_v\Z} \simeq \Z/n_v\Z,
\end{equation}
where $n_v$ is the degree of the place $v\in \Sigma_\K$. The degree
$n_v$ is the same as the cardinality of the orbit of the Frobenius
acting on the fiber of the map \eqref{ptstoplaces} from algebraic
points in $C(\bar \F_q)$ to places in $\Sigma_\K$. Thus, one can
construct in this way an equivariant embedding
\begin{equation}\label{equivembedCXi}
C(\bar \F_q) \hookrightarrow (\A_\K/\K^*)/C_{\K,1}
\end{equation}
obtained, after choosing a point in each orbit, by mapping the orbit
of the integer powers of the Frobenius in $C(\bar\F_q)$ over a place
$v$ to the orbit of $C_\K/C_{\K,1}\sim q^\Z$ on the adele $a^{(v)}$.
\endproof

Modulo the problem created by the fact that the identification above
is non-canonical and relies upon the choice of a point in each
orbit, it is then possible to think of the locus $\Xi_\K$, in the
number field case, as a replacement for $C(\bar\F_q)$ inside the
adeles class space $\A_\K/\K^*$.

\begin{center}
\begin{figure}
\includegraphics[scale=0.9]{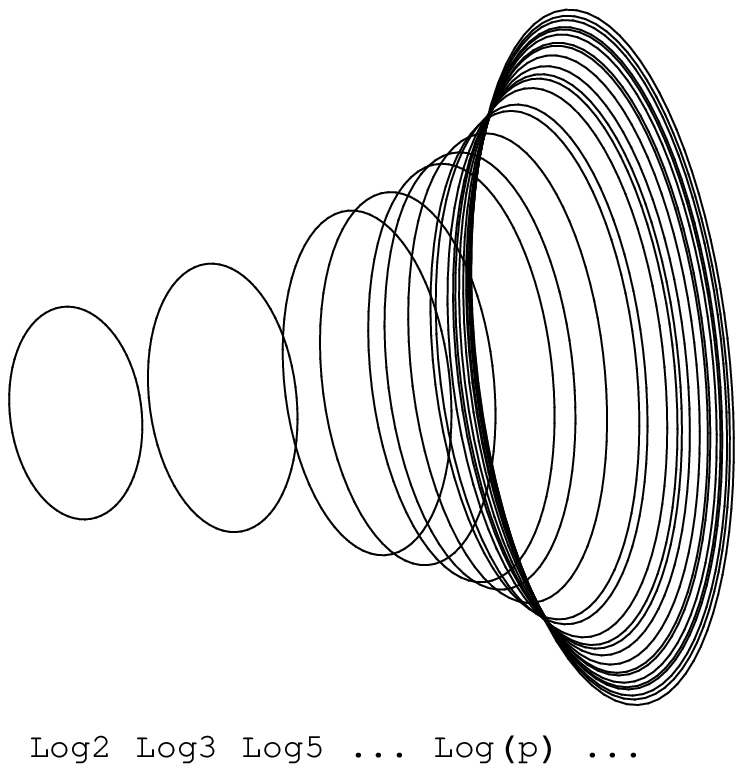}
\caption{The classical points $\Xi_\Q/C_{\Q,1}$ of the adeles class
space $\A_\Q/\Q^*$. \label{Figpoints}}
\end{figure}
\end{center}

In the case of $\K=\Q$, the quotient $\Xi_\Q/C_{\Q,1}$ appears as a
union of periodic orbits of period $\log p$ under the action of
$C_\Q/C_{\Q,1}\sim \R$, as in Figure \ref{Figpoints}. What matters,
however, is not the space $\Xi_\Q/C_{\Q,1}$ in itself but the way it
sits inside $\A_\Q/\Q^*$. Without taking into account the topology
induced by $\A_\K$ the space $\Xi_\K$ would just be a disjoint union
of orbits without any interesting global structure, while it is the
embedding in the adeles class space that provides the geometric
setting underlying the Lefschetz trace formula of \cite{Co-zeta} and
its cohomological formulation of \cite{CCM}.

\subsection{The valuation systems for $\K=\Q$}\label{kqvalSect}

 We concentrate again on the specific case of $\K=\Q$
to understand better the properties of the dynamical systems
$\sigma^p_t$ associated to the finite primes $p\in \Sigma_\Q$.

We know that, in the case of the BC system, the KMS state
at critical temperature $\beta=1$ is given by the
additive Haar measure on finite adeles \cite{BC}. Thus, one expects
that, for the systems associated to the finite primes, the additive
Haar measure of $\A_{\Q,p}$ should play an analogous role.

\begin{defn}\label{Dpdomaindef}
Let  $\A_{\Q,p}^*\subset \A^{(1)}_{\Q,p}$ be the subspace
\begin{equation}\label{Dpstardomain}
\A_{\Q,p}^*=\{ x \in \A_{\Q,p} |\  |x_w|= 1 \ \forall w\neq p,\infty \ \text{
and } p^{-1} \leq |x_\infty|<1 \}.
\end{equation}
\end{defn}
As above $\cG(p)$ denotes the reduction of the groupoid $\cG_{\Q,p}$
by the open subset $\A^{(1)}_{\Q,p}\subset \A_{\Q,p}$, namely
\begin{equation}\label{GQpreduction}
\cG(p)=\{(k,x)\in \cG_{\Q,p}\, | \, x\in \A^{(1)}_{\Q,p}, \, kx\in
\A^{(1)}_{\Q,p}\}.
\end{equation}

Notice that the set $\A^{(1)}_{\Q,p}$ meets all the equivalence
classes in $\A_{\Q,p}$ by the action of $\Q^*$. In fact, given $x\in
\A_{\Q,p}\,$, one can find a representative $y$ with $y\sim x$ in
$\A_{\Q,p}/\Q^*$, such that $y \in\hat \Z\times \R$. Upon
multiplying $y$ by a suitable power of $p$, one can make $y_\infty$
as small as required, and in particular one can obtain in this way a
representative in $\A^{(1)}_{\Q,p}$. Let us assume that $|y_w|=1$
for all finite places $w\neq p$ and that $y_\infty > 0$. Then there
exists a unique $n\in \N\cup\{0\}$ such that $p^n\,y\in
\A_{\Q,p}^*$.

\smallskip
Given a prime $p$ we define the function $f_p(\lambda,\beta)$ for
$\lambda\in (1,p\,]$ and $\beta>1$ by
\begin{equation}\label{ZpbetaDef}
f_p(\lambda,\beta)=\sum c_k\,p^{-k\beta}
\end{equation}
where the $c_k\in\{0,\ldots p-1\}$ are the digits of the expansion
of $\lambda$ in base $p$. There is an ambiguous case where all
digits $c_k$ are equal to $0$ for $k> m$ while $c_m>0$, since the
same number
$$
\lambda =\sum c_k\,p^{-k}
$$
is obtained using the same $c_j$ for $j<m$, $c_m-1$ instead of $c_m$
and $c_j=p-1$ for $j>m$. In that case, for $\beta>1$,
\eqref{ZpbetaDef} gives two different values and we choose the value
coming from the second representation of $\lambda$, \ie the lower of
the two. These coefficients $c_k$ of the expansion of $\lambda$ in
base $p$ are then given by
\begin{equation}\label{coeffceiling}
c_k=\lceil \lambda p^k -1 \rceil - p\,\lceil \lambda p^{k-1} -1
\rceil\,,
\end{equation}
where $\lceil x \rceil=\inf_{n\in\Z} \{ n\geq x \}$ denotes the
ceiling function.

\begin{center}
\begin{figure}
\includegraphics[scale=0.9]{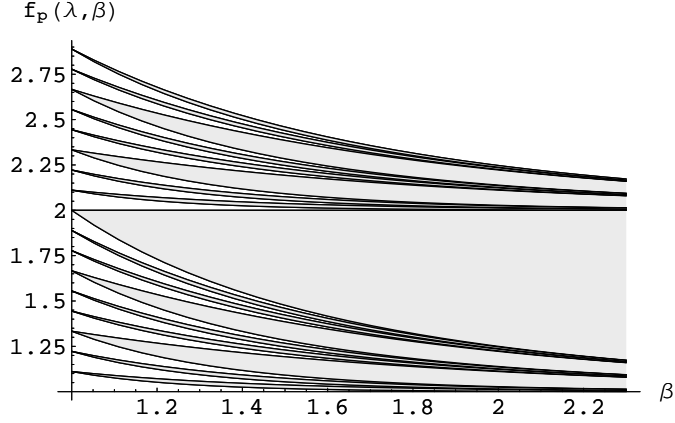}
\caption{Graphs of the functions $f_p(\lambda,\beta)$ as   functions
of $\beta$ for $p=3$, $\lambda=n/27$. The gray regions are the gaps
in the range of $f_p$. \label{Zpfunct}}
\end{figure}
\end{center}

\begin{center}
\begin{figure}
\includegraphics[scale=0.6]{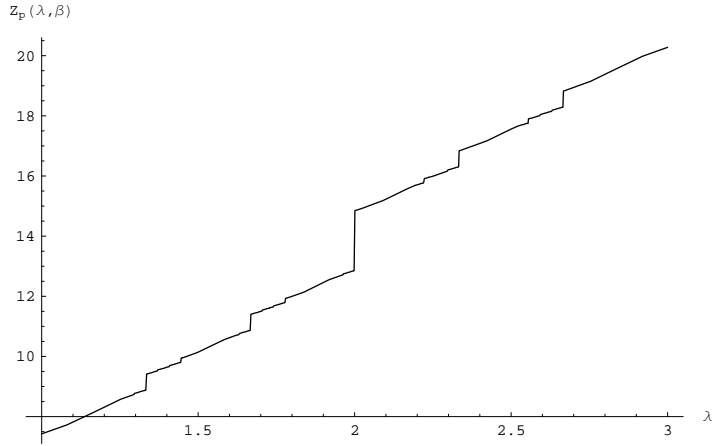}
\caption{Graph of the function $Z_p(\lambda,\beta)$ as a function
of $\lambda$ for $p=3$, $\beta=1.2$. \label{Zpfunct2}}
\end{figure}
\end{center}

\smallskip Note that, for $\beta>1$, the function $f_p(\lambda,\beta)$ is
discontinuous (\cf Figures  \ref{Zpfunct} and \ref{Zpfunct2}) at any
point $(\lambda,\beta)$ where the expansion of $\lambda$ in base $p$
is ambiguous, \ie $\lambda\in \N\,p^{-k}$. Moreover for $\beta=1$
one gets
\begin{equation}\label{ZpbetaDefbone}
f_p(\lambda,1)=\lambda \,,\ \ \forall \lambda\in (1,p\,]\,.
\end{equation}

We then obtain the following result.

\begin{thm}\label{DQprepsHam}
Let $(C^*(\cG(p)),\sigma_t^p)$ be the $C^*$-dynamical system
associated to the groupoid \eqref{GQpreduction} with the time
evolution \eqref{sigmavevol}. Then the following properties hold:
\begin{enumerate}
\item For any $y\in \A_{\Q,p}^*$ the corresponding
representation $\pi_y$ has positive energy.
\item Let $H_y$ denote the Hamiltonian implementing the time evolution
in the representation $\pi_y$, for $y\in \A_{\Q,p}^*$ with $y_\infty
=\lambda^{-1}$ and $\lambda \in (1,p\,]$. Then the partition function
is given by
\begin{equation}\label{Trzp}
Z_p(\lambda,\beta)=\Tr (e^{-\beta H_y}) = 2\,
\frac{1-p^{-\beta}}{1-p^{1-\beta}}\,f_p(\lambda,\beta)\,.
\end{equation}
\item The functionals
\begin{equation}\label{kmszp}
\psi_{\beta,\,y}(a)=\Tr (e^{-\beta H_y}\,\pi_y(a))\,,\ \  \forall a\in
C^*(\cG(p))
\end{equation}
satisfy the KMS$_\beta$ condition for $\sigma_t^p$ and depend weakly
continuously on the parameter $y\in \A_{\Q,p}^*$.
\end{enumerate}
\end{thm}

\proof (1) This follows from Theorem \ref{classpointsperorb}. For
$y\in \A_{\Q,p}^*$ one has
\begin{equation}\label{posener0Qp}
r\in \Q^*,\  ry\in \A^{(1)}_{\Q,p} \Longrightarrow r=p^{-k}m,
\end{equation}
for some $k\geq 0$ and some integer $m$
prime to $p$ and such that $|r\,y_\infty| <1$.
This implies
\begin{equation}\label{posener1Qp}
|m|< p^{k+1}
\end{equation}
and one finds
\begin{equation}\label{posener2Qp}
|r|_p=p^k\geq 1 \ \ \text{ and } \ \ \log |r|_p \geq 0,
\end{equation}
In fact, the argument above shows that the spectrum of the
Hamiltonian $H_y$ implementing the time evolution $\sigma_t^p$ in
the representation $\pi_y$ is given by
\begin{equation}\label{SpHyQp}
\Sp(H_y)=\{ k\,\log p\}_{k\in \N\cup\{ 0 \}},
\end{equation}
hence $\pi_y$ is a positive energy representation.

(2) We begin by the special case with $y_\infty=p^{-1}$. Then
$\lambda=p$ and $f_p(\lambda,\beta)=\frac{p-1}{1-p^{-\beta}}$ since
all digits of $\lambda=p$ are equal to $p-1$. We want to show that
the partition function is given by
\begin{equation}\label{ZbetaDQp}
\Tr (e^{-\beta H_y}) = 2\,\frac{p-1}{1-p^{1-\beta}}.
\end{equation}
The multiplicity of an eigenvalue $k\log p$ of $H_y$ is the number
of integers $m\neq 0 \in \Z$ that are prime to $p$ and such that
$p^{-k}\, |m| \,y_\infty <1$. Since we are assuming that $y_\infty
=\,p^{-1}$, this   gives $|m|< p^{k+1}$. Thus, the multiplicity is
just $2\,(p^{k+1}-p^{k})$. The factor $2$ comes from the sign of the
integer $m$. The factor $(p^{k+1}-p^{k})$ corresponds to subtracting
from the number $p^{k+1}$ of positive integers $m\leq p^{k+1}$ the
number $p^{k}$ of those that are multiples of $p$.

\smallskip
We now pass to the general case. For $x>0$, $\lceil x-1\rceil $  is
the cardinality of $(0,x)\cap \N$. The same argument used above
shows that the multiplicity of the eigenvalue $k\log p$ is given by
the counting
$$2 \left( \lceil \lambda p^k -1 \rceil - \lceil \lambda p^{k-1} -1
\rceil \right).$$ Thus
\begin{equation}\label{Trceiling}
\Tr (e^{-\beta H_y}) = 2\, \sum_{k=0}^\infty \left( \lceil \lambda
p^k -1 \rceil - \lceil \lambda p^{k-1} -1 \rceil \right)\,
p^{-k\beta}\,.
\end{equation}
One has the following equalities of convergent series,
$$
\sum_{k=0}^\infty \left( \lceil \lambda p^k -1 \rceil - \lceil
\lambda p^{k-1} -1 \rceil \right)\,p^{-k\beta}=\sum_{k=0}^\infty
\lceil \lambda p^k -1 \rceil\,(p^{-k\beta}-p^{-(k+1)\beta})
$$
so that,
\begin{equation}\label{Trceiling1}
\Tr (e^{-\beta H_y}) = 2\,(1-p^{-\beta})\sum_{k=0}^\infty \lceil
\lambda p^k -1 \rceil\,p^{-k\beta}\,.
\end{equation}
Similarly
$$
\sum_{k=0}^\infty \left( \lceil \lambda p^k -1 \rceil - p\,\lceil
\lambda p^{k-1} -1 \rceil \right)\,p^{-k\beta}=\sum_{k=0}^\infty
\lceil \lambda p^k -1 \rceil\,(p^{-k\beta}-p\,p^{-(k+1)\beta})
$$
which gives
\begin{equation}\label{Trceiling2}
f_p(\lambda,\beta)=(1-p^{1-\beta})\sum_{k=0}^\infty \lceil \lambda
p^k -1 \rceil\,p^{-k\beta}\,,
\end{equation}
since the coefficients $c_k$ of the expansion of $\lambda$ in base
$p$ are given by \eqref{coeffceiling}. Combining \eqref{Trceiling1}
with \eqref{Trceiling2} gives \eqref{Trzp}.

(3) It follows from \eqref{Hysigmatbis} and the finiteness of the
partition function \eqref{Trzp} that the functionals \eqref{kmszp}
fulfill the KMS$_\beta$ condition. In terms of functions on the
groupoid $\cG(p)$ one has
\begin{equation}\label{kmszpfun}
\psi_{\beta,\,y}(f)=\sum\, f(1,n\,p^{-k}\,y)\,p^{-k\beta}\,,\ \  \forall f\in
C_c(\cG(p))
\end{equation}
where the sum is absolutely convergent. Each of the terms in the sum
gives a weakly continuous linear form thus one obtains the required
continuity.
\endproof

\begin{rem}{\rm The partition function $Z_p(\lambda,\beta)$ is a discontinuous function
of the parameter $\lambda$ and this might seem to contradict the
third statement of Theorem \ref{DQprepsHam}. It would if the algebra
$C^*(\cG(p))$ were unital since, in that case, the partition function
is given by evaluation on the unit and weak continuity implies that
it is continuous. In our case $C^*(\cG(p))$ is {\em not} unital, and
the partition function is expressed as a supremum of the form
$$
Z_p(\lambda,\beta)=\sup\{\psi_{\beta,\,y}(a^*a) |a\in
C^*(\cG(p))\,,\ ||a||\leq 1\}\,.
$$
In particular it shows that $Z_p(\lambda,\beta)$ is lower
semi-continuous as a function of $\lambda$. }
\end{rem}

\smallskip
The precise qualitative properties of the partition functions
$Z_p(\lambda,\beta)$ are described by the following result

\begin{prop} As a function of $\lambda \in (1,\lambda]$ the
partition function $Z_p(\lambda,\beta)$ satisfies for $\beta>1$:
\begin{enumerate}
  \item $Z_p$ is strictly increasing.
    \item $Z_p$ is continuous on the left, and lower semi-continuous.
  \item $Z_p$ is discontinuous at any point of the form $\lambda=m\,p^{-k}$
  with a jump of $2\,p^{-k\beta}$ (for $m$ prime to $p$).
  \item The measure $\frac{\partial Z_p}{\partial \lambda}$ is the
  sum of the Dirac masses at the points $\lambda=m\,p^{-k}$, $m$ prime to
  $p$, with coefficients $2\,p^{-k\beta}$.
  \item The closure of the range of $Z_p$ is a Cantor set.
\end{enumerate}
\end{prop}

\proof (1) This follows from \eqref{Trceiling1} which expresses
$Z_p$ as an absolutely convergent sum of multiples of the functions
$\lceil \lambda p^k -1 \rceil$. The latter
are  non-decreasing and jump by
$1$ at $\lambda \in \N \,p^{-k}\cap (1,p\,]$. The density of the
union of these finite sets for $k \geq 0$ shows that $Z_p$ is
strictly increasing.

\smallskip
(2) This follows as above from \eqref{Trceiling1} and the
semi-continuity properties of the ceiling function.

\smallskip
(3) Let $\lambda=m\,p^{-k}$ with $m$ prime to $p$. Then for any
$j\geq k$ one gets a jump of $2\,(1-p^{-\beta})\,p^{-j\beta}$ coming
from \eqref{Trceiling1} so that their sum gives
$$
2\,(1-p^{-\beta})\,\sum_{j=k}^\infty p^{-j\beta}=\,2\,p^{-k\beta}
$$

\smallskip
(4) This follows as above from \eqref{Trceiling1} and from (3) which
computes the discontinuity at the jumps.

\smallskip
(5) Recall that when writing elements of an interval in base $p$ one
gets a map from the cantor set to the interval. This map is
surjective but fails to be injective due to the identifications
coming from $\sum_0^\infty (p-1)\,p^{-m}=p$. The connectedness of
the interval is recovered from these identifications. In our case
the coefficients $c_k$ of the expansion in base $p$ of elements of
$(1,p\,]$ are such that $c_0\in \{1,\ldots,p-1\}$ while $c_k\in
\{0,\ldots,p-1\}$ for $k>0$. This is a Cantor set $K$ in the product
topology of $K=\{1,\ldots,p-1\}\times \prod_\N\,\{0,\ldots,p-1\}$.
As shown in Figure \ref{Zpfunct2}, the discontinuities of the
function $Z_p(\lambda,\beta)$ as a function of $\lambda$ replace the
connected topology of $(1,p\,]$ by the totally disconnected topology
of $K$.
\endproof

\begin{center}
\begin{figure}
\includegraphics[scale=0.7]{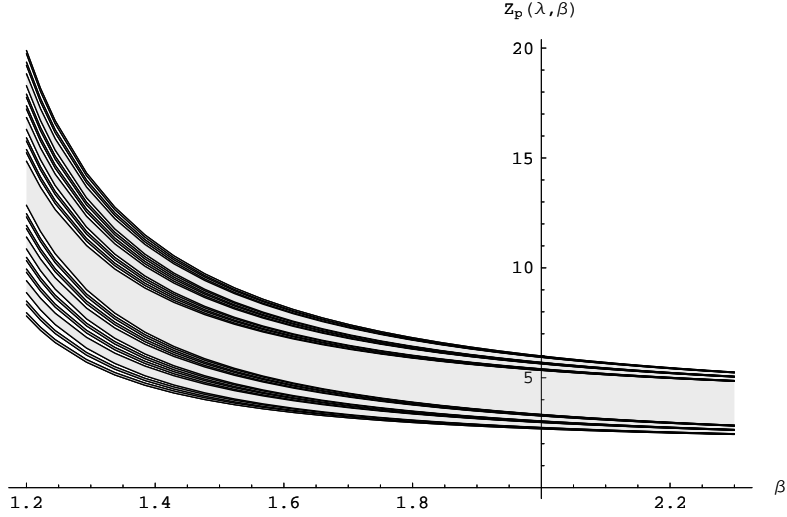}
\caption{Graphs of the functions $Z_p(\lambda,\beta)$ as   functions
of $\beta$ for $p=3$, $\lambda=n/27$. The gray regions are the gaps
in the range. All these functions have a pole at
$\beta=1$.\label{Zpfunct1}}
\end{figure}
\end{center}

\begin{rem}{\rm
One can use  \eqref{Trceiling} to define $Z_p(\lambda,\beta)$  for
any $\lambda>0$, as
\begin{equation}\label{fplambdabeta}
Z_p(\lambda,\beta)=2 \sum_{-\infty}^\infty \left(\lceil \lambda
p^k-1 \rceil - \lceil \lambda p^{k-1} -1 \rceil \right) p^{-k\beta}.
\end{equation}
This makes sense for $\Re(\beta) >1$ since $\lceil \lambda p^k -1
\rceil =0$ for $k\leq -\frac{\log \lambda}{\log p}$. The extended
function \eqref{fplambdabeta} satisfies
$$ Z_p(p\lambda,\beta) = p^\beta Z_p(\lambda,\beta), $$
which suggests replacing $Z_p(\lambda,\beta)$ with
\begin{equation}\label{zetap}
\zeta_p(\lambda,\beta)=\lambda^{-\beta} Z_p(\lambda,\beta)
\end{equation}
so that
\begin{equation}\label{zetap1}
\zeta_p(p\lambda,\beta)=\zeta_p(\lambda,\beta).
\end{equation}
This replacement $Z_p \mapsto \zeta_p$ corresponds to the shift in
the Hamiltonian $H_y$ by
$$
H_y\mapsto H_y - \log |y_\infty|.
$$}
\end{rem}

\smallskip

We can now refine Theorem \ref{classpointsperorb} and consider the
zero temperature KMS state of the system $(C^*(\cG(p)),\sigma_t^p)$
corresponding to the positive energy representation $\pi_y$ for
 $y\in \Xi_{\Q,\,p}\,$.

\begin{prop}\label{vacuumQp}
As $\beta\to \infty$ the vacuum states (zero temperature KMS states)
of the system $(C^*(\cG(p)),\sigma_t^p)$ with Hamiltonian $H_y$ have
a degeneracy of $2\lceil \lambda -1 \rceil$, where $y_\infty
=\lambda^{-1}$.  There is a preferred choice of a vacuum state given
by the evaluation at $y\in \A_{\Q,p}^*\,$.
\end{prop}

\proof When we look at the orbit of $y\in \A_{\Q,p}^*$, \ie at the
intersection $\Q^* y\cap \A^{(1)}_{\Q,p}$ and label its elements by
pairs $(k,m)$ as above, we find that all elements with $k=0$ give a
ground state. This degeneracy of the vacuum reflects the fact that
the limit of the partition function as the temperature goes to $0$
is not in general equal to $1$. For instance, for
$y_\infty =p^{-1}$, one finds
$$
\lim_{\beta\to \infty} \Tr (e^{-\beta H_y})= \lim_{\beta\to \infty}
2\frac{p-1}{1-p^{1-\beta}} =2(p-1).
$$
More generally, one finds similarly the limit
$$
\lim_{\beta\to \infty} \Tr (e^{-\beta H_y})= 2\lceil \lambda-1
\rceil.
$$
Among the $2\lceil \lambda-1 \rceil$ vacuum states, the state given by
evaluation at $y\in \A_{\Q,p}^*$ is singled out, since $m y\notin \A_{\Q,p}^*$
for $m\neq 1$. It is then natural to consider, for each finite place
$p\in \Sigma_\Q$, the section
\begin{equation}\label{sectionp}
s_p(x)=\Q^* x\cap \A_{\Q,p}^*,   \ \ \ \forall x\in C_\Q a^{(p)}\subset \A_\Q/\Q^*
\end{equation}
of the projection from $\A_\Q$ to the orbit $C_\Q a^{(p)}$.
\endproof

Notice that $s_p$ is discontinuous at the boundary of the domain
$\A_{\Q,p}^*$. Indeed when $y_\infty$ crosses the value $p^{-1}$ the
class in $C_\Q a^{(p)}$ varies continuously but the representative
in $\A_{\Q,p}^*$ jumps discontinuously so that its archimedian
component  remains in the interval $[p^{-1},1)$. This suggests to
consider a  cyclic covering of $\Xi_\Q$ which we now discuss in \S
\ref{curvecoverSect}.

\medskip

\subsection{The  cyclic covering $\tilde\Xi_{\Q }$ of $\Xi_{\Q }$}\label{curvecoverSect}

By construction $\Xi_\K$ is a subspace of the adeles class space
$X_\K$. We shall now show, in the case $\K=\Q$, that it admits a
natural lift $\tilde\Xi_{\Q }$ to a subspace of $\A_\Q$ which
reduces the ambiguity group $\Q^*$ to a cyclic group. One thus
obtains a natural cyclic covering $\tilde\Xi_{\Q }\subset \A_\Q$ of
$\Xi_\Q$. We already saw above, in Proposition \ref{vacuumQp}, that
it is natural to choose representatives for the elements of the
orbit $C_\Q a^{(p)}$, for a finite prime $p$, in the subset of
adeles given by
\begin{equation}\label{cqtilde}
\tilde\Xi_{\Q,p}:=\{y\in \A_\Q |\ y_p=0 \text{ and } |y_\ell|=
1 \text{ for } \ell\neq p,\infty \text{ and } y_\infty >0\}.
\end{equation}
We extend this definition at $\infty$ by
\begin{equation}\label{cqtildeinf}
\tilde\Xi_{\Q,\infty}:=\{y\in \A_\Q|\ |y_w|= 1 \ \forall
w\neq \infty \text{ and } y_\infty =0\}.
\end{equation}

\begin{defn}\label{Xitildecover}
The locus $\tilde \Xi_\Q\subset \A_\Q$ is defined as
\begin{equation}\label{XitildeQ}
\tilde \Xi_\Q=\bigcup_{v\in
\Sigma_\Q}\tilde\Xi_{\Q,v} \subset  \A_\Q
\end{equation}
\end{defn}

We then have the following simple fact.

\begin{prop}
Let $\pi$ be the projection from $\tilde \Xi_\Q$ to $\Xi_\Q$, with
$\pi(x)$ the class of $x$ modulo the action of $\Q^*$.
\begin{enumerate}
\item The map $\pi:\tilde \Xi_\Q \to \Xi_\Q$ is surjective.
\item Two elements in $\tilde\Xi_{\Q,v}$ have the same image in
$C_\Q a^{(v)}$ iff they are on the same orbit of the following
transformation $T$
\begin{equation}\label{XitildeQT}
Tx=p\,x \,,\ \  \forall x\in\tilde\Xi_{\Q,p}\,, \ \ Tx=-x \,,\ \
\forall x\in\tilde\Xi_{\Q,\infty}
\end{equation}
\end{enumerate}
\end{prop}

\proof The first statement follows by lifting $C_\Q$ inside
$\A_\Q^*$ as the subgroup $\hat\Z^*\times \R_+^*$. Then any element
of $C_\Q a^{(v)}$ has a representative in $(\hat\Z^*\times
\R_+^*)a^{(v)}$

The proof of the second statement is straightforward, since for a
finite prime $p$ the subgroup $p^\Z\subset \Q^*$ is the group of
elements of $\Q^*$ which leave $\tilde\Xi_{\Q,p}$ globally
invariant.
\endproof

\medskip

\subsection{Arithmetic subalgebra, Frobenius and
monodromy}\label{arithsubalg}

We now describe a natural algebra of coordinates $\cB$ on $\Xi_\Q$.

The BC system of \cite{BC}, as well as its arithmetic generalizations
of \cite{CM} and \cite{CMR}, have the important property that they
come endowed with an arithmetic structure given by an arithmetic
subalgebra. The general framework of endomotives developed in
\cite{CCM} shows a broad class of examples where a similar arithmetic
structure is naturally present. We consider here the issue of
extending the construction of the ``rational subalgebra" of the
BC-system to the algebra $\cS(\cG_\Q^{(1)})$ of \S \ref{globMor}.

In order  to get a good geometric picture it is convenient to think
in terms of $\Q$-lattices rather than of adeles, as in \cite{CM}.
Thus, we let $\cL$ denote the set of 1-dimensional $\Q$-lattices
(as defined in \cite{CM}). We consider the map
\begin{equation}\label{iota}
\iota : \hat \Z\times \R^*_+\to \cL \,, \ \ \
\iota(\rho,\lambda)= (\Lambda,\phi)=
(\lambda^{-1}\Z,\lambda^{-1}\rho)
\end{equation}
which associates to an adele $(\rho,\lambda)\in \hat \Z\times
\R^*_+\subset \A_\Q$ the $\Q$-lattice obtained using $\rho$ to
label the torsion points of $\R/\lambda^{-1}\Z$. Replacing
$(\rho,\lambda)$ by $(n\rho,n\lambda)$, for a positive integer
$n\in \N$, one obtains the pair $(\frac{1}{n}\Lambda,\phi)$, which
is commensurable to $(\Lambda,\phi)$. Thus, the action of $\Q_+^*$
corresponds to commensurability of $\Q$-lattices under the map
$\iota$. Multiplying $\lambda$ by a positive scalar corresponds to
the scaling action of $\R^*_+$ on $\Q$-lattices.

Let us recall the definition of the ``rational algebra" $\cA_\Q$ of
\cite{CM} for the BC system, given in terms of $\Q$-lattices. We let
\begin{equation}\label{epsilonk}
\epsilon_{a} (\Lambda , \phi)= \sum_{y\in  \Lambda +\phi(a)} y^{-1},
\end{equation}
for any $a\in \Q/\Z$. This is well defined, for $\phi(a)\neq 0$,
using the summation $\lim_{N \to \infty}\sum_{-N}^N$, and is zero by
definition for $\phi(a)=0$. The function
\begin{equation}\label{phia}
\varphi_a(\rho,\lambda)=\epsilon_{a}(\iota(\rho,\lambda)), \ \  \ \ \forall
(\rho,\lambda)\in \hat \Z\times \R^*_+ ,
\end{equation}
is well defined and homogeneous of degree $1$ in $\lambda$. Moreover,
for fixed $a\in \Q/\Z$ with denominator $m$, it only depends upon the
projection of $\rho$ on the finite group $\Z/m\Z$, hence it defines
a continuous function on $\hat \Z\times \R^*_+$. Using the degree $1$
homogeneity in $\lambda$, one gets that \eqref{phia} extends by
continuity to $0$ on $\hat \Z\times \{0\}$.

One gets functions that are homogeneous of weight zero by taking
the derivatives of the functions $\varphi_a$. The functions
\begin{equation}\label{psia}
\psi_a(\rho,\lambda)= \frac{1}{2\pi
i}\frac{d}{d\lambda}\varphi_a(\rho,\,\lambda), \ \ \  \forall
(\rho,\lambda)\in \hat \Z\times \R^*_+ ,
\end{equation}
are independent of $\lambda$ hence they define continuous functions
on $\A_\Q^{(1)}$. They are non trivial on
$\tilde\Xi_{\Q,\infty}=\hat\Z^*\times \{0\}\subset \hat\Z\times
\{0\}$ and they agree there with the functions $e_a$ of \cite{CM}.

\begin{prop}\label{Balgphipsi}
Let $\cB$ be the algebra generated by the $\varphi_a$ and $\psi_a$
defined in \eqref{phia} and \eqref{psia} above.
\begin{enumerate}
\item The expression
\begin{equation}\label{monoN}
 N(f)= \frac{1}{2\pi i}\frac{d}{d\lambda}\,f
\end{equation}
defines a derivation $N$ of $\cB$.
\item The algebra $\cB$ is stable under the derivation $Y$ that generates
the $1$-parameter semigroup $F(\mu)$ of endomorphisms of $\cS(\cG_\Q^{(1)})$ of \eqref{Fthetagauto}
and one has, at the formal level, the relation
\begin{equation}\label{mono}
F(\mu) N=\mu N F(\mu) .
\end{equation}
\item For any element $f\in \cB$ one has
\begin{equation}\label{intertwin}
\alpha \circ f(x)= f( \tilde \alpha\,x), \ \ \ \ \forall  x \in
\tilde\Xi_{\Q,\infty}\ \text{ and } \ \forall \alpha \in \Gal(\Q^{cycl}/\Q),
\end{equation}
where $\tilde \alpha \in \hat \Z^*\subset C_\Q$ is the element of
the idele class group associated to $\alpha \in
\Gal(\Q^{cycl}/\Q)$ by the class field theory isomorphism.
\end{enumerate}
\end{prop}

\proof 1) By construction $N$ is a derivation of the algebra of functions.
Moreover \eqref{psia} shows that $N(\varphi_a)= \psi_a$, while
$N(\psi_a)=0$. Thus, the derivation rule shows that $\cB$ is
stable under $N$.

2) The derivation generating the one parameter semigroup $F(\mu)$ is
given, up to sign,  by the grading operator
\begin{equation}\label{Yop}
 Y(f)= \lambda \frac{d}{d\lambda} \, f .
\end{equation}
By construction, any of the $\varphi_a$ is of degree one, \ie
$Y(\varphi_a)= \varphi_a$ and $\psi_a$ is of degree $0$. Thus, again
the derivation rule shows that $\cB$ is stable under $Y$.

3) This only involves the functions $\psi_a$,  since by construction
the restriction of $\varphi_a$ is zero on $\tilde\Xi_{\Q,\infty}$. The
result then follows from the main result of \cite{BC} in the
reformulation given in \cite{CM} (see also \cite{CMbook}, Chapter
3). In fact, all these functions take values in the cyclotomic field
$\Q^{cycl}\subset \C$ and they intertwine the action of the
discontinuous piece $\hat \Z^*$ of $C_\Q$
with the action of the Galois group of $\Q^{cycl}$.
\endproof

This is in agreement with viewing the algebra $\cB$  as the algebra
of coordinates on $\tilde\Xi_\Q$. Indeed, in the case of a global
field $\K$ of positive characteristic, the action of the Frobenius
on the points of $C(\bar \F_q)$ (which have coordinates in $\bar
\F_q$) corresponds to the Frobenius map
\begin{equation}\label{frobend}
\Fr :  u\mapsto u^q, \ \ \ \ \forall u\in  \K
\end{equation}
of the function field $\K$ of the curve $C$. The Frobenius endomorphism
$u\mapsto u^q$ of $\K$ is the operation that replaces a function
$f : C(\bar \F_q)\to \bar \F_q$ by its $q$-th power, \ie the
composition $\Fr\circ f$ with the Frobenius automorphism
$\Fr \in \Gal(\bar \F_q/\F_q)$. For $f\in \K$, one has
\begin{equation}\label{Frfq}
\Fr\circ f= f^q= f\circ \Fr,
\end{equation}
where on the right hand side $\Fr$ is the map that raises
every coordinate to the power $q$. This corresponds to the
interwtining with the Galois action discussed above.

Notice moreover that, as we have seen in Proposition \ref{globMorita},
only the semigroup $C_\Q^+$ acts on the reduced system $\cS(\cG_\Q^{(1)})$
and it acts by endomorphisms. It nevertheless acts in a bijective
manner on the points of $\Xi_\Q$. This is similar to what happens with the
Frobenius endomorphism \eqref{frobend}, which is only an endomorphism
of the field of functions $\bar \K$, while it acts bijectively (as a
Galois automorphism of the coordinates) on the
points of $C(\bar \F_q)$.

\smallskip

Further notice that there is a striking formal analogy between the
operators $F$ and $N$ of Proposition \ref{Balgphipsi} satisfying the
relation \eqref{mono} and the Frobenius and local monodromy operators
introduced in the context of the ``special fiber at arithmetic
infinity'' in Arakelov geometry (see \cite{Cons}, \cite{ConsMar}). In
particular, one should compare \eqref{mono} with \S 2.5 of
\cite{ConsMar} that discusses a notion of Weil--Deligne group at
arithmetic infinity.

\section{Functoriality of the adeles class space}

We investigate in this section the functoriality of the adele
class space $X_\K$ and of its classical subspace $\Xi_\K\subset
X_\K$, for Galois extensions of the global field $\K$.

This issue is related to the question of functoriality. Namely,
given a finite algebraic extension $\bL$ of the global field $\K$,
we want to relate the adele class spaces of both fields. Assume the
extension is a Galois extension. In general, we do not expect the
relation between the adeles class spaces to be canonical, in the
sense that it will involve a symmetry breaking choice on the Galois
group $G=\Gal(\bL/\K)$ of the extension. More precisely, the norm map
\begin{equation}\label{Normmap}
\n(a)=\prod_{\sigma\in G} \sigma(a)\in \A_\K, \ \ \ \forall  a \in \A_\bL
\end{equation}
appears to be the
obvious candidate that relates the two adeles class spaces. In fact,
since $\n(\bL)\subset \K$, the map \eqref{Normmap} passes to the
quotient and gives a natural map from $X_\bL=\A_\bL/\bL^*$ to
$X_\K=\A_\K/\K^*$ that looks like the covering required by
functoriality. However, the problem is that the norm map fails to be
surjective in general, hence it certainly does not qualify as a
covering map. In fact, this problem already occurs at the level of the
idele class group $C_\K$, namely the norm map fails to be a surjection
from $C_\bL$ to $C_\K$.

The correct object to consider is the Weil group $W_{\bL,\K}$. This
is an extension of $C_\bL$ by the Galois group $G=\Gal(\bL/\K)$, which
is not a semi-direct product. The corresponding non-trivial
$2$-cocycle is called the ``fundamental class". One has a natural morphism $t$,
called the {\em transfer}, from $W_{\bL,\K}$ to $C_\K$.
The transfer satisfies the following two properties.
\begin{itemize}
\item The morphism $t$ restricts to the norm map from $C_\bL$ to $C_\K$.
\item The morphism $t$ is surjective on $C_\K$
\end{itemize}

Thus, the correct way to understand the relation between the adeles
class spaces $X_\bL$ and $X_\K$ is by extending the construction of
the Weil group and of the the transfer map.

One obtains in this way $n$ copies of the adele class space $X_\bL$ of
$\bL$ and a map to $X_\K$ which is now a covering from $G\times
\Xi_\bL \to \Xi_\K$. This space has a natural action of the Weil
group. We explain this in more detail in what follows.

\subsection{The norm map}

We begin by recalling the well known properties of the norm map that
are relevant to our set-up. Thus, we let $\bL\supset \K$ be a finite
Galois extension of $\K$ of degree $n$, with $G=\Gal(\bL/\K)$ the
Galois group.

Since the adeles depend naturally on the field, one has a canonical
action of $G$ on $\A_\bL$. If $v\in \Sigma_\K$ is a place of $\K$,
there are $m_v$ places of $\bL$ over $v$ and
they are permuted transitively by the action of $G$. Let $G_w$
be the isotropy subgroup of $w$. Then $G_w$ is the Galois group
$G_w=\Gal(\bL_w/\K_v)$.

One has a canonical embedding of $\A_\K$ as the fixed points of the
action of $G$ on $\A_\bL$ by
\begin{equation}\label{AKGemb}
\A_\K= \A_\bL^G , \ \ \ \  (a_v)\mapsto (a_{\pi(w)}), \ \ \ \text{
with } \ \  \pi : \Sigma_\bL\to \Sigma_\K .
\end{equation}

The norm map $\n: \A_\bL\to \A_\K$ is then defined as in
\eqref{Normmap}. By \cite{Weil} IV 1, Corollary 3, it is given
explicitly by
\begin{equation}\label{explnorm}
\n(x)=z\,,\ \ z_v=\prod_{w|v}\,\n_{\bL_w/\K_v}(x_w)\,,\ \ \forall
x\in \A_\bL\,.
\end{equation}
Here the notation $w|v$ means that $w$ is a place of $\bL$ over the
place $v\in \Sigma_\K$. Also $\n_{\bL_w/\K_v}$ is the norm map of
the extension $\bL_w/\K_v$.
  When restricted to principal adeles of $\bL$ it gives
the norm map from $\bL$ to $\K$. When restricted to the subgroup
$\bL_w^*=(\ldots,1,\ldots, y,\ldots,1,\ldots)\subset \A_\bL^*$, it
gives the norm map of the extension $\bL_w/\K_v$. For nontrivial
extensions this map is never surjective, but its restriction $\n
:\cO(\bL_w)^*\to \cO(\K_v)^*$ is surjective when the extension is
unramified, which is the case for almost all places $v \in
\Sigma_\K$ (\cf \cite{Weil}, Theorem 1 p.153). In such cases, the
module of the subgroup $\n(\bL_w^*)\subset \K_v^*$ is a subgroup of
index the order of the extension $\K_v\subset \bL_w$. The
restriction of the norm map to the idele group $\A_\bL^*$ is very
far from surjective to $\A_\K^*$ and its range is a subgroup of
infinite index. The situation is much better with the idele class
groups since (\cf \cite{Weil}, Corollary  p.153) the norm map is an
open mapping $\n :C_\bL\to C_\K$ whose range is a subgroup of finite
index.

\subsection{The Weil group and the transfer map}\label{nmap}

The Weil group $W_{\bL,\K}$ associated to the Galois extension
$\K\subset \bL$ is an extension
\begin{equation}\label{CWGext}
1\to C_\bL\to W_{\bL,\K}\to G\to 1
\end{equation}
of $C_\bL$ by the Galois group $G$.
One chooses a section $s$ from $G$ and lets $a\in Z^2(G,C_\bL)$ be
the corresponding $2$-cocycle so that
\begin{equation}
a_{\alpha,\beta}= s_{\alpha\beta}^{-1}\,s_\alpha\,s_\beta, \ \ \ \
\forall \alpha,\beta\in G .
\end{equation}

The algebraic rules in $W_{\bL,\K}$ are then given by
\begin{equation}\label{rulesW}
s_\alpha\,s_\beta= s_{\alpha\beta}\,a_{\alpha,\beta} , \ \ \ \  \forall
\alpha,\beta\in G
\end{equation}
and
\begin{equation}
s_\alpha\,x\, s_\alpha^{-1}= \alpha(x), \ \ \ \ \forall \alpha\in G, \
\ \forall x\in C_\bL .
\end{equation}

The transfer homomorphism
\begin{equation}\label{transferhom}
t:  W_{\bL,\K}\to C_\K
\end{equation}
is then given by
\begin{equation}
t(x)=\n(x), \ \ \  \forall  x\in C_\bL \ \ \ \text{ and } \ \ \
t(s_\alpha)=\prod_\beta a_{\alpha,\beta}, \ \ \ \forall \alpha\in G .
\end{equation}
Its main properties are the following (see \cite{Weil-cc}).
\begin{itemize}
\item $t$ is a surjective group morphism $W_{\bL,\K}\to C_\K$.
\item Let $W_{\bL,\K}^{ab}$ be the abelian quotient of $W_{\bL,\K}$ by
the closure of its commutator subgroup $W_{\bL,\K}^c$. Then
$t$ induces an isomorphism of $W_{\bL,\K}^{ab}$ with $C_\K$.
\end{itemize}

\subsection{The covering}\label{acover}

We finally describe the resulting functoriality of the adeles class
spaces in terms of a covering map obtained by extending the Weil group
and transfer map described above.
Let, as above, $\bL\supset \K$ be a finite Galois extension of $\K$.

\begin{lem}\label{tautransf}
The transfer map extends to a map
\begin{equation}\label{tXmap}
\tau : G \times X_\bL \to X_\K
\end{equation}
of the adele class spaces.
\end{lem}

\proof
We endow $G \times X_\bL$ with a two sided
action of $G$ compatible with $\tau$. By construction the norm map
$\n$ of \eqref{Normmap} is well defined on $\A_\bL$. Since it is
multiplicative and we have $\n(\bL^*)\subset \K^*$, it induces a map
of quotient spaces $\n : X_\bL\to X_\K$.
By construction $C_\bL$ acts on $X_\bL$ and the actions by left and
right  multiplication coincide, so we use both notations.
We define the map $\tau$ as
\begin{equation}\label{taumap}
\tau : G \times X_\bL \to X_\K ,\ \ \ \
\tau(\alpha,x)= t(s_\alpha)\n(x), \ \ \ \forall  x\in X_\bL,\ \forall
\alpha \in G.
\end{equation}
This makes sense since $t(s_\alpha)\in C_\K$ and $C_\K$ acts on
$X_\K$.
By construction, the restriction of $\tau$ to $G \times C_\bL$
is the transfer map.
\endproof

One identifies $G \times C_\bL$ with $W_{\bL,\K}$ by
the map which to $(\alpha,g)\in G \times C_\bL$ associates the
element $s_\alpha\,g$ of $W_{\bL,\K}$.

In the following we use the notation
\begin{equation}\label{xg}
x^g=g^{-1}(x).
\end{equation}

We have the following result.

\begin{lem}\label{WtauLem}
Let $\bL\supset \K$ be a finite Galois extension of $\K$.
\begin{enumerate}
\item The expressions
\begin{equation}\label{leftrightact}
s_\alpha g (\beta, x)= (\alpha\beta, a_{\alpha,\beta} g^\beta x), \ \
\ \text{ and } \ \ \  (\alpha,x) s_\beta g= (\alpha\beta,
a_{\alpha,\beta} x^\beta g)
\end{equation}
define a left and a right action of $W_{\bL,\K}$ on $G \times X_\bL$.
\item The map $\tau$ of \eqref{taumap} satisfies the equivariance
property
\begin{equation}\label{tauequiv}
\tau(g x k)= t(g) \tau(x) t(k), \ \ \ \ \forall  x\in G \times X_\bL,
\ \ \text{ and } \ \ \forall g,k\in W_{\bL,\K}.
\end{equation}
\end{enumerate}
\end{lem}

\proof 1) We defined the rules \eqref{rulesW} as the natural extension
of the multiplication in $W_{\bL,\K}$ using
\begin{equation}\label{rulesW2}
s_\alpha g s_\beta h= s_\alpha s_\beta g^\beta h=
 s_{\alpha\beta} a_{\alpha,\beta} g^\beta h .
\end{equation}
Thus, the proof of associativity in the group $W_{\bL,\K}$ extends and
it implies that \eqref{leftrightact} defines a left and a right action of
$W_{\bL,\K}$ and that these two actions commute.

2) The proof that the transfer map $t$ is a group homomorphism extends
to give the required equality, since the norm map is a bimodule
morphism when extended to $X_\bL$.
\endproof

At the level of the classical points, we can then describe the
covering map in the following way.

\begin{prop}
Let $\bL$ and $\K$ be as above.
\begin{enumerate}
\item The restriction of $\tau$ to $G\times\Xi_\bL
\subset G \times X_\bL$ defines a surjection
\begin{equation}\label{Xitau}
\tau : G\times\Xi_\bL \to \Xi_\K .
\end{equation}
\item The map $\tau$ induces a surjection
\begin{equation}\label{tauXiC1}
\tau: G\times(\Xi_\bL/C_{\bL,1}) \to \Xi_\K/C_{\K,1} .
\end{equation}
\end{enumerate}
\end{prop}

\proof 1) By construction $\Xi_\bL=\cup_{w\in \Sigma_\bL} C_\bL
a^{(w)}$, where $a^{(w)}\in X_\bL$ is the class, modulo the action of
$\bL^*$, of the adele with all entries equal to $1$ except for a zero
at $w$ as in \eqref{avwadele}. Let $\pi$ denote the natural surjection
from $\Sigma_\bL$ to $\Sigma_\K$. One has
\begin{equation}\label{awapiw}
\tau(1,a^{(w)})=a^{(\pi(w))}, \ \ \ \ \forall  w\in \Sigma_\bL .
\end{equation}
In fact, one has $\tau(1,a^{(w)})= \n(a^{(w)})$. Moreover, by
\eqref{explnorm}, the adele $a=\n(a^{(w)})$ has components  $a_z=1$
for all $z\neq \pi(w)$ and $a_{\pi(w)}=0$. Thus $a=a^{(\pi(w))}$.
The equivariance of the map $\tau$ as in Lemma \ref{WtauLem}
together with the surjectivity of the transfer map from $W_{\bL,\K}$
to $C_\K$ then show that we have
$$
\tau(W_{\bL,\K} (1,a^{(w)}))= C_\K\, a^{(\pi(w))}, \ \ \ \ \forall
w\in \Sigma_\bL .
$$
For $s_\alpha g\in W_{\bL,\K}$, one has
$$ s_\alpha g (1,a^{(w)})= (\alpha,g a^{(w)}), $$
since $a_{\alpha,1}=1$. Thus, $W_{\bL,\K}(1,a^{(w)})= G\times C_\bL
a^{(w)}$ and one gets
$$
\tau( G\times C_\bL a^{(w)})= C_\K  a^{(\pi(w))}, \ \ \ \ \forall w\in
\Sigma_\bL .
$$
Since the map $\pi$ is surjective we get the conclusion.

2) The transfer map satisfies $t(C_{\bL,1})\subset C_{\K,1}$.
When restricted to the subgroup $C_\bL$ the transfer coincides
with the norm map $\n$ and in particular if $|g|=1$ one has
$|\n(g)|=1$. Thus one obtains a surjection of the quotient spaces
$$
\tau : (G\times\Xi_\bL)/C_{\bL,1} \to \Xi_\K/C_{\K,1} .
$$
Moreover, the right action of the subgroup $C_{\bL,1}\subset W_{\bL,\K}$
is given by
$$
(\alpha,x)g=(\alpha,x g).
$$
This means that we can identify
$$
(G\times\Xi_\bL)/C_{\bL,1}\sim G\times(\Xi_\bL/C_{\bL,1}).
$$
\endproof

\subsection{The function field case}\label{ff}

Let $\K=\F_q(C)$ be a global field of positive characteristic, identified
with the field of rational functions on a nonsingular curve $C$ over
$\F_q$. We consider the extensions
\begin{equation}\label{basicext}
\bL=\K\otimes_{\F_q} \F_{q^n} .
\end{equation}
The Galois group $G$ is the cyclic group of order $n$ with generator
$\sigma\in \Gal(\bL/\K)$ given by
$\sigma= {\rm id}\otimes \Fr$, where $\Fr \in \Gal(\F_{q^n}/\F_q)$ is
the Frobenius automorphism.
Given a point $x\in C(\bar\F_q)$ we let $n$ be the order of its
orbit under the Frobenius. One then has $x\in C(\F_{q^n})$ and
evaluation at $x$ gives a well defined place
$w(x)\in \Sigma_\bL$.
The projection $\pi(w(x))\in \Sigma_\K$ is a well defined place of
$\K$ which is invariant under $x\mapsto \Fr(x)$.

In the isomorphism of $\Z$-spaces
$$
\vartheta_\bL : C(\bar\F_q)\to \Xi_\bL/C_{\bL,1}
$$
described in \S \ref{curveinSect}, we have
no ambiguity for places corresponding to points $x\in C(\F_{q^n})$.
To such a point we assign simply
$$
\vartheta_\bL(x)=a^{(w(x))}\in \Xi_\bL/C_{\bL,1}.
$$

We now describe what happens with these points of $\Xi_\bL/C_{\bL,1}$
under the covering map $\tau$. We first need to see explicitly why the
surjectivity only occurs after crossing by $G$.

\begin{prop} Let $\K$ and $\bL=\K\otimes_{\F_q}\,\F_{q^n}$ be as above.
\begin{enumerate}
\item The image $\n(C_\bL)\subset C_\K$ is the kernel of the morphism
from $C_\K$ to $G=\Z/n\Z$ given by
$$ g\mapsto \rho(g)=\log_q |g| \mod n . $$
\item One has $\rho(t(s_\sigma))= 1 \mod n$, where
$\sigma\in \Gal(\bL/\K)$ is the Frobenius generator.
\end{enumerate}
\end{prop}

\proof
Since $\bL$ is an abelian extension of $\K$, one has the inclusions
\begin{equation}\label{inclusions}
\K\subset \bL \subset \K^{ab}\subset \bL^{ab},
\end{equation}
where $\K^{ab}$ is the maximal abelian extension of $\K$. Using the
class field theory isomorphisms
$$
C_\K\sim W(\K^{ab}/\K)\ \ \  \text{ and } \ \  \  C_\bL\sim
W(\bL^{ab}/\bL),
$$
one can translate the proposition in terms of Galois groups. The
result then follows using \cite{Weil-cc} p.502.
\endproof

\section{Vanishing cycles: an analogy}

We begin by considering some simple examples that illustrate some
aspects of the geometry of the adeles class space, by restricting to
the semilocal case of a finite number of places. This will also
illustrate more explicitly the idea of considering the adeles class
space as a noncommutative compactification of the idele class group.

We draw an analogy between the complement
of the idele classes in the adele classes and the
singular fiber of a degeneration. This analogy should be taken with a
big grain of salt, since this complement is a highly singular space
and it really makes sense only as a noncommutative space in the
motivic sense described in sections \ref{adclspSect}
and \ref{motivicSect} above.

\subsection{Two real places}

We first consider the example of the real quadratic field
$\K=\Q(\sqrt{2})$ and we restrict to its two real places $v_1$ and
$v_2$. Thus, we replace the adeles $\A_\K$ simply by the
product $\K_{v_1}\times \K_{v_2}$ over the real places, which is just
the product of two copies of $\R$. The ideles $\A_\K^*$
are correspondingly replaced by $\K_{v_1}^*\times \K_{v_2}^*$ and the
inclusion of ideles in adeles is simply given by the inclusion
\begin{equation}\label{incl1}
 (\R^*)^2\subset \R^2.
\end{equation}

The role of the action of the group $\K^*$
by multiplication is now replaced by the
action by multiplication
of the group $U$ of units of $\K=\Q(\sqrt{2})$. This
group is
$$
U= \Z/2\Z\times \Z
$$
where the $\Z/2\Z$ comes from $\pm 1$ and the $\Z$
is generated by the unit $u=3-2 \sqrt{2}$.
Its action on $\R^2$ is given by
the transformation
\begin{equation}\label{Sxy}
S(x,y)= (ux,u^{-1}y).
\end{equation}
Thus, in this case of two real places the semi-local version of the
adeles class space is the quotient
\begin{equation}\label{R2Uquot}
X_{v_1,v_2} := \R^2/U
\end{equation}
of $\R^2$ by the symmetry $(x,y)\mapsto (-x,-y)$ and the
transformation $S$.

Both of these transformations preserve
the function
\begin{equation}\label{fxyinv}
\tilde f: \R^2\to \R, \ \ \ \  \tilde f(x,y)= xy ,
\end{equation}
which descends to a function
\begin{equation}\label{fxyquot}
f : X_{v_1,v_2}\to \R .
\end{equation}

Moreover one has
$$
(x,y)\in (\R^*)^2\subset \R^2 \Leftrightarrow f(x,y)\neq 0
$$
and the fiber of $f$ over any non zero $\varepsilon \in \R$
is easily identified with a one dimensional torus
\begin{equation}\label{ftorus}
f^{-1}(\varepsilon)\sim \R_+^*/u^\Z, \ \ \ \forall \varepsilon\neq 0
\end{equation}
where one can use the map
$(x,y)\mapsto |x|$
to obtain the required isomorphism.

The fiber $f^{-1}(0)$ of $f$ over the point $\varepsilon =0$, on the
other hand, is no longer a one dimensional torus and it is
singular. It is the union of three pieces
\begin{equation}\label{singpieces}
f^{-1}(0)= T_1 \cup T_2\cup \{0\}
\end{equation}
corresponding respectively to
\begin{itemize}
\item $T_1$ is the locus $x=0$, $y\neq 0$, which is a torus $T_1 \sim
\R_+^*/u^\Z$ under the identification given by
the map $(x,y)\mapsto |y|$.
\item $T_2$ is the locus $x\neq 0$, $y=0$, which is also identified
with a torus $T_2 \sim \R_+^*/u^\Z$ under the analogous map
$(x,y)\mapsto |x|$.
\item The last piece is the single point $x=0$, $y=0$.
\end{itemize}
One can see that at the naive level that the quotient topology on
the singular fiber \eqref{singpieces} looks as follows. For any
point $x\in T_j$ its closure is $\bar x =\{x\}\cup \{0\}$. Moreover
the point $0$ is closed and the induced topology on its complement
is the same as the disjoint union of two one dimensional tori $T_j$.
In fact one can be more precise and see what happens by analyzing
the $C^*$-algebras involved. The $C^*$-algebra $A$ associated to the
singular fiber is by construction the crossed product
\begin{equation}\label{Acstar}
A=C_0(\tilde f^{-1}(0))\rtimes U
\end{equation}
with $\tilde f$ as in \eqref{fxyinv}. One lets
\begin{equation}\label{Aj}
A_j=C_0(V_j)\rtimes U
\end{equation}
where we use the restriction of the action of $U$ to the subsets
$$ V_j = \{(x_1,x_2)\,|\, x_j= 0 \}\sim  \R. $$
Evaluation at $0\in \R$ gives a homomorphism
$$
\epsilon_j : A_j\to C^*(U) \, .
$$

\begin{lem}\label{quottoplem}
One has an exact sequence of the form
$$
0\to C(T_j)\otimes \cK \to A_j \stackrel{\epsilon_j}{\to} C^*(U) \to
0,
$$
where $\cK$ is the algebra of compact operators.

The $C^*$-algebra $A$ is the fibered product of the $A_j$ over
$C^*(U)$ using the morphisms $\epsilon_j$.
\end{lem}

\proof The first statement follows using the fact that the action of
$U$ on $\R^*$ is free. Since the decomposition of  $\tilde
f^{-1}(0)$ as the union of the $V_j$ over their common point $0$ is
$U$-equivariant one gets the second statement.
\endproof

After collapsing the spectrum of $C^*(U)$ to a point, the topology
of the spectrum of $A_j$ is the topology of $T_j\cup \{0\}$
described above. The topology of the spectrum of $A$ is the topology
of $f^{-1}(0)$ of \eqref{singpieces} described above.

\smallskip

\subsection{A real and a non-archimedean place}

We now consider another example, namely the case of $\K=\Q$ with two
places $v_1$, $v_2$, where $v_1=p$ is a non-archimedean place
associated to a prime $p$ and $v_2=\infty$ is the real place.
Again, we replace adeles by the
product $\K_{v_1}\times \K_{v_2}$ over the two places, which in this
case is just the product
\begin{equation}\label{Ainftyp}
\K_{v_1}\times \K_{v_2} = \Q_p\times \R.
\end{equation}
The ideles are correspondingly replaced by $\K_{v_1}^*\times
\K_{v_2}^*=\Q_p^*\times \R^*$ and the inclusion is given by
\begin{equation}\label{incl2}
\Q_p^*\times  \R^*\subset \Q_p\times \R .
\end{equation}

The role of the action of the group $\K^*$
by multiplication is now replaced by the
action by multiplication
by the group $U$ of elements of $\K^*=\Q^*$
which are units outside the above two places. This
group is
\begin{equation}\label{Ugr2}
U= \Z/2\Z\times \Z ,
\end{equation}
where the $\Z/2\Z$ comes from $\pm 1$ and the cyclic group is $p^\Z$
generated by $p\in \K^*=\Q^*$.

The action of $U$ of \eqref{Ugr2} on $\R\times \Q_p$ is given by
the transformation
\begin{equation}\label{Sxy2}
S(x,y)= (px,py).
\end{equation}

By comparison with the previous case of $\K=\Q(\sqrt{2})$, notice how
in that case (\cf \eqref{Sxy}) the pair $(u,u^{-1})$ was just the image
of the element $3-2 \sqrt{2}$ under the diagonal embedding of $\K$ in
$\K_{v_1}\times \K_{v_2}$.

In the present case, the role of the adeles class space
$X_\K=\A_\K/\K^*$ is then played by its semi-local version
\begin{equation}\label{Xinftyp}
X_{p,\infty}= (\Q_p\times \R)/U
\end{equation}
quotient of $\Q_p\times \R$ by the symmetry $(x,y)\to (-x,-y)$ and
the transformation $S$. Both of these transformations preserve the
function
\begin{equation}\label{fxyinv2}
\tilde f: \Q_p\times \R \to \R_+, \ \ \ \ \tilde f(x,y)= \vert
x\vert_p\,\vert y\vert \in \R_+,
\end{equation}
which descends to a function
\begin{equation}\label{fxyquot2}
f: X_{p,\infty}\to \R_+ .
\end{equation}

Moreover, one has
$$
(x,y)\in \Q_p^*\times \R^*\subset \Q_p\times \R \Leftrightarrow f(x,y)\neq 0
$$
and the fiber of $f$ over any non zero $\varepsilon \in \R_+$
is easily identified with $\Z_p^*$
$$
f^{-1}(\varepsilon)\sim \Z_p^*,  \ \ \  \forall \varepsilon\neq 0 .
$$
In fact, one can use the fundamental domain
$$
\Z_p^*\times \R_+^*
$$
for the action of $U$ on  $\Q_p^*\times \R^*$ to obtain the required
isomorphism.

The fiber $f^{-1}(0)$ of $f$ over the point $\varepsilon =0$ is no longer
$\Z_p^*$ and once again it is singular. It is again described as the
union of three pieces
\begin{equation}\label{singpieces2}
f^{-1}(0)= T_p \cup  T_\infty \cup \{0\},
\end{equation}
which have, respectively, the following description.
\begin{itemize}
\item $T_p$ is the locus $x=0$, $y\neq 0$, which is identified with a
torus $T_p \sim \R_+^*/p^\Z$, using the map $(x,y)\mapsto |y|$.
\item $T_\infty$ is the locus $x\neq 0$, $y= 0$, which gives the
compact space $T_\infty \sim \Z_p^*/\pm 1$ obtained
as quotient of $\Q_p^*$ by the action of $U$.
\item The remaining piece is the point $x=0$, $y=0$.
\end{itemize}
The description of the topology of $f^{-1}(0)$ is
similar to what happens in the case of $\Q(\sqrt{2})$
analyzed above.

What is not obvious in this case is how the totally disconnected fiber
$f^{-1}(\varepsilon)\sim \Z_p^*$ can tie in with the torus
$T_p \sim \R_+^*/p^\Z$ when $\varepsilon\to 0$.

To see what happens, we use the map
\begin{equation}\label{gxymap}
X_{p,\infty} \ni (x,y) \mapsto g(x,y)=\text{class\,of}\  |y|\in
\R_+^*/p^\Z.
\end{equation}
This is well defined on the open set $y\neq 0$. It is continuous and
passes to the quotient. Thus, when a sequence $(x_n,y_n)\in
X_{p,\infty}$ converges to a point $(0,y)\in T_p$, $y\neq 0$, one
has $g(0,y)=\lim_n g(x_n,y_n)$.

The point then is simply that we have the relation
\begin{equation}\label{gxyfxy}
g(x,y)= f(x,y)\in \R_+^*/p^\Z .
\end{equation}
In other words, $g(x_n,y_n)=\varepsilon_n$ with $(x_n,y_n)$
in the fiber $f^{-1}(\varepsilon_n)$
and the point of the singular fiber $T_p$ towards which
$(x_n,y_n)\in X_{p,\infty}$ converges depends only on the value of
$\varepsilon_n$ in $\R_+^*/p^\Z$.

\begin{figure}
\begin{center}
\includegraphics[scale=0.6]{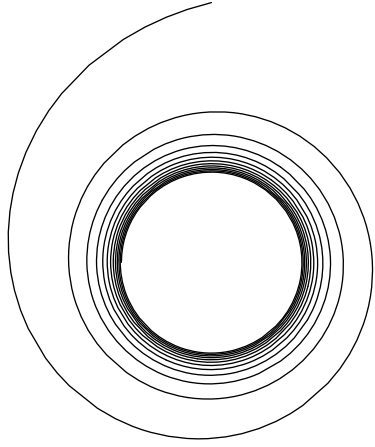}
\caption{The  limit cycle of a foliation. \label{hol}}
\end{center}
\end{figure}

This phenomenon is reminiscent of the behavior of holonomy in
the context of foliations, using a logarithmic
scale $\R_+^*/p^\Z \sim \R/(\Z \log p)$.
It corresponds to what happens
in the limit cycle of the foliation associated to
a flow as depicted in Figure \ref{hol}.

 As we argued in \cite{CCM} (see also \S \ref{Frob0sect}
and \S \ref{ScaleCorrSect} here above), the role of Frobenius in
characteristic zero is played by the one parameter group $\Fr(t)$
with $t\in \R$ which corresponds to the action of $\R$ on the adele
class space $X_\Q=\A_\Q/\Q^*$ given in the above logarithmic scale,
namely
\begin{equation}\label{Frtlog}
\Fr(t)(a)=e^t\,a, \ \ \ \  \forall a \in X_\Q .
\end{equation}
Its orbit over $p\in \Sigma_\Q$ is of length $\log p$ and it
corresponds, in the simplified picture of $X_{p,\infty}$, to the
component $T_p$ of the singular fiber $f^{-1}(0)$.

\subsection{Singularities of maps}

The simple examples described above illustrate how
one can use the function $f(x)=|x|$ in general, and see the place
where it vanishes as the complement of $C_\K$ in the adeles class
space $X_\K$. This provides a way of thinking of the inclusion
of $C_\K$ in $X_\K$ in terms of the notions of ``singular fiber" and
``generic fiber" as seen in the examples above. The generic fiber
appears to be typically identified with $C_{\K,1}$, with the union of
the generic fibers giving $C_\K$ as it should.
This suggests the possibility of adapting to our noncommutative
geometry context some aspects of the well developed theory of nearby
and vanishing cycles. A brief dictionary summarizing this analogy is
given here below.

\medskip

\begin{center}
\begin{tabular}{|c|c|}
\hline &  \\
 Total space &  Adele class space $X_\K=\A_\K/\K^*$ \\
&\\ \hline &  \\
The map $f$ &  $f(x)=|x|$ \\
&\\
\hline &  \\
 Singular fiber & $X_\K\smallsetminus C_\K= f^{-1}(0)$ \\
&\\ \hline &  \\
Union of generic fibers &  $C_\K=\,f^{-1}(\{0\}^c)$ \\
&\\ \hline
\end{tabular}
\end{center}


\medskip


\begin{thebibliography}{99}

\bibitem{Andre} Y.~Andr\'e, {\em Une introduction aux motifs},
Panoramas et Synth\`eses, Vol.17, Soci\'et\'e math\'ematique de
France, 2005.

\bibitem{AB} M.F.~Atiyah, R.~Bott, {\em A Lefschetz fixed point
formula for elliptic complexes: I}, Annals of Math. Vol. 86 (1967)
374--407.

\bibitem{BEK} S.~Bloch, H.~Esnault, D.~Kreimer, {\em On motives
associated to graph polynomials}, Commun. Math. Phys. Vol. 267
(2006) 181--225.

\bibitem{EB} E.~Bombieri, {\em Problems of the Millenium:
 The Riemann Hypothesis}, Clay mathematical Institute (2000).


\bibitem{BC} J.B.~Bost, A.~Connes, {\em Hecke algebras, type III
factors and phase transitions with spontaneous symmetry breaking in
number theory}.  Selecta Math. (N.S.)  1  (1995),  no. 3, 411--457.

\bibitem{Co-th} A.~Connes, {\em Une classification des facteurs de
type ${\rm III}$}.  Ann. Sci. \'Ecole Norm. Sup. (4)  6  (1973),
133--252.

\bibitem{CoExt} A.~Connes, {\em  Cohomologie cyclique et foncteurs
${\rm Ext}\sp n$}.  C. R. Acad. Sci. Paris S\'er. I Math. 296
(1983), no. 23, 953--958.

\bibitem{CoIHES} A.~Connes, {\em Noncommutative differential
geometry}.  Inst. Hautes \'Etudes Sci. Publ. Math.  No. 62 (1985),
257--360.

\bibitem{Co94} A.~Connes, {\em Noncommutative geometry}. Academic
Press, 1994.

\bibitem{Co-zeta} A.Connes {\em Trace formula in noncommutative
geometry and the zeros of the Riemann zeta function}.
Selecta Math. (N.S.)  5  (1999),  no. 1, 29--106.

\bibitem{CCM} A.~Connes, C.~Consani, M.~Marcolli, {\em
Noncommutative geometry and motives: the thermodynamics of
endomotives}, preprint math.QA/0512138.

\bibitem{CM} A.~Connes, M.~Marcolli, {\em From physics to number
theory via noncommutative geometry. Part I. Quantum statitical
mechanics of $\Q$-lattices}, preprint math.NT/0404128.

\bibitem{CMbook} A.~Connes, M.~Marcolli, {\em Noncommutative geometry
from quantum fields to motives} (tentative title), book in preparation.

\bibitem{CMR} A.~Connes, M.~Marcolli, N.~Ramachandran, {\em KMS states
and complex multiplication}, Selecta Math. (New Ser.)  Vol.11  (2005),
no. 3-4, 325--347.

\bibitem{CMR2} A.~Connes, M.~Marcolli, N.~Ramachandran, {\em KMS states
and complex multiplication, II}, in ``Operator Algebras - The Abel
Symposium 2004", pp.15--60, Springer Verlag, 2006.

\bibitem{CoSka} A.~Connes, G.~Skandalis, {\em The longitudinal index
theorem for foliations}, Publ. RIMS Kyoto Univ. 20 (1984)
1139--1183.

\bibitem{ct} A.~Connes, M.~Takesaki, {\it The flow of weights on
factors of type III}. Tohoku Math. J., 29, (1977) 473--575.

\bibitem{Cons} C.~Consani, {\em Double complexes and Euler
$L$-factors}.  Compositio Math.  111  (1998),  no. 3, 323--358.

\bibitem{ConsMar} C.~Consani, M.~Marcolli, {\em Archimedean cohomology
revisited}, in ``Noncommutative Geometry and Number Theory"
pp.109--140. Vieweg Verlag, 2006.

\bibitem{ConsMar-ff} C.~Consani, M.~Marcolli, {\em Quantum statistical
mechanics over function fields}, math.QA/0607363, to appear in Journal
of Number Theory.

\bibitem{De} P.~Deligne, {\em Th\'eorie de Hodge III},  Publ. Math. IHES Vol. 44 (1974)
5--78.

\bibitem{GM} A.~Goncharov, Yu.~Manin {\em Multiple zeta motives and
moduli spaces $\bar M_{0,n}$}, Compos. Math. Vol. 140 no. 1 (2004)
1--14.

\bibitem{GS} V.~Guillemin, S.~Sternberg, {\em Geometric asymptotics},
Math. Surveys Vol. 14, American Mathematical Society, 1977.

\bibitem{gui}  V. Guillemin,  {\em Lectures on spectral theory of
elliptic operators},   Duke Math. J.,  Vol. 44, 3 (1977), 485-517.

\bibitem{HaPau} E.~Ha, F.~Paugam, {\em Bost-Connes-Marcolli systems
for Shimura varieties. I. Definitions and formal analytic
properties}, IMRP Int. Math. Res. Pap.  2005,  no. 5, 237--286.

\bibitem{Jacob} B.~Jacob, {\em Bost-Connes type systems for function
fields}, math.OA/0602554, to appear in Journal of Noncommutative
Geometry.

\bibitem{Man-mot}  Yu.I.~Manin, {\em Correspondences, motifs and
monoidal transformations}, Mat. Sb. (N.S.) 77 (119) 1968, 475--507.

\bibitem{Man-zetas} Yu.I.~Manin, {\em Lectures on zeta functions
and motives (according to Deninger and Kurokawa)}. Columbia
University Number Theory Seminar (New York, 1992). Ast\'erisque  No.
228  (1995), 4, 121--163.

\bibitem{Man1} Yu.I.~Manin, {\em Real Multiplication and
noncommutative geometry (ein Alterstraum)}, in  ``The legacy of
Niels Henrik Abel'',  pp.685--727, Springer, Berlin, 2004.

\bibitem{Man2} Yu.I.~Manin, {\em Von Zahlen und Figuren}, preprint
arXiv math.AG/0201005.

\bibitem{Mar} M.~Marcolli, {\em Arithmetic noncommutative
geometry}, University Lectures Series, Vol.36, American Mathematical
Society, 2005.

\bibitem{Meyer} R.~Meyer, {\em On a representation of the idele class group
related to primes and zeros of $L$-functions}.  Duke Math. J.  127
(2005),  no. 3, 519--595.

\bibitem{Riemann} B.~Riemann, {\it Mathematical Works}, Dover,
New York, 1953.

\bibitem{SD} H.P.F.~Swinnerton-Dyer, {\em Applications of
Algebraic Geometry to Number Theory}, Proceedings Symposia in Pure
Math. Vol. XX (1969) 21--26.

\bibitem{tt}  M.~Takesaki, {\it Tomita's theory of modular
Hilbert algebras and its applications}. Lecture Notes in Math., 28,
Springer, 1970.

\bibitem{Tak} M.~Takesaki, {\em  Duality for crossed products and the structure
of von Neumann algebras of type III},  Acta Math. (131) (1973),
249-310.

\bibitem{Weil-lett} A.~Weil, {\em Letter to Artin},
Collected Papers, Vol.I (1980) 280--298.

\bibitem{Weil-RH} A.~Weil, {\em On the Riemann hypothesis in
function-fields},  Proc. Nat. Acad. Sci. U.S.A.  Vol.27 (1941)
345--347.

\bibitem{Weil-cc} A.~Weil, {\em Sur la th\'eorie du corps de classes},
J. Math. Soc. Japan, Vol.3 (1951) 1--35.

\bibitem{Weil-distr} A.~Weil, {\em
Fonction zeta et distributions},
S\'eminaire Bourbaki, Vol. 9 (1966), Exp. No. 312, 523--531.

\bibitem{weilpos} A.~Weil, {\it Sur les formules explicites de la
th\'eorie des nombres premiers}, Oeuvres complètes, Vol. 2, 48--62.

\bibitem{weilexplicit} A.~Weil, {\it Sur les formules explicites de la
th\'eorie des nombres}, Izv. Mat. Nauk., (Ser. Mat.) Vol.36 (1972)
3--18.

\bibitem{Weil} A.~Weil, {\em Basic Number Theory}, Reprint of the
second (1973) edition. Classics in Mathematics. Springer-Verlag,
1995.

\end{thebibliography}
\end{document}